\documentclass[a4paper,11pt]{article}

\usepackage{amssymb,amsmath,amsthm,latexsym}
\usepackage[colorlinks=true,linkcolor=blue,citecolor=blue]{hyperref}
\usepackage{xcolor}

\theoremstyle{plain}
\newtheorem{thm}{Theorem}[section]
\newtheorem{lem}[thm]{Lemma}
\newtheorem{pro}[thm]{Proposition}
\newtheorem{cor}[thm]{Corollary}

\theoremstyle{definition}
\newtheorem{defn}[thm]{Definition}
\newtheorem{eg}[thm]{Example}
\newtheorem{assum}[thm]{Assumption}

\theoremstyle{remark}
\newtheorem{rem}[thm]{Remark}

\numberwithin{equation}{section}

\newcommand{\R}{\mathbb{R}}

\newcommand{\E}{\mathbb{E}}
\newcommand{\Prob}{\mathbb{P}}

\newcommand{\F}{\mathbb{F}}
\newcommand{\Sbb}{\mathbb{S}}
\newcommand{\Hbb}{\mathbb{H}}
\newcommand{\Lbb}{\mathbb{L}}
\newcommand{\cF}{\mathcal{F}}
\newcommand{\cM}{\mathcal{M}}
\newcommand{\cT}{\mathcal{T}}
\newcommand{\cP}{\mathcal{P}}
\newcommand{\cO}{\mathcal{O}}
\newcommand{\cL}{\mathcal{L}}
\newcommand{\cE}{\mathcal{E}}
\newcommand{\cK}{\mathcal{K}}

\newcommand{\cB}{\mathcal{B}}
\newcommand{\cQ}{\mathcal{Q}}
\newcommand{\cS}{\mathcal{S}}
\newcommand{\cN}{\mathcal{N}}
\newcommand{\Dp}{\Delta^{+}}
\newcommand{\Dm}{\Delta}
\newcommand{\cadlag}{c\`adl\`ag}
\newcommand{\ladlag}{l\`adl\`ag}
\newcommand{\trans}{^{\scriptscriptstyle\top}}
\newcommand{\ind}{\mathbf{1}}
\DeclareMathOperator*{\esssup}{ess\,sup}
\DeclareMathOperator*{\essinf}{ess\,inf}

\begin{document}

\title{Minimal Reflected BSDEs and Superhedging\\
of Irregular American Claims}

\author{%
Siham BOUHADOU\thanks{Cadi Ayyad University, Faculty of Sciences Semlalia,
Department of Mathematics, B.P.~2390, Marrakesh 40.000, Morocco.
E-mail: sihambouhadou@gmail.com}\and
Imade FAKHOURI\thanks{Ibnou Zohr University, National School of Applied Sciences,
LISAD Laboratory, \'Equipe Aide \`a la D\'ecision, B.P.~1136, Agadir, Morocco.
E-mail: i.fakhouri@uiz.ac.ma}\,\,\,\thanks{Corresponding author.}\and
Youssef OUKNINE\thanks{Cadi Ayyad University, Marrakesh, and Mohammed VI
Polytechnic University, Africa Business School, Rabat, Morocco.
E-mails: ouknine@uca.ma, youssef.ouknine@um6p.ma. Supported by the Hassan~II
Academy of Sciences and Technology, Rabat, Morocco.}}

\date{\today}
\maketitle

\begin{abstract}
We study the superhedging of American claims with an irregular optional payoff
process $\xi$ over the full set $\cM(X)$ of equivalent local martingale
measures of a locally bounded market.  Combining the Mertens decomposition with
the unrestricted optional decomposition of F\"ollmer and Kabanov, we show that a
process which is a strong optional supermartingale of class~(D) under every
$Q\in\cM(X)$ --- with no bound on it from either side ---
decomposes as $S=S_0+\int H\,dX-A-C_-$, with a common integrand $H$, a
right-continuous increasing process $A$ and a purely discontinuous increasing
process $C$ given explicitly by the right jumps of $S$.  Applied to the robust
Snell envelope, this yields the superhedging duality
\[
p_0(\xi)=\sup_{Q\in\cM(X)}\ \sup_{\tau\in\cT}\ \E^{Q}\bigl[\xi_\tau\bigr],
\]
with the attaining gross wealth equal to $V+A+C_-$,
where $V$ is the robust Snell envelope; no regularity of $\xi$ in time is
needed for this.  The payoff need not be bounded: it suffices that
$\{\xi_\tau\}$ be uniformly integrable jointly over measures and stopping times,
for which a uniform $L^{p}$ bound with $p>1$ is sufficient, or that $|\xi|$ be
dominated by a traded uniformly integrable martingale.  For payoffs bounded
below the infimum is attained in the usual class of strategies with wealth
bounded below by a constant.  For payoffs unbounded below this can fail,
already in a complete market; the infimum is then attained once the constant
lower bound is replaced by $-M$, where $M$ is the traded martingale dominating
$|\xi|$.  In
an incomplete Brownian market we then study the
associated reflected BSDE with a bounded optional obstacle, with no path
regularity assumed, dominated at every stopping time, and a Lipschitz driver
depending on the value.  We prove
existence and uniqueness of the minimal supersolution, a representation at
stopping times, comparison, and an equivalent dynamic minimum condition, and
we extend them to unbounded data dominated by a traded martingale --- in
particular, for drivers with $f(\cdot,0)$ bounded, to the integrability
conditions of the linear theory --- with minimality among all supersolutions of
class~(D), by measuring the solution in a norm weighted by that martingale.  For
a European-type obstacle, the classical contact conditions on the total
regulator hold if and only if the terminal claim is attainable by a bounded
wealth process; hence in an incomplete market there are data for which no
supersolution satisfies them, while the right-jump regulator $C$ satisfies its
contact condition for every obstacle.  We relate these statements to
constrained reflected BSDEs and to second-order reflected BSDEs, where the
failure of the Skorokhod condition for a total increasing process is already
recognised.
\end{abstract}



\section{Introduction}\label{sec:intro}

\subsection{Motivation}

Two classical decomposition theorems sit at the foundation of optimal stopping
and of hedging in incomplete markets.

The first is the \emph{Mertens decomposition} \cite{Mertens1972,
Dellacherie-Meyer-1980}: a strong optional supermartingale of class~(D) --- a
supermartingale indexed by stopping times, with no right-continuity assumed ---
splits as
\begin{equation}\label{eq:intro-mertens}
S_t=M_t-A_t-C_{t-},
\end{equation}
with $M$ a \cadlag{} uniformly integrable martingale, $A$ predictable
non-decreasing and $C$ adapted, purely discontinuous and non-decreasing.  The
third term is what distinguishes the optional theory from the classical one: it
records the \emph{right} jumps of $S$, which appear as soon as the underlying
reward is not right-continuous.  This is the decomposition that governs optimal
stopping with irregular payoffs, and it is the backbone of the theory of
reflected BSDEs with non-right-continuous obstacles developed by Grigorova,
Imkeller, Offen, Ouknine and Quenez \cite{Grigorovaetal2017}; see also
\cite{Grigorova2020} for optimal stopping with $f$-expectations in the
irregular case and \cite{QuenezSulem2014} for the connection with dynamic risk
measures.

The second is Kramkov's \emph{optional decomposition} \cite{Kramkov96} (see also
\cite{ElKarouiQuenez1995,FollmerKabanov1998,FollmerKramkov1997,Jacka2012}): if a
non-negative \cadlag{} process $Y$ is a supermartingale under \emph{every}
equivalent local martingale measure $Q\in\cM(X)$ of a locally bounded process
$X$, then
\begin{equation}\label{eq:intro-kramkov}
Y_t=Y_0+\int_0^tH_s\,dX_s-A_t
\end{equation}
with $H$ predictable and $X$-integrable and $A$ adapted non-decreasing.  The
finite-variation part cannot be taken predictable in general --- this is exactly
why the decomposition is called \emph{optional} --- and $A$ has a concrete
financial meaning: it is the cumulative consumption of a superhedging strategy.
Equation~\eqref{eq:intro-kramkov} is the analytic content of the superhedging
duality in incomplete markets.

The two theorems have complementary weaknesses.  Mertens' decomposition handles
irregular processes but is tied to a single measure: both $M$ and $A$ in
\eqref{eq:intro-mertens} depend on $Q$.  Kramkov's decomposition is uniform in
$Q$ but concerns right-continuous processes, and therefore does not see the
right jumps that irregular American-type payoffs create.  A common refinement,
for a non-linear expectation and a continuous integrator, is due to Bouchard,
Possama\"i and Tan \cite{BouchardPossamaiTan2016}; we give the linear version
for a general locally bounded integrator under class~(D), in a form that
separates the right-jump part.  Section~\ref{subsec:literature} compares the
two precisely.

The refinement is a tool, and the destination is a financial theorem: the
superhedging price of an American claim whose payoff process is merely
optional, with no time regularity and no two-sided bound, the infimum being
attained for payoffs bounded below, together with a strategy and
the increasing processes that go with it (Theorem~\ref{thm:robust-american}).
In the Brownian part of the paper the same objects are studied as a reflected
BSDE with an optional obstacle, dominated at every stopping time, and with the
constraint that incompleteness imposes on the integrand.  For that equation we
prove existence, uniqueness and comparison for the minimal supersolution, and,
for European-type obstacles, we determine exactly when the classical Skorokhod
conditions can be imposed on the total regulator.

\subsection{Main results}

The paper proves four things, and they do not stand in the same relation to
the literature; Section~\ref{subsec:literature} gives the detail.

\paragraph{A uniform optional decomposition (Section~\ref{sec:decomposition}).}
Theorem~\ref{thm:uniform} states that if an optional process $S$
is a strong optional supermartingale of class~(D) under every $Q\in\cM(X)$ ---
no bound on $S$ being assumed --- then
\begin{equation}\label{eq:intro-main}
S_t=S_0+\int_0^tH_s\,dX_s-A_t-C_{t-},\qquad 0\le t\le T,
\end{equation}
where $H$, $A$ and $C$ do \emph{not} depend on $Q$, $\int H\,dX$ is a uniformly
integrable martingale under every $Q$, and the purely discontinuous part is
given explicitly and pathwise by
\[
C_t=-\sum_{s\le t,\ s<T}\Dp S_s .
\]
The mechanism is the observation (Lemma~\ref{lem:aggreg-C}) that in the
Mertens decomposition under a fixed $Q$ the right-jump part is determined by
the pathwise relation $\Dp S=-\Dm C^{Q}$, and is therefore the same for all
$Q$; once it is removed, $S+C_-$ is right-continuous and an optional
decomposition theorem applies.  The route --- regularise on the right, then
apply an optional decomposition --- is that of \cite{BouchardPossamaiTan2016}.
The theorem applied is the unrestricted one of F\"ollmer and Kabanov
\cite[Theorem~1]{FollmerKabanov1998}, which needs no bound on the process being
decomposed; the uniformly integrable martingale property of $\int H\,dX$, which
a lower bound would have delivered through the Ansel--Stricker lemma, is
recovered instead from the class~(D) hypothesis by way of
\cite[Lemma~A1]{BiaginiSirbu2012} (Remark~\ref{rem:no-lower-bound}).  What is
added here is an arbitrary locally bounded \cadlag{} integrator, the class~(D)
regime and the explicit separation of the right-jump part.
Corollary~\ref{cor:uniform-submartingale} states the sign counterpart, which is
what the companion papers \cite{BA,P4} use; the companion paper on game options
\cite{GO} uses both forms.  Section~\ref{sec:decomposition}
also records when $A$ can be chosen predictable: a sufficient condition is
Jacka's scalar closure of the localised density logarithms \cite{Jacka2012},
and Example~\ref{eg:counterexample-H1} shows that closure under scalars in
$[0,1]$, which convexity supplies, does not suffice.  Removing the lower bound
makes Jacka's theorem inapplicable, since it takes a non-negative process as
input, and Theorem~\ref{thm:predictable} is proved instead by compensating the
increasing part under $\Prob$ and transferring it to every $Q$ through the
continuity of the density processes.

\paragraph{The robust Snell envelope (Section~\ref{sec:snell}).}
For an optional reward $\eta$ satisfying either of the two integrability
conditions of Section~\ref{sec:snell} --- uniform integrability of
$\{\eta_\tau\}$ over measures and stopping times (Assumption~\ref{assum:E}), or
domination of $|\eta|$ by a traded uniformly integrable martingale
(Assumption~\ref{assum:G}) --- the family
\[
V_\theta:=\esssup_{Q\in\cM(X)}\ \esssup_{\tau\in\cT_\theta}\
\E^{Q}\!\left[\eta_\tau\mid\cF_\theta\right],\qquad\theta\in\cT ,
\]
is a strong $Q$-supermartingale system for every $Q\in\cM(X)$
(Theorem~\ref{thm:snell-supermart}), aggregates into an optional \ladlag{}
process (Theorem~\ref{thm:snell-aggreg}), and is the smallest process which is
a strong $Q$-supermartingale for every $Q$ and dominates $\eta$
(Theorem~\ref{thm:snell-minimal}).  None of these statements is new in
substance: the characterisation is proved for a general dominated family of
priors by Arharas, Bouhadou, Hilbert and Ouknine \cite{Arharas2025}, who also
aggregate the family, and the aggregation of supermartingale systems is
classical \cite{Dellacherie-Lenglart-1982,Kobylanski2012}.  We include the
short proofs because the order of the argument (supermartingale property first,
aggregation second) is what allows a reward with no semicontinuity at all, and
because the unbounded regime needs the estimates of
Proposition~\ref{pro:envelope-UI}.  The two integrability conditions are
compared in Corollary~\ref{cor:ui-traded-domination} and
Example~\ref{eg:G-not-UI}: the first implies the second, and not conversely.

\paragraph{Robust superhedging (Section~\ref{sec:financial}).}
Theorem~\ref{thm:robust-american} identifies the superhedging price of the
American claim with payoff $\xi$ with $V_0(\xi)=\sup_Q\sup_\tau\E^{Q}[\xi_\tau]$
in the general locally bounded market, for any optional $\xi$ which is bounded
below and satisfies one of the two integrability conditions.  The
infimum is attained, by the strategy $H$ of the decomposition of $V(\xi)$, and
the attaining gross wealth is $V+A+C_-$.  For payoffs unbounded below the
same value formula holds, and the infimum is attained once the constant lower
bound on wealth is replaced by $-M$, where $M$ is the traded martingale
dominating $|\xi|$ (Theorem~\ref{thm:robust-american-signed});
Example~\ref{eg:no-attainment} shows in a \emph{complete} market that
attainment can fail in the constant-lower-bound class --- so the admissibility convention, and not
only the payoff hypothesis, has to change.  For a right-continuous payoff the
duality is Kramkov's \cite[Theorem~3.3]{Kramkov96}.  For an irregular payoff
the value process is in general not the capital of any portfolio, and
Remark~\ref{rem:kramkov-american} states precisely where the boundary lies.
With a Lipschitz driver describing non-linear wealth dynamics, the superhedging
price is $Y_0$, the initial value of the minimal solution below
(Theorem~\ref{thm:nonlinear-wealth}).

\paragraph{Reflected BSDEs and the Skorokhod conditions (Section~\ref{sec:bsde}).}
In an incomplete It\^o market --- $n$ Brownian motions and $d$ traded assets,
with $\operatorname{rank}\sigma<n$ on a non-null set --- we study
\begin{equation}\label{eq:intro-rbsde}
Y_\tau=\xi_T+\int_\tau^T f(s,Y_s)\,ds-\int_\tau^T H_s\,dX_s
+A_T-A_\tau+C_{T-}-C_{\tau-},\qquad Y_\tau\ge\xi_\tau,\quad\tau\in\cT,
\end{equation}
for a bounded optional obstacle, with no left or right limits assumed, and a
driver Lipschitz in $y$.  We define a
solution as the \emph{minimal supersolution} (Definition~\ref{def:rbsde}),
following the constrained-BSDE literature
\cite{CvitanicKaratzasSoner1998,Peng1999,PengXu2010}, and prove existence,
uniqueness, the representation \eqref{eq:representation} at stopping times and
comparison (Theorems~\ref{thm:existence} and~\ref{thm:comparison}).
Corollary~\ref{cor:dyn-minimality} restates minimality as an equivalent
condition at each stopping time, of the shape of the minimum condition of
second-order BSDEs.  Theorem~\ref{thm:weighted} removes the boundedness: for
data dominated by a traded martingale $M$, the same map is a contraction in a
norm weighted by $M$, and the solution is unique and minimal among all
supersolutions of class~(D) under every $Q$, each of which is itself dominated
by a traded martingale (Lemma~\ref{lem:classD-tr}).  Comparison and the dynamic
characterisation carry over (Corollary~\ref{cor:weighted-comparison}).  For
drivers with $f(\cdot,0)$ bounded this covers exactly the integrability
conditions of the linear theory (Corollary~\ref{cor:weighted-EG}).  No maximal
inequality is used, and Example~\ref{eg:weighted-not-L2} shows that the
square-integrable estimates of the bounded case have no counterpart.

Proposition~\ref{pro:no-skorokhod} determines when the classical Skorokhod
conditions can be imposed on the total regulator of \eqref{eq:intro-rbsde}.
For a European-type obstacle they hold for some supersolution if and only if
the terminal claim is attainable by a bounded wealth process.  Consequently
(Corollary~\ref{cor:dichotomy}) for obstacles with finite left limits they hold
for every admissible datum when $\cM(X)$ is a singleton, and fail for some data
as soon as it is not.  On the
other hand the right-jump regulator $C$ satisfies its contact condition for
every obstacle (Proposition~\ref{pro:C-skorokhod}), so the obstruction concerns
the predictable part $A$ only.  The mechanism --- a single increasing process
that carries both reflection and the premium for non-attainability cannot be
flat off the contact set --- is not new: Matoussi, Possama\"i and Zhou
\cite[Remark~3.2, p.~2429]{MatoussiPossamaiZhou2013} observe that the total increasing
process of a second-order reflected BSDE fails the Skorokhod condition for this
reason, and constrained reflected BSDEs \cite{PengXu2010} keep a Skorokhod
condition by carrying the reflection in a separate process.  What
Proposition~\ref{pro:no-skorokhod} adds is an exact criterion in the dominated
setting, and Problem~P3 asks for the corresponding split of $A$.

\subsection{Relation to the literature}\label{subsec:literature}

\paragraph{Decomposition.}
Theorem~\ref{thm:uniform} reduces to Kramkov's theorem \cite{Kramkov96} for
right-continuous processes satisfying the class~(D) requirement of
\eqref{eq:H-super} (then $C\equiv0$), and to Mertens' theorem
\cite{Mertens1972} when $\cM(X)=\{\Prob\}$, so that $\int H\,dX$ is an
arbitrary martingale.  It does not contain Kramkov's theorem: class~(D) is not
needed there, and Example~\ref{eg:bessel} is a process to which his theorem
applies and ours does not.

The closest antecedent is Bouchard, Possama\"i and Tan
\cite{BouchardPossamaiTan2016}.  Their Theorem~3.1 is a Mertens-type
decomposition, with aggregation, for $\cE^{g}$-supermartingale \emph{systems}
which assumes no right-continuity of the system, and their Theorem~4.1 deduces
from it an optional decomposition, uniform over the measures equivalent to
$\Prob$ under which the integrator is a martingale.  For a null driver their
Theorem~4.1 and Theorem~\ref{thm:uniform} concern the same object and are
reached by the same route.  Three differences remain.
\begin{itemize}\itemsep1pt
\item \emph{The integrator.}  Their filtration, like ours, is general; it is
the integrator that is restricted.  Theirs is a continuous martingale
$\int\alpha^{\top}dW$ with $\alpha$ bounded with bounded inverse.  Ours is an
arbitrary locally bounded \cadlag{} semimartingale, which may jump and whose
quadratic variation need not be absolutely continuous.
\item \emph{The integrability.}  They assume
$\esssup_{\tau\in\cT}|S(\tau)|\in L^{p}(Q)$ with $p>1$, a condition on the
running maximum above the $L^{1}$ level; we assume class~(D) under each $Q$,
which is a condition at stopping times at the $L^{1}$ level, and no bound on
$S$.  Neither hypothesis implies the other.
\item \emph{The form of the conclusion.}  Their Theorem~4.1 asserts that
$S-\int Z\,dM$ is non-increasing.  We split that process into $A$ and the
purely discontinuous $C_-$, with $C$ given by the right jumps of $S$.
\end{itemize}
Their framework also covers questions outside the scope of this paper: their
decomposition is stated for a non-linear expectation with a Lipschitz driver,
and their Theorem~4.2 gives a dual formulation for minimal supersolutions of
BSDEs with a constraint on the gains process.

\paragraph{Robust optimal stopping.}
Arharas, Bouhadou, Hilbert and Ouknine \cite{Arharas2025} study the value
family
$R(v)=\esssup_{\Prob'\in\cP}\esssup_{\tau\ge v}\E^{\Prob'}[Y(\tau)\mid\cF_v]$
for a reward family indexed by stopping times and a general dominated family
of priors $\cP$.  They characterise $R$ as the smallest $\cP$-supermartingale
family dominating the reward, and in their Theorem~5.1 they aggregate a
$\cP$-supermartingale family into an optional process, through
\cite[Theorem~3.1]{BouchardPossamaiTan2016}, and decompose it with a single
optional non-decreasing process and an abstract $\cP$-martingale, under a
subspace hypothesis on the density-generating martingales and square
integrability.  Counterparts of Theorems~\ref{thm:snell-supermart},
\ref{thm:snell-aggreg} and~\ref{thm:snell-minimal} for a general family $\cP$
are proved there under these hypotheses; for a single measure the characterisation and the aggregation go back to Kobyla\'nski
and Quenez \cite[\S1 and Appendix~B, Proposition~B.1]{Kobylanski2012}.  The
present paper differs from \cite{Arharas2025} in taking $\cP=\cM(X)$, which is
what makes the martingale part a stochastic integral of $X$ and yields a
hedging strategy, and in splitting the non-decreasing part into $A$ and $C$.
The price of the first point is that $\cM(X)$ must be the \emph{entire} set of
equivalent local martingale measures of $X$ (Remark~\ref{rem:full-family}).
Concretely, \cite[Theorem~5.1]{Arharas2025} does not directly identify its
martingale component as a stochastic integral of the traded process $X$, nor
does it separate the right-jump regulator.  For the full family $\cM(X)$ an
abstract uniform local martingale bounded below can be represented as an
$X$-integral, by Lemma~\ref{lem:unif-loc-mart}; the argument here uses the
pathwise identification of the right-jump part together with the unrestricted
optional decomposition theorem of \cite{FollmerKabanov1998}
to obtain the trading integrand and the two increasing processes under the
class~(D) assumption alone.  Lemma~\ref{lem:unif-loc-mart} is retained only for
this comparison; it is not used in any proof below, its role in
Theorem~\ref{thm:predictable} having been taken over by
Corollary~\ref{cor:uniform-ui-representation}.

For a \emph{single} non-linear expectation and a completely irregular payoff,
Grigorova, Imkeller, Ouknine and Quenez \cite{Grigorova2020} aggregate the
value family, characterise it as the $\cE^{f}$-Snell envelope, give the
infinitesimal characterisation through an irregular reflected BSDE with
Skorokhod conditions, and apply it to the superhedging of American options
with irregular payoff in an imperfect market model with jumps
\cite[Section~11]{Grigorova2020}.  The present paper replaces one non-linear
expectation by the family $\cM(X)$ of an incomplete market; that replacement
is what gives the regulator a component that is not flat off the contact set
and gives Proposition~\ref{pro:no-skorokhod} its content.

\paragraph{Constrained and second-order reflected BSDEs.}
Selecting the solution of a constrained backward equation by minimality is
standard: Cvitani\'c, Karatzas and Soner \cite{CvitanicKaratzasSoner1998}
treat constraints on the gains process, Peng \cite{Peng1999} obtains smallest
constrained $g$-supersolutions from his monotonic limit theorem, and Peng and
Xu \cite{PengXu2010} treat reflected BSDEs with a constraint on the integrand,
carrying the reflection in a separate increasing process that satisfies a
generalised Skorokhod condition (Remark~\ref{rem:pengxu}).  In the
non-dominated theory, Matoussi, Possama\"i and Zhou
\cite{MatoussiPossamaiZhou2013} construct second-order reflected BSDEs with a
common value and integrand and a family of increasing processes $K^{\Prob}$
indexed by the priors; they point out that $K^{\Prob}$ does not satisfy the
Skorokhod condition, since it combines a second-order and a reflection
component \cite[Remark~3.2, p.~2429]{MatoussiPossamaiZhou2013}, and they leave its
splitting open \cite[Remark~3.5, p.~2433]{MatoussiPossamaiZhou2013}.  Their minimality condition must
be read with the corrigendum \cite{MPZcorrigendum}, which also shows that, under
an assumption on the oscillation of the obstacle, it is equivalent to a
Skorokhod-type condition in which an essential infimum over priors is taken.
Our Definition~\ref{def:rbsde} is a dominated counterpart, with a single
aggregated regulator; Remark~\ref{rem:2bsde} states the comparison.

\paragraph{Interpretation.}
The family $\cM(X)$ is a duality object generated by market incompleteness;
the functional $\sup_Q\sup_\tau$ is the classical superhedging price of an
American claim \cite{ElKarouiQuenez1995,Kramkov96}, and the supremum over
measures is there because the seller must be covered in every scenario, not
because a decision maker is ambiguity averse.  We do not treat Knightian
uncertainty in the sense of Chen and Epstein \cite{CE2002}, Riedel
\cite{Riedel2009} or Nutz \cite{Nutz2013}, where the set of priors is a
primitive of the preferences; see Remark~\ref{rem:knightian}.

\subsection{Organisation}

Section~\ref{sec:prelim} fixes notation and recalls the Mertens decomposition
and the $m$-stability of $\cM(X)$.  Section~\ref{sec:decomposition} contains the
uniform optional decomposition and its predictable refinement.
Section~\ref{sec:snell} develops the generalised Snell envelope.
Section~\ref{sec:bsde} treats reflected BSDEs under the family $\cM(X)$.
Section~\ref{sec:financial} gives the application to robust superhedging, and
Section~\ref{sec:conclusion} lists the open problems that the present techniques
leave untouched.

\section{Preliminaries}\label{sec:prelim}

\subsection{The probabilistic setting}

Throughout, $T>0$ is a fixed horizon and
$(\Omega,\cF,\F=(\cF_t)_{0\le t\le T},\Prob)$ is a filtered probability space
satisfying the usual conditions, with $\cF_0$ trivial up to null sets.  We denote by $\cP$, $\mathrm{Prog}$ and $\cO$
the predictable, progressive and optional $\sigma$-fields on
$\Omega\times[0,T]$, by $\cT$ the set of stopping times with values in $[0,T]$
and, for $\theta\in\cT$, by $\cT_\theta$ the set of $\tau\in\cT$ with
$\tau\ge\theta$ a.s.

A process $y$ is \ladlag{} if it admits finite left and right limits at every
point; we write $y_{t-}$ and $y_{t+}$ for these limits (with $y_{0-}:=y_0$ and
$y_{T+}:=y_T$; the increasing processes $C$ of
Theorem~\ref{thm:uniform} and Definition~\ref{def:rbsde} are the one
exception, being normalised by $C_{0-}=0$), and
\[
\Dm y_t:=y_t-y_{t-},\qquad \Dp y_t:=y_{t+}-y_t
\]
for the left and right jumps.  Recall that a \ladlag{} path has at most
countably many discontinuities.  Processes are identified up to
indistinguishability, inequalities between optional processes hold outside an
evanescent set, and integrands $H$ are identified when their stochastic
integrals $\int H\,dX$ coincide.

Let $X=(X^1,\dots,X^d)$ be a $d$-dimensional \cadlag{} \emph{locally bounded}
semimartingale, and let
\[
\cM(X):=\bigl\{Q\ \text{probability measure on }\cF_T:\ Q\sim\Prob\ \text{and}\
X\ \text{is a }Q\text{-local martingale}\bigr\}.
\]
We assume throughout that
\begin{equation}\label{eq:standing}
\Prob\in\cM(X),
\end{equation}
which is a normalisation (it says that $X$ is already a $\Prob$-local
martingale, as is the case for discounted prices under a reference martingale
measure) and in particular $\cM(X)\neq\emptyset$, which by the fundamental
theorem of asset pricing \cite{DelbaenSchachermayer1994} is the absence of
arbitrage in the sense of NFLVR.  Since all the elements of $\cM(X)$ are equivalent to $\Prob$, the classes of
null sets and of evanescent sets do not depend on which $Q\in\cM(X)$ is used,
and we write ``almost surely'' without naming the measure.  We do, however,
distinguish an identity holding almost surely at each fixed time from one
holding simultaneously at all times outside an evanescent set; the difference
is the subject of Section~\ref{sec:bsde}.

For $Q\in\cM(X)$ we write $\Lambda^{Q}$ for the density process
$\Lambda^{Q}_t=\frac{dQ}{d\Prob}\big|_{\cF_t}$ and
\[
\lambda^{Q}:=\int_0^{\cdot}\frac{d\Lambda^{Q}_s}{\Lambda^{Q}_{s-}}
\]
for its stochastic logarithm, so that $\Lambda^{Q}=\cE(\lambda^{Q})$, where
$\cE$ denotes the Dol\'eans-Dade exponential; here $\Lambda^{Q}_0=1$ because
$\cF_0$ is trivial.  By \cite[Lemma~2.1]{Jacka2012},
$\lambda^{Q}$ is a $\Prob$-local martingale whose jumps are strictly bounded
below by $-1$.  We set
\begin{equation}\label{eq:defL}
\cL:=\bigl\{\lambda^{Q}\ :\ Q\in\cM(X)\bigr\}.
\end{equation}

\begin{defn}\label{def:strong-supermart}
Let $Q\in\cM(X)$.  An optional process $S=(S_t)_{0\le t\le T}$ is a
\emph{strong optional $Q$-supermartingale} if $S_\tau\in L^1(Q)$ for every
$\tau\in\cT$ and
\[
S_\sigma\ \ge\ \E^{Q}\bigl[S_\tau\mid\cF_\sigma\bigr]\quad Q\text{-a.s.}
\qquad\text{for all }\sigma,\tau\in\cT\ \text{with}\ \sigma\le\tau\ Q\text{-a.s.}
\]
It is of \emph{class~(D)} under $Q$ if the family $\{S_\tau:\tau\in\cT\}$ is
uniformly $Q$-integrable.
\end{defn}

Every strong optional supermartingale is \ladlag{} up to evanescence
(see \cite[Appendix~A]{Grigorovaetal2017} or \cite{Dellacherie-Meyer-1980}); we
use this repeatedly.

\begin{thm}[Mertens' decomposition; {\cite{Mertens1972}},
{\cite[Theorem~A.1]{Grigorovaetal2017}}]\label{thm:mertens}
Let $Q\in\cM(X)$ and let $S$ be a strong optional $Q$-supermartingale of
class~(D).  There exist a unique \cadlag{} uniformly integrable
$Q$-martingale $M^{Q}$, a unique predictable right-continuous non-decreasing
process $A^{Q}$ with $A^{Q}_0=0$ and $\E^{Q}[A^{Q}_T]<\infty$, and a unique
right-continuous adapted purely discontinuous non-decreasing process $C^{Q}$
with $C^{Q}_{0-}=0$, $C^{Q}_T=C^{Q}_{T-}$ and $\E^{Q}[C^{Q}_T]<\infty$, such
that
\begin{equation}\label{eq:mertens}
S_t=M^{Q}_t-A^{Q}_t-C^{Q}_{t-},\qquad 0\le t\le T,\quad Q\text{-a.s.}
\end{equation}
\end{thm}

Because $M^{Q}$ and $A^{Q}$ are right-continuous while $t\mapsto C^{Q}_{t-}$ is
left-continuous, \eqref{eq:mertens} gives at once the identity
\begin{equation}\label{eq:rightjump}
\Dp S_t=-\Dm C^{Q}_t=-\bigl(C^{Q}_t-C^{Q}_{t-}\bigr),
\qquad 0\le t<T,\quad Q\text{-a.s.},
\end{equation}
which is the starting point of Section~\ref{sec:decomposition}.  The
normalisation $C^{Q}_T=C^{Q}_{T-}$ is needed for uniqueness: \eqref{eq:mertens}
involves $C^{Q}$ only through $C^{Q}_{t-}$, $t\le T$, so a terminal atom of
$C^{Q}$ would be invisible to it.  It is distinct from the convention
$S_{T+}:=S_T$.  We impose the same normalisation on every purely discontinuous
increasing process below.

\subsection{Stability under pasting}

The following classical property of $\cM(X)$ is the only structural input needed
for the $\esssup$ over measures to behave well, and is used throughout
Sections~\ref{sec:snell}--\ref{sec:financial}.

\begin{defn}\label{def:mstable}
A family $\cQ$ of probability measures equivalent to $\Prob$ is
\emph{$m$-stable} (stable under pasting) if for all $Q_1,Q_2\in\cQ$ and all
$\tau\in\cT$, the measure $Q^{*}$ defined by
\begin{equation}\label{eq:pasting}
Q^{*}(B):=\E^{Q_1}\bigl[\,Q_2(B\mid\cF_\tau)\,\bigr],\qquad B\in\cF_T,
\end{equation}
belongs to $\cQ$.  Equivalently, the density process of $Q^*$ is
$\Lambda^{Q_1}_{t\wedge\tau}\cdot
\Lambda^{Q_2}_{t\vee\tau}/\Lambda^{Q_2}_{\tau}$.
\end{defn}

\begin{pro}\label{pro:mstable}
$\cM(X)$ is $m$-stable and convex.  In particular $\cM(X)$ is either a singleton
or infinite: there is no family with finitely many, but more than one, element.
\end{pro}

\begin{proof}
Convexity is immediate, a convex combination of densities of martingale measures
being the density of a martingale measure; and if $Q_1\neq Q_2$ both lie in
$\cM(X)$, convexity places the whole segment
$\{(1-\alpha)Q_1+\alpha Q_2:\alpha\in[0,1]\}$ in $\cM(X)$, which is therefore
infinite.  This proves the last assertion.  For $m$-stability, let
$Q_1,Q_2\in\cM(X)$, $\tau\in\cT$ and let $Q^{*}$ be given
by~\eqref{eq:pasting}.  Its density process is
$\Lambda_t:=\Lambda^{Q_1}_{t\wedge\tau}\Lambda^{Q_2}_{t\vee\tau}
/\Lambda^{Q_2}_{\tau}$, which is a strictly positive
$\Prob$-martingale with $\E^{\Prob}[\Lambda_T]=1$; hence $Q^{*}\sim\Prob$.
Let $(\rho_n)$ localise $X$ to a bounded process.  For $\sigma\le\sigma'$ in
$\cT$, the identity $\E^{Q^*}[X_{\sigma'\wedge\rho_n}\mid\cF_\sigma]
=X_{\sigma\wedge\rho_n}$ holds on $\{\sigma'\le\tau\}$ because $Q^{*}=Q_1$ on
$\cF_\tau$, and on $\{\sigma\ge\tau\}$ because
$\E^{Q^*}[\,\cdot\mid\cF_\tau]=\E^{Q_2}[\,\cdot\mid\cF_\tau]$; the general case
follows by splitting at $\tau$.  Hence $X$ is a $Q^{*}$-local martingale.
\end{proof}

Pasting at a stopping time implies pasting on events, which is the form used
for directedness.

\begin{lem}[Eventwise pasting]\label{lem:event-paste}
Let $Q_0,Q_1,Q_2\in\cM(X)$, $\theta\in\cT$ and $B\in\cF_\theta$.  The
measure $Q^{*}$ with density
\begin{equation}\label{eq:event-paste}
\frac{dQ^{*}}{d\Prob}
=\Lambda^{Q_0}_\theta\Bigl(\ind_B\,\frac{\Lambda^{Q_1}_T}{\Lambda^{Q_1}_\theta}
+\ind_{B^{c}}\,\frac{\Lambda^{Q_2}_T}{\Lambda^{Q_2}_\theta}\Bigr)
\end{equation}
belongs to $\cM(X)$.  It agrees with $Q_0$ on $\cF_\theta$, and
$\E^{Q^{*}}[\,\cdot\mid\cF_\theta]$ equals $\E^{Q_1}[\,\cdot\mid\cF_\theta]$
on $B$ and $\E^{Q_2}[\,\cdot\mid\cF_\theta]$ on $B^{c}$.
\end{lem}

\begin{proof}
Let $\theta_B:=\theta\ind_B+T\ind_{B^{c}}$ and
$\theta_{B^c}:=\theta\ind_{B^{c}}+T\ind_{B}$, which are stopping times because
$B\in\cF_\theta$.  Let $Q^{a}$ paste $Q_0$ and $Q_1$ at $\theta_B$, and let
$Q^{*}$ paste $Q^{a}$ and $Q_2$ at $\theta_{B^c}$, in the sense of
Definition~\ref{def:mstable}.  Two applications of
Proposition~\ref{pro:mstable} give $Q^{a},Q^{*}\in\cM(X)$.  On $B^{c}$ one has
$\theta_B=T$, so the density process of $Q^{a}$ equals $\Lambda^{Q_0}$ there;
on $B$ one has $\theta_{B^c}=T$, so the density process of $Q^{*}$ equals that
of $Q^{a}$ there.  Reading off the terminal densities on $B$ and on $B^{c}$
gives \eqref{eq:event-paste}, and the conditional statements follow from the
form of the density after $\theta$.
\end{proof}

We shall use $m$-stability through the following standard consequence.

\begin{lem}\label{lem:directed}
Let $\theta\in\cT$ and let $\{Z^{Q}\}_{Q\in\cM(X)}$ be a family of
$\cF_\theta$-measurable random variables of the form
$Z^{Q}=\E^{Q}[\zeta\mid\cF_\theta]$ for a fixed $\zeta\in L^1(Q)$, all
$Q\in\cM(X)$.  Then the family is directed upward, and consequently there is a
sequence $(Q_n)\subset\cM(X)$ with
$Z^{Q_n}\uparrow\esssup_{Q\in\cM(X)}Z^{Q}$ a.s.  Moreover, for every
$Q\in\cM(X)$,
\begin{equation}\label{eq:tower-sup}
\E^{Q}\Bigl[\esssup_{Q'\in\cM(X)}\E^{Q'}[\zeta\mid\cF_\theta]\Bigr]
=\sup\bigl\{\E^{Q'}[\zeta]:\ Q'\in\cM(X),\ Q'=Q\ \text{on}\ \cF_\theta\bigr\}
\ \le\ \sup_{Q'\in\cM(X)}\E^{Q'}[\zeta],
\end{equation}
whenever $\zeta\ge0$; taking the supremum of the left-hand side over
$Q\in\cM(X)$ returns the right-hand side.
\end{lem}

\begin{proof}
Given $Q_1,Q_2$, let $B:=\{Z^{Q_1}\ge Z^{Q_2}\}\in\cF_\theta$ and let $Q^{*}$
be the measure of Lemma~\ref{lem:event-paste} with $Q_0=\Prob$.  Then
$Z^{Q^{*}}=Z^{Q_1}\vee Z^{Q_2}$, which is the directedness.
The existence of the increasing sequence is the standard property of essential
suprema of upward directed families \cite[Appendix~A]{Dellacherie-Meyer-1980}.
For \eqref{eq:tower-sup}, monotone convergence gives
$\E^{Q}[\esssup_{Q'}Z^{Q'}]=\lim_n\E^{Q}[\E^{Q_n}[\zeta\mid\cF_\theta]]
=\lim_n\E^{Q^{*}_n}[\zeta]$,
where $Q^{*}_n$ pastes $Q$ and $Q_n$ at $\theta$ and therefore agrees with $Q$
on $\cF_\theta$; this is ``$\le$''.  For ``$\ge$'', let $Q'\in\cM(X)$ agree
with $Q$ on $\cF_\theta$; then
$\E^{Q'}[\zeta]=\E^{Q}\bigl[\E^{Q'}[\zeta\mid\cF_\theta]\bigr]
\le\E^{Q}[\esssup_{Q'}Z^{Q'}]$.  Finally, the last assertion follows on taking
$Q'=Q$ in the middle term of \eqref{eq:tower-sup}.

The restriction to measures agreeing with $Q$ on $\cF_\theta$ cannot be
removed: at $\theta=T$ the left-hand side is $\E^{Q}[\zeta]$, which is not
$\sup_{Q'}\E^{Q'}[\zeta]$ in an incomplete market.
\end{proof}

\begin{rem}\label{rem:knightian}
Proposition~\ref{pro:mstable} is the reason we work with the \emph{whole} set of
equivalent local martingale measures.  A family of priors in the sense of
Knightian uncertainty --- for instance the family
$\{Q^{\vartheta}:\ dQ^{\vartheta}/d\Prob=\cE(\int\vartheta\,dW)_T,\
|\vartheta|\le\kappa\}$ of Chen and Epstein \cite{CE2002} --- is also
$m$-stable, so Section~\ref{sec:snell} would go through verbatim for it; what
fails for such a family is Theorem~\ref{thm:uniform}, because the
representation of uniform local martingales as stochastic integrals of $X$ is
specific to the full set $\cM(X)$.  See Remark~\ref{rem:full-family} and
Section~\ref{sec:conclusion}.
\end{rem}

\section{The uniform optional decomposition}\label{sec:decomposition}

Throughout this section $S=(S_t)_{0\le t\le T}$ is an optional process and we
assume
\begin{equation}\label{eq:H-super}
S\ \text{is a strong optional }Q\text{-supermartingale of class~(D)
for every }Q\in\cM(X).
\end{equation}

\subsection{Aggregation of the right-jump part}

The following lemma is what makes the decomposition uniform.  In
\eqref{eq:mertens} the martingale $M^{Q}$ and the predictable part $A^{Q}$
genuinely depend on $Q$, whereas the purely discontinuous part is determined by
$\Dp S=-\Dm C^{Q}$, an identity involving only the paths of $S$; it is therefore
the same for all $Q$.  Some care is needed, a pathwise sum of right jumps not
being a priori convergent: we obtain convergence from the decomposition under
$\Prob$ alone, and only then let $Q$ vary.

\begin{lem}[Aggregation of the purely discontinuous part]\label{lem:aggreg-C}
Assume \eqref{eq:H-super}.  There exists a right-continuous adapted purely
discontinuous non-decreasing process $C$ with $C_{0-}=0$ and $C_T=C_{T-}$ such
that, for every $Q\in\cM(X)$,
\[
C=C^{Q}\quad\text{up to indistinguishability},\qquad
\E^{Q}\bigl[C_T\bigr]<\infty,
\]
where $C^{Q}$ is the process appearing in the Mertens decomposition of $S$ under
$Q$.  Moreover
\begin{equation}\label{eq:C-explicit}
C_t=-\sum_{s\le t,\ s<T}\Dp S_s,\qquad 0\le t\le T,
\end{equation}
the sum being over the at most countably many $s\le t$, $s<T$, with
$\Dp S_s\neq0$, and it converges absolutely.
\end{lem}

\begin{proof}
Since $S$ is a strong optional supermartingale it is \ladlag{} up to
evanescence, so $S_{+}$ and hence $\Dp S=S_{+}-S$ are optional processes, and
for every $\omega$ the set $\{s\in[0,T]:\Dp S_s(\omega)\neq0\}$ is at most
countable.

Apply Theorem~\ref{thm:mertens} under $\Prob\in\cM(X)$ (recall
\eqref{eq:standing}) and set $C:=C^{\Prob}$.  Then $C$ is right-continuous,
adapted, non-decreasing, purely discontinuous, $C_{0-}=0$ and
$\E^{\Prob}[C_T]<\infty$, $C_T=C_{T-}$.  By \eqref{eq:rightjump},
$\Dm C_t=-\Dp S_t$ for all $t<T$, $\Prob$-a.s., and $\Dm C_T=0$; since $C$ is
purely discontinuous with $C_{0-}=0$, it is the sum of its jumps, so that
$\Prob$-a.s.
\[
C_t=\sum_{s\le t}\Dm C_s=-\sum_{s\le t,\ s<T}\Dp S_s ,
\]
the convergence being absolute because $C_T<\infty$ and every term is
non-negative.  This proves \eqref{eq:C-explicit}.

Now fix an arbitrary $Q\in\cM(X)$ and let $C^{Q}$ be the purely discontinuous
part of the Mertens decomposition of $S$ under $Q$, which exists by
\eqref{eq:H-super}.  By \eqref{eq:rightjump} applied under $Q$,
$\Dm C^{Q}_t=-\Dp S_t=\Dm C_t$ for all $t<T$, $Q$-a.s., and both processes
have no jump at $T$.  Both $C^{Q}$ and $C$ are purely discontinuous,
non-decreasing and vanish at $0-$, hence both are the sums of their jumps and therefore $C^{Q}_t=C_t$ for all $t$, $Q$-a.s.  As $Q\sim\Prob$
and both processes are right-continuous, they are indistinguishable.  Finally
$\E^{Q}[C_T]=\E^{Q}[C^{Q}_T]<\infty$ by Theorem~\ref{thm:mertens}.
\end{proof}

Lemma~\ref{lem:aggreg-C} gives $\E^{Q}[C_T]<\infty$ for each $Q$ separately.
For the \emph{full} family of a locally bounded market this already upgrades to
a uniform bound on first moments, by the following elementary observation.

\begin{lem}[Countable mixtures]\label{lem:countable-mixtures}
$\cM(X)$ is closed under countable convex mixtures: if $Q_n\in\cM(X)$,
$a_n>0$ and $\sum_na_n=1$, then $\overline Q:=\sum_na_nQ_n\in\cM(X)$.
Consequently, if $Z\ge0$ satisfies $\E^{Q}[Z]<\infty$ for every $Q\in\cM(X)$,
then
\begin{equation}\label{eq:uniform-first-moment}
\sup_{Q\in\cM(X)}\E^{Q}[Z]<\infty .
\end{equation}
\end{lem}

\begin{proof}
The density $d\overline Q/d\Prob=\sum_na_n\,dQ_n/d\Prob$ is strictly positive
with mean one, so $\overline Q\sim\Prob$.  Let $(\rho_k)$ be a localising
sequence for which $X^{\rho_k}$ is bounded, which exists because $X$ is locally
bounded and may be chosen independently of the measure.  For each $n$, the
bounded process $X^{\rho_k}$ is a $Q_n$-local martingale, hence a true
$Q_n$-martingale; so for $s\le t$ and $B\in\cF_s$,
$\E^{\overline Q}[\ind_B(X^{\rho_k}_t-X^{\rho_k}_s)]
=\sum_na_n\E^{Q_n}[\ind_B(X^{\rho_k}_t-X^{\rho_k}_s)]=0$, the interchange being
legitimate by boundedness.  Thus $X^{\rho_k}$ is a $\overline Q$-martingale for
every $k$ and $X$ is a $\overline Q$-local martingale.

For \eqref{eq:uniform-first-moment}, suppose the supremum were infinite and
choose $Q_n\in\cM(X)$ with $\E^{Q_n}[Z]\ge2^{n}$.  Then
$\overline Q:=\sum_{n\ge1}2^{-n}Q_n\in\cM(X)$ and
$\E^{\overline Q}[Z]=\sum_n2^{-n}\E^{Q_n}[Z]\ge\sum_n1=\infty$, contradicting
integrability under $\overline Q$.
\end{proof}

\begin{rem}[First moments and tails are different things]\label{rem:no-uniform-C}
Applying \eqref{eq:uniform-first-moment} to $C_T$, and after
Theorem~\ref{thm:uniform} to $A_T+C_T$, gives
$\sup_{Q}\E^{Q}[A_T+C_T]<\infty$ with no bound on $S$.  The deterministic lower
bound $S\ge-K$ of the earlier form of that theorem supplies the explicit
estimate $S_0+K$ through \eqref{eq:total-L1}, but it is not needed for
finiteness.

Two warnings.  First, a uniform bound on first moments is \emph{not} a uniform
tail estimate: \eqref{eq:uniform-first-moment} says nothing about the size of
$\sup_{Q}\E^{Q}[Z\ind_{\{Z>L\}}]$ as $L\to\infty$, and the distinction is
exactly the one between Assumptions~\ref{assum:E} and~\ref{assum:G} in
Section~\ref{sec:snell}.  Second, the same argument applied to a payoff shows
that class~(D) under every $Q$ separately implies
$\sup_{Q,\tau}\E^{Q}|\xi_\tau|<\infty$ --- if not, choose $(Q_n,\tau_n)$ with
$\E^{Q_n}|\xi_{\tau_n}|\ge4^{n}$, and the mixture $\sum_n2^{-n}Q_n$ has
$\E^{\overline Q}|\xi_{\tau_n}|\ge2^{n}$, contradicting even
$L^{1}$-boundedness under $\overline Q$ --- but this does \emph{not} establish
class~(D) of the robust envelope of Section~\ref{sec:snell}, whose uniform
integrability needs the separate argument given there.  Both statements use
that $\cM(X)$ is the full family of a locally bounded market; an arbitrary
stable subfamily need not be closed under countable mixtures.
\end{rem}

\subsection{The two inputs from the optional decomposition theory}

Recall that $X$ is locally bounded and $\cM(X)\neq\emptyset$.  Kramkov's
theorem is the classical form of the optional decomposition, and we state it
for orientation; the input actually used below is the unrestricted version of
F\"ollmer and Kabanov, which requires no bound on the process being
decomposed.

\begin{thm}[Kramkov \cite{Kramkov96}; see also
\cite{StrickerYan1998}]\label{thm:kramkov}
Let $Y$ be a non-negative \cadlag{} process which is a $Q$-supermartingale for
every $Q\in\cM(X)$.  There exist a predictable $X$-integrable process $H$ and an
adapted right-continuous non-decreasing process $A$ with $A_0=0$ such that
\[
Y_t=Y_0+\int_0^tH_s\,dX_s-A_t,\qquad 0\le t\le T,\quad\text{a.s.}
\]
\end{thm}

\begin{thm}[F\"ollmer--Kabanov {\cite[Theorem~1]{FollmerKabanov1998}}]
\label{thm:fk}
Let $Y$ be a right-continuous process which is a $Q$-local supermartingale for
every $Q\in\cM(X)$.  There exist a predictable $X$-integrable process $H$ and an
adapted right-continuous non-decreasing process $A$ with $A_0=0$ such that
\[
Y_t=Y_0+\int_0^tH_s\,dX_s-A_t,\qquad 0\le t\le T,\quad\text{a.s.}
\]
\end{thm}

No lower bound on $Y$ is assumed in Theorem~\ref{thm:fk}, and this is what
makes Theorem~\ref{thm:uniform} below available at class~(D) alone.  What is
lost with the bound is the \emph{admissibility} of $H$: \cite{FollmerKabanov1998}
records that $H$ is admissible when $Y$ is bounded below, and without that
hypothesis the uniformly integrable martingale property of $\int H\,dX$ has to
be established separately.  Lemma~\ref{lem:BS} is the tool for that, and
Step~4 of the proof of Theorem~\ref{thm:uniform} is where it is used.  The
process $A$ is optional and is not asserted to be unique, in
Theorem~\ref{thm:fk} any more than in Theorem~\ref{thm:kramkov}.

\begin{lem}[Biagini--S\^irbu {\cite[Lemma~A1]{BiaginiSirbu2012}}]\label{lem:BS}
Let $N$ be a $Q$-sigma-martingale.  Then $N$ is a true $Q$-supermartingale if
and only if the family $\{N^{-}_\tau:\tau\in\cT\}$ is uniformly
$Q$-integrable, and in that case $N$ is also a $Q$-local martingale.  (The
horizon $T$ is finite throughout, which is what makes the two conditions
equivalent: a supermartingale on $[0,T]$ satisfies
$N^{-}_\tau\le\E^{Q}[N^{-}_T\mid\cF_\tau]$.)
\end{lem}

\begin{proof}[Proof in the case used below]
The lemma is used only for $N=\int_0^{\cdot}\vartheta\,dX'$, with $X'$ a
vector-valued $Q$-local martingale and $\vartheta\in L(X')$ --- for $X'=X$,
and, in the proof of Theorem~\ref{thm:predictable}, for $X'=(X,D)$ and
$\vartheta=(H,-1)$ --- and we include the short argument in that case.  If $N$ is a $Q$-supermartingale, then
$N^{-}_\tau\le\E^{Q}[N^{-}_T\mid\cF_\tau]$ for every $\tau\in\cT$, so the family
is uniformly integrable.  Conversely, let $\rho_n:=\inf\{t:|N_t|\ge n\}\wedge T$.
On $]\!]0,\rho_n]\!]$ we have $|N_{-}|\le n$, so
$\Delta N^{\rho_n}\ge-(2n+N^{-}_{\rho_n})$ and $N^{\rho_n}\ge-(n+N^{-}_{\rho_n})$,
integrable lower bounds; by
\cite[Proposition~3.3]{AnselStricker1994} the stopped integral $N^{\rho_n}$ is a
$Q$-local martingale, and hence so is $N$.

Choose an increasing localising sequence $(\pi_k)$ such that each
$N^{\pi_k}$ is a true $Q$-martingale and $Q(\pi_k=T)\to1$.  Since
$(\pi_k)$ is non-decreasing and bounded by $T$, the events $\{\pi_k=T\}$
increase with $k$, and their union has $Q$-probability
$\lim_kQ(\pi_k=T)=1$; thus $\pi_k=T$ eventually, $Q$-a.s.  Uniform integrability gives
$K_Q:=\sup_{\theta\in\cT}\E^{Q}[N^{-}_\theta]<\infty$.  Since $N_0=0$,
optional sampling for $N^{\pi_k}$ and Fatou's lemma give, for every
$\tau\in\cT$,
\[
\E^{Q}[N^{+}_\tau]
\le\liminf_k\E^{Q}[N^{+}_{\tau\wedge\pi_k}]
=\liminf_k\E^{Q}[N^{-}_{\tau\wedge\pi_k}]\le K_Q .
\]
Consequently $N_\tau\in L^1(Q)$.  Moreover,
$N^{-}_{\tau\wedge\pi_k}\to N^{-}_\tau$ in $L^1(Q)$, by eventual equality
and uniform integrability.

Fix $\sigma\le\tau$ in $\cT$.  The conditional expectations of these
negative parts converge in $L^1(Q)$, so choose a subsequence along which they
converge $Q$-a.s.  Along that subsequence, conditional Fatou applied only to
the non-negative positive parts, followed by optional sampling, yields
\[
\begin{aligned}
\E^{Q}[N^{+}_\tau\mid\cF_\sigma]
&\le\liminf_k\E^{Q}[N^{+}_{\tau\wedge\pi_k}\mid\cF_\sigma]\\
&=\liminf_k\Bigl(N_{\sigma\wedge\pi_k}
+\E^{Q}[N^{-}_{\tau\wedge\pi_k}\mid\cF_\sigma]\Bigr)
=N_\sigma+\E^{Q}[N^{-}_\tau\mid\cF_\sigma] .
\end{aligned}
\]
Subtracting the last conditional expectation gives
$\E^{Q}[N_\tau\mid\cF_\sigma]\le N_\sigma$, as required.
\end{proof}

A stochastic integral of a local martingale is a sigma-martingale, so
Lemma~\ref{lem:BS} applies to $\int H\,dX$ under every $Q\in\cM(X)$.  It
generalises the Ansel--Stricker lemma \cite{AnselStricker1994}, which is the
case of a deterministic lower bound; \cite[Lemma~A2]{BiaginiSirbu2012} gives
the intermediate form in which $N\ge-W$ for an integrable $W>0$, with the
conclusion that $N$ is a local martingale and a supermartingale.  Both lemmas
are stated in \cite{BiaginiSirbu2012} under a single probability measure, which
is how they are applied here, under each $Q\in\cM(X)$ separately.

\begin{lem}[Uniform local martingales are stochastic
integrals]\label{lem:unif-loc-mart}
Let $N$ be a \cadlag{} process with $N_0=0$, bounded below by a constant, which
is a $Q$-local martingale for every $Q\in\cM(X)$.  Then there exists a
predictable $X$-integrable process $H$ with $N=\int_0^{\cdot}H_s\,dX_s$.
\end{lem}

\begin{proof}
Being a local martingale bounded below, $N$ is a $Q$-supermartingale for every
$Q$, so Theorem~\ref{thm:fk} yields
$N=\int_0^{\cdot}H\,dX-A$ with $A$ adapted, non-decreasing, $A_0=0$.  Then
$\int_0^{\cdot}H\,dX=N+A\ge -K$ is a stochastic integral of a locally bounded
local martingale which is bounded below, hence a $Q$-local martingale by the
Ansel--Stricker lemma \cite{AnselStricker1994}.  Consequently
$A=\int_0^{\cdot}H\,dX-N$ is a $\Prob$-local martingale; being non-decreasing
and null at $0$, it satisfies $\E^{\Prob}[A_{\rho_n}]=0$ along a localising
sequence $(\rho_n)$, whence $A_{\rho_n}=0$ a.s.\ for every $n$ and $A\equiv0$.
\end{proof}

\subsection{The decomposition}

\begin{thm}[Uniform optional decomposition]\label{thm:uniform}
Let $S$ be an optional process satisfying \eqref{eq:H-super}.  No bound on $S$,
from either side, is assumed.  Then there exist
\begin{itemize}\itemsep2pt
\item a predictable $X$-integrable process $H=(H^{i})_{1\le i\le d}$,
\item an adapted right-continuous non-decreasing process $A$ with $A_0=0$,
\item the right-continuous adapted purely discontinuous non-decreasing process
      $C$ of Lemma~\ref{lem:aggreg-C}, with $C_{0-}=0$ and $C_T=C_{T-}$,
\end{itemize}
none of which depends on $Q$, such that
\begin{equation}\label{eq:uniform-decomp}
S_t=S_0+\int_0^tH_s\,dX_s-A_t-C_{t-},\qquad 0\le t\le T,\quad\text{a.s.}
\end{equation}
Moreover, for every $Q\in\cM(X)$, $\int_0^{\cdot}H\,dX$ is a uniformly
integrable $Q$-martingale, $\E^{Q}[A_T+C_T]<\infty$, and
\begin{equation}\label{eq:total-L1}
\E^{Q}\bigl[A_T+C_T\bigr]=S_0-\E^{Q}[S_T].
\end{equation}
No uniqueness of the adapted process $A$ is asserted; if $A$ can be chosen
predictable, then $A$ and $\int H\,dX$ are unique
(Theorem~\ref{thm:predictable}).
\end{thm}

The proof applies Theorem~\ref{thm:fk} to $Y:=S+C_-$, and we record at the
outset where each of its hypotheses is checked: $X$ is a locally bounded
\cadlag{} semimartingale with $\cM(X)\neq\emptyset$, which holds because
$\Prob\in\cM(X)$ by \eqref{eq:standing}, and
$\cM(X)$ is the \emph{entire} set of equivalent local martingale measures of
$X$, which is what the representation of uniform local martingales as
stochastic integrals requires (Remark~\ref{rem:full-family}); $Y$ is
right-continuous (Step~1); and $Y$ is a $Q$-supermartingale, hence a fortiori a
$Q$-local supermartingale, for every $Q\in\cM(X)$ (Step~2).  Class~(D) is not
needed by Theorem~\ref{thm:fk} itself; it is used in Step~4, through
Lemma~\ref{lem:BS}, and through Lemma~\ref{lem:aggreg-C}.

\begin{proof}
Let $C$ be the process of Lemma~\ref{lem:aggreg-C} and set
\[
Y_t:=S_t+C_{t-},\qquad 0\le t\le T .
\]

\emph{Step 1: $Y$ is \cadlag.}  By
Theorem~\ref{thm:mertens} under $\Prob$ and Lemma~\ref{lem:aggreg-C},
$Y=M^{\Prob}-A^{\Prob}$ up to indistinguishability, and the right-hand side is
\cadlag.

\emph{Step 2: $Y$ is a $Q$-supermartingale of class~(D) for every $Q$.}  Fix
$Q\in\cM(X)$.  By Lemma~\ref{lem:aggreg-C} and Theorem~\ref{thm:mertens},
$Y=M^{Q}-A^{Q}$ $Q$-a.s., which is a \cadlag{} $Q$-supermartingale; it is of
class~(D) because $M^{Q}$ is a uniformly integrable $Q$-martingale and
$0\le A^{Q}\le A^{Q}_T\in L^1(Q)$.

\emph{Step 3: decomposition.}  By Theorem~\ref{thm:fk} there
are a predictable $X$-integrable $H$ and an adapted right-continuous
non-decreasing $A$ with $A_0=0$ such that
$Y_t=Y_0+\int_0^tH\,dX-A_t$ for all $t$, a.s.  Since $C_{0-}=0$ we have
$Y_0=S_0$, and substituting $Y=S+C_{-}$ gives \eqref{eq:uniform-decomp}.

\emph{Step 4: integrability.}  Fix $Q\in\cM(X)$ and put
$N:=\int_0^{\cdot}H\,dX$, a $Q$-sigma-martingale.  Here and below we use that
the class $L(X)$ of $X$-integrable predictable processes and the value of
$\int H\,dX$ do not change under an equivalent change of measure, so that $N$
is one and the same process under every $Q\in\cM(X)$.  Since $A\ge0$ and $C_-\ge0$,
\eqref{eq:uniform-decomp} gives $N\ge S-S_0$, hence
\begin{equation}\label{eq:neg-parts}
N^{-}_\tau\ \le\ (S_\tau-S_0)^{-}\qquad\text{for every }\tau\in\cT .
\end{equation}
The family on the right is uniformly $Q$-integrable, because $S$ is of
class~(D) under $Q$ by \eqref{eq:H-super} and $S_0$ is a constant.  By
Lemma~\ref{lem:BS}, $N$ is therefore a $Q$-local martingale and a true
$Q$-supermartingale; in particular $N_T\in L^{1}(Q)$.  It follows that
$A_T+C_T=S_0+N_T-S_T$ is $Q$-integrable, $S_T$ being integrable by class~(D).
This order of argument avoids taking an expectation of $A_T$ before its
integrability is known.  For every $\tau\in\cT$,
\[
|N_\tau|\ \le\ |S_\tau|+|S_0|+A_T+C_T ,
\]
so $N$ is of class~(D) under $Q$; a local martingale of class~(D) is a
uniformly integrable martingale.  Taking expectations in
$A_T+C_T=S_0+N_T-S_T$ gives \eqref{eq:total-L1}, and $\E^{Q}[C_T]<\infty$ by
Lemma~\ref{lem:aggreg-C}.
\end{proof}

\begin{rem}[Where the bound was used, and why it can go]
\label{rem:no-lower-bound}
A deterministic lower bound $S\ge-K$ would serve two distinct purposes in the
proof above, and neither is needed.  It would make $Y$ bounded below, allowing
Kramkov's theorem to be reached through a non-negative shift; here
Theorem~\ref{thm:fk} needs no shift.  And it would make $\int H\,dX$ a local
martingale bounded below, so that the Ansel--Stricker lemma applied;
\eqref{eq:neg-parts} and Lemma~\ref{lem:BS} do that work instead, using
exactly the class~(D) hypothesis already present in \eqref{eq:H-super}.  What
such a bound does buy is an \emph{explicit} estimate: it gives
$\sup_{Q}\E^{Q}[A_T+C_T]\le S_0+K$ through \eqref{eq:total-L1}, whereas
Lemma~\ref{lem:countable-mixtures} gives finiteness of that supremum with no
numerical value attached.
\end{rem}

\begin{cor}[Uniform martingales are stochastic integrals]
\label{cor:uniform-ui-representation}
Let $N$ be a \cadlag{} process with $N_0=0$ which is a uniformly integrable
$Q$-martingale for every $Q\in\cM(X)$.  Then $N=\int_0^{\cdot}H\,dX$ for some
predictable $X$-integrable $H$.
\end{cor}

\begin{proof}
$N$ satisfies \eqref{eq:H-super}: a uniformly integrable martingale is a
supermartingale of class~(D).  Theorem~\ref{thm:uniform} gives
$N=\int_0^{\cdot}H\,dX-A-C_-$.  As $N$ is \cadlag, $\Dp N\equiv0$, so $C\equiv0$
by Lemma~\ref{lem:aggreg-C}.  Then \eqref{eq:total-L1} under $\Prob$ gives
$\E^{\Prob}[A_T]=N_0-\E^{\Prob}[N_T]=0$, and $A$ is non-decreasing with
$A_0=0$, so $A\equiv0$.
\end{proof}

Corollary~\ref{cor:uniform-ui-representation} removes the lower bound from
Lemma~\ref{lem:unif-loc-mart} at the cost of assuming a uniformly integrable
martingale rather than a local one, which is the form needed in
Theorem~\ref{thm:predictable}.

\begin{cor}[Uniform decomposition of submartingales]
\label{cor:uniform-submartingale}
Let $S$ be an optional process which is a strong optional $Q$-submartingale of
class~(D) for every $Q\in\cM(X)$.  Then $\Dp S\ge0$ and
\begin{equation}\label{eq:uniform-sub}
S_t=S_0+\int_0^tH_s\,dX_s+A_t+C_{t-},\qquad
C_t=\sum_{s\le t,\ s<T}\Dp S_s ,
\end{equation}
with $H,A,C$ not depending on $Q$, $A$ adapted right-continuous non-decreasing
with $A_0=0$, $C$ as in Lemma~\ref{lem:aggreg-C} applied to $-S$, and, for every
$Q\in\cM(X)$, $\int_0^{\cdot}H\,dX$ a uniformly integrable $Q$-martingale with
$\E^{Q}[A_T+C_T]=\E^{Q}[S_T]-S_0$.  Under Assumption~\ref{assum:H1} the process
$A$ may be taken predictable and is then unique, together with the gains
process.  No bound on $S$ is assumed.
\end{cor}

\begin{proof}
Apply Theorem~\ref{thm:uniform}, and Theorem~\ref{thm:predictable} for the last
assertion, to $-S$, which is a strong optional $Q$-supermartingale of class~(D)
for every $Q$, and change signs.
\end{proof}

This is the form in which the decomposition is used for upper obstacles and for
the buyer's problem in the companion papers \cite{BA,P4}, and, together with
Theorem~\ref{thm:uniform}, for game options in \cite{GO}; those papers cite it
rather than reproving it.

\begin{rem}[What Theorem~\ref{thm:uniform} adds, statement by statement]
\label{rem:what-adds}
The comparison with the antecedents is best made hypothesis by hypothesis
rather than in prose; Section~\ref{subsec:literature} gives the references.
\begin{itemize}\itemsep1pt
\item \emph{The integrator.}  Here $X$ is an arbitrary locally bounded
\cadlag{} semimartingale, so it may jump and its quadratic variation need not
be absolutely continuous.  In \cite[Theorem~4.1]{BouchardPossamaiTan2016} it is
a continuous martingale $\int\alpha^{\top}dW$ with $\alpha$ bounded with
bounded inverse.
\item \emph{The regularity of $S$.}  Right-continuity is assumed in
\cite{Kramkov96} and is not assumed here, in \cite{BouchardPossamaiTan2016} or
in \cite{Arharas2025}.
\item \emph{The integrability.}  We assume class~(D) under each $Q$ and nothing
else: no bound on $S$ from either side.  By contrast
\cite{BouchardPossamaiTan2016} assume
$\esssup_{\tau}|S(\tau)|\in L^{p}(Q)$ with $p>1$, and
\cite{Arharas2025} assume square integrability.  These are hypotheses on the
running maximum, at a level above $L^{1}$; class~(D) is a condition at stopping
times, at the $L^{1}$ level.  They are different integrability levels rather
than one being stronger throughout, and we claim no superiority over every
hypothesis of those theorems --- only that class~(D) is the regime in which the
present formulation works.
\item \emph{The martingale part.}  For a general family of priors the
martingale part is abstract \cite[Theorem~5.1]{Arharas2025}.  Taking the
\emph{entire} $\cM(X)$ makes it the stochastic integral $\int H\,dX$
\cite{Kramkov96,FollmerKabanov1998}, which is what turns the decomposition into
a trading strategy in Theorem~\ref{thm:robust-american}.
\item \emph{The decomposition theorem used.}  The càdlàg optional decomposition
without a bound on the process is due to F\"ollmer and Kabanov
\cite[Theorem~1]{FollmerKabanov1998} and we claim nothing for it.  What is
formulated here is the strong optional version: the intrinsic right-jump term,
the class~(D) hypothesis at stopping times, and the uniformly integrable
martingale conclusion for $\int H\,dX$, which does not come with that theorem
and is supplied by \cite[Lemma~A1]{BiaginiSirbu2012}.
\item \emph{The finite-variation part.}  The antecedents produce one
non-decreasing process.  Here it is split as $A+C_-$, with $C$ given pathwise
by the right jumps of $S$, which is what connects the statement to Mertens'
theorem and makes the reflected BSDE \eqref{eq:rbsde} and the consumption
stream of Theorem~\ref{thm:robust-american} available.
\end{itemize}
In one sentence: Theorem~\ref{thm:uniform} combines these ingredients into a
single triple $(H,A,C)$ for a general locally bounded \cadlag{} integrator and
a process of class~(D) under every $Q\in\cM(X)$, with the martingale part
represented as a trading gain and the right-jump regulator identified
explicitly; that formulation is the input used in the superhedging and BSDE
results below.  A
further difference, which belongs to Section~\ref{sec:bsde} rather than here,
is that the obstacle is dominated at every stopping time and not
$d\Prob\otimes dt$-a.e.\ (Remark~\ref{rem:pengxu}).
\end{rem}

The converse of Theorem~\ref{thm:uniform} requires an additional uniform
integrability requirement, and is \emph{false} without it.

\begin{pro}[Converse]\label{pro:converse}
Let $H$ be predictable and $X$-integrable, let $A$ and $C$ be adapted
right-continuous non-decreasing processes with $A_0=C_{0-}=0$, and define $S$ by
\eqref{eq:uniform-decomp}.  Assume that for every $Q\in\cM(X)$
\begin{equation}\label{eq:converse-hyp}
\E^{Q}\bigl[A_T+C_T\bigr]<\infty
\qquad\text{and}\qquad
\textstyle\int_0^{\cdot}H\,dX\ \text{is a uniformly integrable }
Q\text{-martingale}.
\end{equation}
Then $S$ satisfies \eqref{eq:H-super}.
\end{pro}

\begin{proof}
Write $N:=\int_0^{\cdot}H\,dX$ and fix $Q\in\cM(X)$.  For $\sigma\le\tau$ in
$\cT$, optional sampling for the uniformly integrable martingale $N$ gives
$\E^{Q}[N_\tau\mid\cF_\sigma]=N_\sigma$, and since $A$ and $C_{-}$ are
non-decreasing,
\[
\E^{Q}[S_\tau\mid\cF_\sigma]
=S_0+N_\sigma-\E^{Q}\bigl[A_\tau+C_{\tau-}\mid\cF_\sigma\bigr]
\le S_0+N_\sigma-A_\sigma-C_{\sigma-}=S_\sigma .
\]
For class~(D), note $|S_\tau|\le|S_0|+|N_\tau|+A_T+C_T$ for every $\tau\in\cT$;
the family $\{N_\tau\}_{\tau\in\cT}$ is uniformly $Q$-integrable because $N$ is
a uniformly integrable $Q$-martingale, and $A_T+C_T\in L^1(Q)$ dominates the
remaining terms.
\end{proof}

\begin{eg}[The converse fails without \eqref{eq:converse-hyp}]
\label{eg:bessel}
Let $R$ be a three-dimensional Bessel process with $R_0=1$ under $\Prob$ and set
$X:=1/R$, a continuous --- hence locally bounded --- strictly positive
$\Prob$-local martingale which is not a martingale:
$\E^{\Prob}[X_T]<X_0$.  Then $\Prob\in\cM(X)$.  Take $S:=X$, $H\equiv1$ and
$A\equiv C\equiv0$, so that \eqref{eq:uniform-decomp} holds identically and
$S\ge0$ is bounded below.  For every $Q\in\cM(X)$, $S$ is a non-negative
$Q$-local martingale and therefore a $Q$-supermartingale.  However $S$ is
\emph{not} of class~(D) under $\Prob$: a non-negative local martingale of
class~(D) is a uniformly integrable martingale, while $\E^{\Prob}[X_T]<X_0$.
Thus \eqref{eq:H-super} fails although the right-hand side of
\eqref{eq:uniform-decomp} is available.  This shows that
Theorem~\ref{thm:uniform} cannot be stated as an equivalence without
\eqref{eq:converse-hyp}.
\end{eg}

\begin{rem}\label{rem:full-family}
Theorem~\ref{thm:uniform} uses that $\cM(X)$ is the \emph{entire} set of
equivalent local martingale measures of $X$, through
Theorem~\ref{thm:fk}.  For a strict subfamily $\cP\subsetneq\cM(X)$ the
hypothesis ``$S$ is a $\cP$-uniform supermartingale'' is weaker and the
conclusion generally fails: a $\cP$-uniform local martingale need not be a
stochastic integral of $X$.  This is why the decomposition of
\cite{Arharas2025}, which allows an arbitrary family of priors, is stated with
an abstract martingale part.  See Section~\ref{sec:conclusion}.
\end{rem}

\subsection{When is the finite-variation part predictable?}

We now give a sufficient condition under which $A$ in
\eqref{eq:uniform-decomp} can be chosen predictable.  Following Jacka
\cite{Jacka2012}, it is stated on the localisation of the set $\cL$ of
\eqref{eq:defL} of stochastic logarithms of the density processes:
\begin{equation}\label{eq:Lloc}
\cL^{\mathrm{loc}}:=\Bigl\{\lambda\ \text{with}\ \lambda_0=0:\
\exists\ \text{stopping times}\ \rho_n\uparrow,\ \Prob(\rho_n=T)\to1,\
\lambda^{\rho_n}\in\cL\ \text{for every }n\Bigr\},
\end{equation}
where $\lambda^{\rho_n}$ denotes the stopped process; clearly
$\cL\subset\cL^{\mathrm{loc}}$.

\begin{assum}[{\bf H1}]\label{assum:H1}
$\cL^{\mathrm{loc}}$ is stable under multiplication by real scalars: for every
$\lambda\in\cL^{\mathrm{loc}}$ and every $\delta\in\R$, one has
$\delta\lambda\in\cL^{\mathrm{loc}}$.
\end{assum}

This is the second hypothesis of \cite[Theorem~2.2]{Jacka2012}, whose first
hypothesis, $\Prob\in\cM(X)$, is our standing assumption \eqref{eq:standing};
$\cL^{\mathrm{loc}}$ is the class denoted $\cL_{\rm loc}$ there.  It is a
sufficient condition, and we make no claim that it is necessary.

\begin{lem}[{\bf H1} forces continuous densities]\label{lem:H1-continuous}
Under Assumption~\ref{assum:H1} every element of $\cL^{\mathrm{loc}}$ is
continuous; consequently every density process $Z^{Q}=\cE(\lambda^{Q})$,
$Q\in\cM(X)$, is continuous.
\end{lem}

\begin{proof}
Let $\lambda\in\cL^{\mathrm{loc}}$.  For every $\delta\in\R$ the stopped
processes $(\delta\lambda)^{\rho^{\delta}_n}$ lie in $\cL$ for a localising
sequence depending on $\delta$, and elements of $\cL$ have jumps strictly
greater than $-1$ \cite[Lemma~2.1]{Jacka2012}.  Hence
$\delta\Dm\lambda>-1$ on $[0,\rho^{\delta}_n]$ and, since
$\Prob(\rho^{\delta}_n=T)\to1$, on $[0,T]$ almost surely, for every
$\delta\in\R$; letting $\delta\to\pm\infty$ forces $\Dm\lambda\equiv0$.  The
stochastic exponential of a continuous process is continuous.
\end{proof}

This observation is due to Jacka \cite{Jacka2012}, and it is what
Theorem~\ref{thm:predictable} now rests on; we state it separately because it
has become load-bearing.

\begin{rem}[What {\bf H1} really says]\label{rem:H1-sharp}
Two comments are essential for the correct use of Assumption~\ref{assum:H1}.

\emph{(i) It is a restriction on the model.}  By
Lemma~\ref{lem:H1-continuous} it forces the density processes to be
continuous.  Assumption~\ref{assum:H1} is therefore a
genuine structural restriction on the model, satisfied for instance in the
It\^o setting of Example~\ref{eg:H1-brownian} and, trivially, when
$\cM(X)=\{\Prob\}$ and $\cL^{\mathrm{loc}}=\{0\}$, in which case
Theorem~\ref{thm:predictable} reduces to Mertens' decomposition.

\emph{(ii) Stability under multiplication by scalars in $[0,1]$ does not
suffice.}  The weaker requirement ``$\alpha\lambda\in\cL^{\mathrm{loc}}$ for
$\alpha\in[0,1]$'' places no restriction on the jumps of $\lambda$, since
$\alpha\Dm\lambda>-1$ holds automatically for $\alpha\in[0,1]$ whenever
$\Dm\lambda>-1$; in the elementary model below it is supplied by convexity of
$\cM(X)$.  Example~\ref{eg:counterexample-H1} exhibits a model in which
it holds and in which no predictable decomposition exists.
\end{rem}

\begin{thm}[Predictable uniform decomposition]\label{thm:predictable}
Let $S$ be as in Theorem~\ref{thm:uniform} and assume in addition
Assumption~\ref{assum:H1}.  Then the process $A$ in \eqref{eq:uniform-decomp}
may be chosen \emph{predictable}, and it is then unique.  The process $C$ is
unchanged; it is adapted and purely discontinuous, and is not predictable in
general.
\end{thm}

Without a lower bound on $S$ the route through \cite[Theorem~2.2]{Jacka2012}
is no longer available: that theorem takes a non-negative process as its input,
and the non-negative shift $Y+K$ was exactly what the deleted hypothesis
supplied.  The proof below compensates the finite-variation part directly, and
uses Assumption~\ref{assum:H1} only through the continuity of the density
processes established in Lemma~\ref{lem:H1-continuous}.

\begin{proof}
Let $Y=S+C_{-}=S_0+N-A$ be as in the proof of Theorem~\ref{thm:uniform}, with
$N=\int_0^{\cdot}H\,dX$ a uniformly integrable $Q$-martingale for every
$Q\in\cM(X)$ and $A$ adapted, right-continuous, non-decreasing, $A_0=0$, and
$A_T\in L^1(Q)$ for every $Q$ by Step~4 there.

\emph{Compensating $A$.}  Let $A^{p}$ be the $\Prob$-predictable compensator of
the integrable increasing process $A$, and put $D:=A-A^{p}$.  Then $D$ is a
uniformly integrable $\Prob$-martingale with paths of finite variation, since
$|D_t|\le A_T+A^{p}_T\in L^{1}(\Prob)$.

\emph{$D$ is a $Q$-local martingale for every $Q$.}  Let $Q\in\cM(X)$ and let
$Z^{Q}$ be its density process.  By Lemma~\ref{lem:H1-continuous}, $Z^{Q}$ is continuous.  As $D$ has finite
variation, $[D,Z^{Q}]=\sum_{s\le\cdot}\Dm D_s\,\Dm Z^{Q}_s=0$, because
$Z^{Q}$ has no jumps.  Integration by parts therefore gives
\[
DZ^{Q}=\int_0^{\cdot}D_{s-}\,dZ^{Q}_s+\int_0^{\cdot}Z^{Q}_{s-}\,dD_s
+[D,Z^{Q}]
=\int_0^{\cdot}D_{s-}\,dZ^{Q}_s+\int_0^{\cdot}Z^{Q}_{s-}\,dD_s ,
\]
and both integrals are $\Prob$-local martingales, the integrands $D_{-}$ and
$Z^{Q}_{-}$ being predictable and locally bounded and the integrators $Z^{Q}$
and $D$ being $\Prob$-martingales.  Hence $DZ^{Q}$ is a $\Prob$-local
martingale, and by the Girsanov--Lenglart criterion --- $D$ is a $Q$-local
martingale if and only if $DZ^{Q}$ is a $\Prob$-local martingale --- the
process $D$ is a $Q$-local martingale.

\emph{The compensated integrand.}  Put $\widetilde N:=N-D$, so that
\[
Y=S_0+N-A=S_0+\widetilde N-A^{p} .
\]
As a difference of a sigma-martingale and a $Q$-local martingale,
$\widetilde N$ is a $Q$-sigma-martingale.  Since $A^{p}\ge0$ we have
$\widetilde N\ge Y-S_0$, hence
$\widetilde N^{-}_\tau\le(Y_\tau-S_0)^{-}$ for every $\tau\in\cT$, and the
right-hand family is uniformly $Q$-integrable because $Y$ is of class~(D) under
$Q$ (Step~2 of Theorem~\ref{thm:uniform}).  By Lemma~\ref{lem:BS},
$\widetilde N$ is a $Q$-local martingale and a true $Q$-supermartingale, so
$\widetilde N_T\in L^{1}(Q)$ and $A^{p}_T=S_0+\widetilde N_T-Y_T\in L^{1}(Q)$.
Then $|\widetilde N_\tau|\le|Y_\tau|+|S_0|+A^{p}_T$ makes $\widetilde N$ of
class~(D) under $Q$, hence a uniformly integrable $Q$-martingale.  As this holds
for every $Q\in\cM(X)$ and $\widetilde N$ is \cadlag{} with
$\widetilde N_0=0$, Corollary~\ref{cor:uniform-ui-representation} gives
$\widetilde N=\int_0^{\cdot}H'\,dX$ for a predictable $X$-integrable $H'$.
Substituting $Y=S+C_{-}$ and using $C_{0-}=0$ gives \eqref{eq:uniform-decomp}
with the predictable process $A^{p}$ in place of $A$.

Uniqueness: if $A$ and $A'$ are two predictable such processes, then
$\int_0^{\cdot}H\,dX-\int_0^{\cdot}H'\,dX=A-A'$ is, for each fixed
$Q\in\cM(X)$, a $Q$-martingale, and is predictable and of
finite variation, hence constant, hence $\equiv0$.  That both integrals are
uniformly integrable $Q$-martingales is not read off
Theorem~\ref{thm:uniform}, which asserts it only for the triple it constructs;
it follows by rerunning Step~4 of that proof, which applies to \emph{any}
decomposition \eqref{eq:uniform-decomp} with a non-decreasing $A'$: then
$N':=\int H'\,dX\ge S-S_0$, so $(N')^{-}_\tau\le(S_\tau-S_0)^{-}$ is uniformly
$Q$-integrable, Lemma~\ref{lem:BS} makes $N'$ a $Q$-local martingale and a
supermartingale, and $|N'_\tau|\le|S_\tau|+|S_0|+A'_T+C_T$ makes it of
class~(D).  Thus $A$ and the gain
process $\int H\,dX$ are unique; $H$ itself is unique only up to the
identification of integrands with equal stochastic integrals, which matters
when $\sigma$ below does not have full rank.
\end{proof}

\begin{rem}
Under Assumption~\ref{assum:H1} the process $A$ of \eqref{eq:uniform-decomp}
coincides, for \emph{every} $Q\in\cM(X)$ simultaneously, with the predictable
part $A^{Q}$ of the Mertens decomposition of $S$ under $Q$.  Indeed
$M^{Q}-A^{Q}=Y=Y_0+\int H\,dX-A$ with $\int H\,dX$ a $Q$-local martingale, and
uniqueness in the Doob--Meyer decomposition of the $Q$-supermartingale $Y$
identifies $A^{Q}=A$.
\end{rem}

\begin{eg}[{Stability under $[0,1]$-scaling is not enough}]
\label{eg:counterexample-H1}
Let $\varepsilon$ be a random variable on $(\Omega,\cF,\Prob)$ with values in
$\{-1,0,1\}$, with $\Prob(\varepsilon=i)>0$ for $i=-1,0,1$ and
$\E^{\Prob}[\varepsilon]=0$.  Put
\[
X_t:=\varepsilon\,\ind_{\{t\ge1/2\}},\qquad
\cF_t:=\{\emptyset,\Omega\}\ (t<1/2),\qquad
\cF_t:=\sigma(\varepsilon)\vee\cN\ (t\ge1/2),
\]
on $[0,T]$ with $T=1$.  Then $X$ is bounded and $\Prob\in\cM(X)$.  A measure
$Q\sim\Prob$ belongs to $\cM(X)$ if and only if its density
$\Lambda:=dQ/d\Prob$ satisfies $\Lambda>0$, $\E^{\Prob}[\Lambda]=1$ and
$\E^{\Prob}[\Lambda\varepsilon]=0$; the density process is $1$ before $1/2$ and
$\Lambda$ afterwards, so that $\lambda^{Q}=(\Lambda-1)\ind_{[1/2,1]}$.

\emph{The $[0,1]$-scaling property holds.}  For $\alpha\in[0,1]$, the measure
with density $1+\alpha(\Lambda-1)$ is the convex combination
$(1-\alpha)\Prob+\alpha Q$, which lies in $\cM(X)$ by
Proposition~\ref{pro:mstable}, is equivalent to $\Prob$, and has stochastic
logarithm exactly $\alpha\lambda^{Q}$.

\emph{No predictable decomposition exists.}  Let
$S_t:=-|X_t|=-|\varepsilon|\ind_{\{t\ge1/2\}}$.  Then $-1\le S\le0$, so $S$ is
bounded below and of class~(D), and for every $Q\in\cM(X)$,
$\E^{Q}[S_{1/2}]=-\E^{Q}|\varepsilon|\le0=S_0$, so $S$ is a strong optional
$Q$-supermartingale.  Also $\Dp S\equiv0$, so $C\equiv0$ by
Lemma~\ref{lem:aggreg-C}.  Suppose $S_t=S_0+\int_0^tH\,dX-A_t$ with $H$
predictable and $A$ predictable non-decreasing, $A_0=0$.  Since $1/2$ is a
deterministic --- hence predictable --- time with $\cF_{(1/2)-}$ trivial, both
$h:=H_{1/2}$ and $a:=A_{1/2}\ge0$ are deterministic, and
$h\varepsilon-a=-|\varepsilon|$ must hold for each of the three values of
$\varepsilon$.  Taking $\varepsilon=0$ forces $a=0$; then $\varepsilon=1$ forces
$h=-1$ and $\varepsilon=-1$ forces $h=1$, a contradiction.

The optional decomposition of Theorem~\ref{thm:uniform} does of course exist:
take $H\equiv0$ and $A_t=|\varepsilon|\ind_{\{t\ge1/2\}}$, which is adapted but
not predictable.  Consistently with Remark~\ref{rem:H1-sharp}(i),
Assumption~\ref{assum:H1} fails here: if $Q\neq\Prob$, then for $\delta$ large
enough the jump $\delta(\Lambda-1)$ of $\delta\lambda^{Q}$ at the deterministic
time $1/2$ is $\le-1$ on an event of positive probability, and this event meets
$\{\rho_n\ge1/2\}$ for $n$ large whenever $\Prob(\rho_n=T)\to1$; so
$\delta\lambda^{Q}\notin\cL^{\mathrm{loc}}$.
\end{eg}

\begin{eg}[An It\^o market satisfying {\bf H1}]\label{eg:H1-brownian}
Let $W$ be an $n$-dimensional Brownian motion, let $\F$ be its usual
augmentation, and let
\[
X_t=\int_0^t\sigma_s\,dW_s
\]
for a predictable $d\times n$ matrix process $\sigma$ with
$\int_0^T\|\sigma_s\|^2ds<\infty$ a.s., $\|\cdot\|$ denoting the Frobenius
norm.  Then $X$ is a continuous local martingale, hence locally bounded, as
Section~\ref{sec:prelim} requires: local boundedness is asked of the integrator
$X$, not of the coefficient $\sigma$.  No bound on $\sigma$, no local
boundedness of $\sigma$ and no relation between $d$ and $n$ is needed anywhere
in what follows, so the geometric case $\sigma_t=(\varsigma_t X_t)$ is covered,
also with a volatility that blows up in time, such as
$\lvert\varsigma_t\rvert=t^{-1/4}$, and so are markets with redundant assets.  A measure $Q\sim\Prob$ has density process
$\cE\bigl(\int\lambda\trans\,dW\bigr)$ for a predictable $\R^{n}$-valued
$\lambda$ with $\int_0^T|\lambda_s|^2ds<\infty$ a.s., and by Girsanov's
theorem $X$ is a $Q$-local martingale if and only if $\sigma_s\lambda_s=0$
$ds\otimes d\Prob$-a.e.; the drift $\int_0^{\cdot}\sigma_s\lambda_s\,ds$ of $X$
under $Q$ is well defined, since
$\int_0^T|\sigma_s\lambda_s|\,ds\le(\int_0^T\|\sigma_s\|^2ds)^{1/2}
(\int_0^T|\lambda_s|^2ds)^{1/2}<\infty$.  We claim
\[
\cL^{\mathrm{loc}}
=\Bigl\{\textstyle\int\lambda\trans\,dW\ :\ \lambda\ \text{predictable},\
\sigma\lambda=0,\ \int_0^T|\lambda_s|^2ds<\infty\ \text{a.s.}\Bigr\}.
\]
The inclusion ``$\subset$'' follows from the kernel condition applied to each
stopped process $\lambda^{\rho_n}\in\cL$.  Conversely, given $\lambda$ in the
right-hand side, let $\rho_n:=\inf\{t:\int_0^t|\lambda_s|^2ds\ge n\}\wedge T$;
then $\Prob(\rho_n=T)\to1$, $\cE(\int\lambda\trans\ind_{[0,\rho_n]}dW)$ is a
strictly positive true martingale by Novikov's criterion, and by Girsanov's
theorem $X$ is a local martingale under the corresponding measure, whose
stochastic logarithm is $(\int\lambda\trans dW)^{\rho_n}$.  The right-hand
side is a vector space, so Assumption~\ref{assum:H1} holds, and all its
elements are continuous, as Remark~\ref{rem:H1-sharp}(i) requires.  This is
the situation of the closing remark of \cite{Jacka2012}: when the filtration
carries only continuous martingales, $\lambda\in\cL^{\mathrm{loc}}$ if and only
if $\langle\lambda,X\rangle=0$, which is stable under real scalars.  The family
$\cM(X)$ is a singleton if and only if $\sigma_s$ has rank $n$
$ds\otimes d\Prob$-a.e., which requires $d\ge n$, i.e.\ if and only if the
market is complete; it is infinite as soon as $\mathrm{rank}\,\sigma_s<n$ on a set of positive
$ds\otimes d\Prob$-measure.  This is the setting of
Section~\ref{sec:bsde}.

In this setting Assumption~\ref{assum:H1} does no work.  In the Brownian
filtration the martingale part $M^{\Prob}$ of the Mertens decomposition and the
integral $\int H\,dX$ are continuous, so the difference $A-A^{\Prob}$ of the
two increasing parts is a continuous local martingale of finite variation and
vanishes; thus the process $A$ of \eqref{eq:uniform-decomp} is the predictable
part $A^{\Prob}$, without any hypothesis, and Theorem~\ref{thm:predictable}
adds nothing to Theorem~\ref{thm:uniform} in Section~\ref{sec:bsde}.  Assumption~\ref{assum:H1} earns its place in
Section~\ref{sec:decomposition}, where $\F$ is arbitrary, and
Example~\ref{eg:counterexample-H1} shows that it is not implied by the convexity-type
closure one might hope to use instead.
\end{eg}

\section{The generalised Snell envelope}\label{sec:snell}

In this section $\eta=(\eta_t)_{0\le t\le T}$ is an optional process, thought of
as the reward of an optimal stopping problem, and we study
\begin{equation}\label{eq:snell-def}
V_\theta:=\esssup_{Q\in\cM(X)}\ \esssup_{\tau\in\cT_\theta}\
\E^{Q}\bigl[\eta_\tau\mid\cF_\theta\bigr],\qquad\theta\in\cT .
\end{equation}
No right-continuity, and no right upper semicontinuity, is assumed on $\eta$;
this is the point of working with strong optional supermartingales.

\begin{assum}[{\bf E}]\label{assum:E}
$\eta$ is optional and
\begin{equation}\label{eq:UI}
\varepsilon(L):=\sup_{Q\in\cM(X)}\ \sup_{\tau\in\cT}\
\E^{Q}\bigl[|\eta_\tau|\ind_{\{|\eta_\tau|>L\}}\bigr]\ \longrightarrow\ 0
\qquad(L\to\infty).
\end{equation}
\end{assum}

Condition \eqref{eq:UI} is uniform integrability of $\{\eta_\tau\}$ jointly over
measures and stopping times.  It is strictly stronger than class~(D) separately
under each measure --- by Remark~\ref{rem:no-uniform-C} the latter gives only
$B:=\sup_{Q,\tau}\E^{Q}|\eta_\tau|<\infty$ --- and no time regularity of $\eta$
is imposed.  Bounded payoffs satisfy it with $\varepsilon(L)=0$ for
$L$ at least the bound on $|\eta|$.  Two sufficient conditions stated directly
on the data are worth
recording.

\begin{lem}[A power moment suffices]\label{lem:Lp}
If $K_p:=\sup_{Q\in\cM(X)}\sup_{\tau\in\cT}\E^{Q}\bigl[|\eta_\tau|^{p}\bigr]
<\infty$ for some $p>1$, then \eqref{eq:UI} holds, with
$\varepsilon(L)\le K_pL^{1-p}$.
\end{lem}

\begin{proof}
On $\{|\eta_\tau|>L\}$ one has $|\eta_\tau|\le|\eta_\tau|^{p}L^{1-p}$.
\end{proof}

The supremum over stopping times in $K_p$ sits \emph{outside} the expectation,
so Lemma~\ref{lem:Lp} is weaker than a bound on
$\E^{Q}\bigl[(\sup_t|\eta_t|)^{p}\bigr]$, and it permits payoffs unbounded on
both sides.  It is not the maximal inequality of Problem~P1.

\begin{assum}[{\bf G}]\label{assum:G}
$\eta$ is optional and $|\eta_t|\le M_t$ for all $t$, where
$M=m+\int_0^{\cdot}G\,dX\ge0$ is a uniformly integrable $Q$-martingale for
every $Q\in\cM(X)$.
\end{assum}

Assumption~\ref{assum:G} is the growth condition of \cite{BA}, and it is an
\emph{alternative} to Assumption~\ref{assum:E}, not equivalent to it: it gives
class~(D) of the envelope directly
(Corollary~\ref{cor:envelope-common-domination}) without giving \eqref{eq:UI},
and Example~\ref{eg:G-not-UI} shows that the implication
{\bf(G)}$\Rightarrow${\bf(E)} genuinely fails.  The converse implication does
hold, and is proved once the envelope theory is available
(Corollary~\ref{cor:ui-traded-domination}).  The resulting hierarchy is
\[
\text{power moment (Lemma~\ref{lem:Lp})}\ \Longrightarrow\
{\bf(E)}\ \Longrightarrow\ {\bf(G)},\qquad
{\bf(G)}\ \not\Longrightarrow\ {\bf(E)} .
\]

\begin{rem}\label{rem:E-relax}
What the envelope theory needs of the data is exactly that $\{V_\theta\}$ be of
class~(D) under every $Q$; the deleted lower bound of
Theorem~\ref{thm:uniform} is no longer required.  Under
Assumption~\ref{assum:E} this is Proposition~\ref{pro:envelope-UI}, and under
Assumption~\ref{assum:G} it is
Corollary~\ref{cor:envelope-common-domination}.  Both \eqref{eq:UI} and the
$L^{p}$ bound of Lemma~\ref{lem:Lp} are conditions on the data alone, which is
why the corresponding part of Problem~P1 is settled.  What remains open there
is the stronger maximal estimate; the non-linear equation of
Section~\ref{sec:bsde} is treated for such data in
Section~\ref{subsec:weighted}.
\end{rem}

\begin{pro}[The envelope is uniformly integrable]\label{pro:envelope-UI}
Let Assumption~\ref{assum:E} hold and put, for $\theta\in\cT$,
\[
R_\theta:=\esssup_{Q\in\cM(X)}\ \esssup_{\tau\in\cT_\theta}\
\E^{Q}\bigl[|\eta_\tau|\mid\cF_\theta\bigr],\qquad
R^{L}_\theta:=\esssup_{Q\in\cM(X)}\ \esssup_{\tau\in\cT_\theta}\
\E^{Q}\bigl[|\eta_\tau|\ind_{\{|\eta_\tau|>L\}}\mid\cF_\theta\bigr].
\]
Then $|V_\theta|\le R_\theta$, $R_\theta\le L+R^{L}_\theta$, and for every
$Q\in\cM(X)$ and $\theta\in\cT$,
\begin{equation}\label{eq:envelope-bounds}
\E^{Q}[R_\theta]\le B,\qquad \E^{Q}[R^{L}_\theta]\le\varepsilon(L),\qquad
\sup_{Q,\theta}\E^{Q}\bigl[R_\theta\ind_{\{R_\theta>k\}}\bigr]
\ \le\ \frac{LB}{k}+\varepsilon(L).
\end{equation}
Consequently $\{R_\theta\}_{\theta\in\cT}$, and hence
$\{V_\theta\}_{\theta\in\cT}$, is uniformly integrable under every
$Q\in\cM(X)$, uniformly in $Q$.  If moreover
$K_p<\infty$ for some $p>1$, then $\sup_{Q,\theta}\E^{Q}[R_\theta^{p}]\le K_p$.
\end{pro}

\begin{proof}
$|V_\theta|\le R_\theta$: the upper bound is termwise, and
$V_\theta\ge\E^{Q}[\eta_\theta\mid\cF_\theta]=\eta_\theta\ge-R_\theta$.  The
inequality $R_\theta\le L+R^{L}_\theta$ follows from
$|\eta_\tau|\le L+|\eta_\tau|\ind_{\{|\eta_\tau|>L\}}$.

For the first two bounds in \eqref{eq:envelope-bounds}, fix $Q$ and $\theta$.
The family $\{\E^{Q'}[|\eta_\tau|\mid\cF_\theta]\}_{Q',\tau}$ is directed upward:
the maximum of two of its elements is again of the same form, obtained by
pasting the two measures and the two stopping times on the
$\cF_\theta$-measurable set where the first dominates the second, by
Lemma~\ref{lem:event-paste}.  So there
are pairs $(Q_n,\tau_n)$ whose conditional values $r_n$ increase to $R_\theta$;
a second application of monotone convergence gives the $L^{p}$ claim, since
$r_n^{p}\uparrow R_\theta^{p}$.
Let $\widehat Q_n:=Q\otimes_\theta Q_n\in\cM(X)$ be the pasting of
Definition~\ref{def:mstable}; it agrees with $Q$ on
$\cF_\theta$ and has the same conditional expectations after $\theta$ as $Q_n$,
so $r_n=\E^{\widehat Q_n}[|\eta_{\tau_n}|\mid\cF_\theta]$ and
\[
\E^{Q}[r_n]=\E^{\widehat Q_n}[r_n]
=\E^{\widehat Q_n}\bigl[|\eta_{\tau_n}|\bigr]\le B .
\]
Monotone convergence gives $\E^{Q}[R_\theta]\le B$, and the same argument with
$|\eta_\tau|\ind_{\{|\eta_\tau|>L\}}$ in place of $|\eta_\tau|$ gives
$\E^{Q}[R^{L}_\theta]\le\varepsilon(L)$.  For the $L^{p}$ statement, conditional
Jensen gives $r_n^{p}\le\E^{\widehat Q_n}[|\eta_{\tau_n}|^{p}\mid\cF_\theta]$ and
hence $\E^{Q}[r_n^{p}]\le K_p$.

For the third bound, $R_\theta\le L+R^{L}_\theta$ gives
\[
R_\theta\ind_{\{R_\theta>k\}}\ \le\ L\ind_{\{R_\theta>k\}}+R^{L}_\theta ,
\qquad
\E^{Q}\bigl[L\ind_{\{R_\theta>k\}}\bigr]\ \le\
\frac{L\,\E^{Q}[R_\theta]}{k}\ \le\ \frac{LB}{k},
\]
the second by Markov's inequality.
Given $\epsilon>0$, choose $L$ with $\varepsilon(L)<\epsilon/2$ and then $k$
with $LB/k<\epsilon/2$; this proves the uniform integrability, uniformly in
$Q$ and $\theta$.
\end{proof}

\begin{cor}[Envelope under the growth condition]
\label{cor:envelope-common-domination}
Under Assumption~\ref{assum:G}, $-M_\theta\le V_\theta\le M_\theta$ for every
$\theta\in\cT$, and $\{V_\theta\}_{\theta\in\cT}$ is of class~(D) under every
$Q\in\cM(X)$.
\end{cor}

\begin{proof}
For $\tau\in\cT_\theta$ and $Q\in\cM(X)$, optional sampling for the uniformly
integrable $Q$-martingale $M$ gives
$\E^{Q}[\eta_\tau\mid\cF_\theta]\le\E^{Q}[M_\tau\mid\cF_\theta]=M_\theta$, and
likewise $\ge-M_\theta$; take the essential suprema.  The family
$\{M_\theta\}_{\theta\in\cT}$ is uniformly $Q$-integrable, and a family
dominated in absolute value by a uniformly integrable one is uniformly
integrable.
\end{proof}

\subsection{\texorpdfstring{$\cT$}{T}-systems and aggregation}

\begin{defn}\label{def:Tsystem}
A family $\{G(\theta)\}_{\theta\in\cT}$ of random variables is a
\emph{$\cT$-system} if $G(\theta)$ is $\cF_\theta$-measurable for every
$\theta\in\cT$ and $G(\theta)=G(\theta')$ a.s.\ on $\{\theta=\theta'\}$, for all
$\theta,\theta'\in\cT$.  It is a \emph{strong $Q$-supermartingale system} if
moreover $G(\theta)\in L^1(Q)$ for all $\theta$ and
$G(\sigma)\ge\E^{Q}[G(\tau)\mid\cF_\sigma]$ $Q$-a.s.\ whenever
$\sigma\le\tau$.  An optional process $Z$ \emph{aggregates} the system if
$Z_\theta=G(\theta)$ a.s.\ for every $\theta\in\cT$.
\end{defn}

We use the classical aggregation theorem for supermartingale systems, which
requires no semicontinuity hypothesis whatsoever.  It is a special case of
\cite[Th\'eor\`eme~6 and Remarque~7(a)]{Dellacherie-Lenglart-1982}.  That
theorem aggregates a $\cT$-system which is \emph{s.c.s.\ \`a droite en
esp\'erance} --- right upper semicontinuous in expectation at stopping times,
that is, $\E^{Q}[G(\theta)]\ge\limsup_n\E^{Q}[G(\theta_n)]$ whenever
$\theta_n\downarrow\theta$ --- and which \emph{majore une martingale
\cadlag}, i.e.\ dominates one.  A supermartingale system satisfies both
hypotheses for free: it dominates the \cadlag{} $Q$-martingale
$\E^{Q}[G(T)\mid\cF_\cdot]$, and its supermartingale property gives the
expectation inequality along any decreasing sequence of stopping times.  This
is precisely the specialisation recorded in Remarque~7(a), ``toute
$\tau$-surmartingale s'agr\`ege en une surmartingale forte''.  Class~(D) is what makes
the aggregator a strong optional supermartingale \emph{of class~(D)}, hence
\ladlag{} with a Mertens decomposition under each $Q$, which is what
Section~\ref{sec:decomposition} consumes; $L^{1}$ at stopping times, already
part of Definition~\ref{def:Tsystem}, suffices for the aggregation itself,
which already delivers a strong optional supermartingale, hence \ladlag.  See
also \cite[Ch.~VI]{Dellacherie-Meyer-1980} and, for the
aggregation of value families, \cite[Appendix~B, Proposition~B.1]{Kobylanski2012}.

\begin{thm}[Aggregation of supermartingale
systems]\label{thm:aggregation}
Let $Q$ be a probability measure and let $\{G(\theta)\}_{\theta\in\cT}$ be a
strong $Q$-supermartingale system of class~(D), i.e.\ such that
$\{G(\theta)\}_{\theta\in\cT}$ is uniformly $Q$-integrable.  Then
there exists an optional \ladlag{} process $Z$, unique up to
indistinguishability, which aggregates the system.  Uniqueness holds by the
optional section theorem.  A uniformly bounded system is a special case.
\end{thm}

\begin{rem}\label{rem:aggregation-hyp}
We state Theorem~\ref{thm:aggregation} at class~(D) because that is what is
available for unbounded data: the system aggregated in
Theorem~\ref{thm:snell-aggreg} satisfies $|V_\theta|\le R_\theta$ with
$\{R_\theta\}$ uniformly integrable (Proposition~\ref{pro:envelope-UI}), or
$|V_\theta|\le M_\theta$ under Assumption~\ref{assum:G}
(Corollary~\ref{cor:envelope-common-domination}).  Aggregation is carried out
under the single reference measure $\Prob$; equivalence of the measures and the
stopping-time identities then give the same optional process under every other
$Q\in\cM(X)$.
\end{rem}

\subsection{The envelope is a uniform supermartingale}

\begin{thm}\label{thm:snell-supermart}
Under Assumption~\ref{assum:E} or Assumption~\ref{assum:G},
$\{V_\theta\}_{\theta\in\cT}$ is a $\cT$-system which is of class~(D) under every
$Q\in\cM(X)$, and for every such $Q$ it is a strong $Q$-supermartingale system.
Moreover $V_\theta\ge\eta_\theta$ a.s.\ for every $\theta\in\cT$, and
$V_T=\eta_T$.
\end{thm}

\begin{proof}
\emph{$\cT$-system.}  $V_\theta$ is $\cF_\theta$-measurable by construction.  If
$\theta,\theta'\in\cT$ and $B:=\{\theta=\theta'\}$, then $B\in\cF_\theta\cap
\cF_{\theta'}$ and, for $\tau\in\cT_\theta$, the time
$\tau^{B}:=\tau\ind_B+\theta'\ind_{B^{c}}$ belongs to $\cT_{\theta'}$ and
coincides with $\tau$ on $B$; since $\cF_\theta$ and $\cF_{\theta'}$ agree on
$B$, taking essential suprema gives $V_\theta=V_{\theta'}$ a.s.\ on $B$.  The
class~(D) property is Proposition~\ref{pro:envelope-UI} under
Assumption~\ref{assum:E} and Corollary~\ref{cor:envelope-common-domination}
under Assumption~\ref{assum:G}; in particular $V_\theta\in L^{1}(Q)$ for every
$Q\in\cM(X)$ and $\theta\in\cT$.  Next, $V_\theta\ge\eta_\theta$
follows by taking $\tau=\theta$ and $Q=\Prob$.  Finally $\cT_T=\{T\}$, so
$V_T=\eta_T$.

\emph{Supermartingale property.}  Fix $Q\in\cM(X)$ and $\sigma\le\tau$ in $\cT$.
By Lemma~\ref{lem:directed} applied at the stopping time $\tau$ to the family
$\bigl\{\E^{Q'}[\eta_{\tau'}\mid\cF_\tau]:\,Q'\in\cM(X),\,\tau'\in\cT_\tau
\bigr\}$ --- which is upward directed, the maximum of two of its elements being
again of the same form after pasting the two measures and the two stopping times
on the $\cF_\tau$-measurable set where the first dominates the second, by
Lemma~\ref{lem:event-paste} --- there
are sequences $(Q_n)\subset\cM(X)$ and $(\theta_n)$ with
$\theta_n\in\cT_\tau$ such that
\begin{equation}\label{eq:Xn}
\Xi_n:=\E^{Q_n}\bigl[\eta_{\theta_n}\mid\cF_\tau\bigr]\ \uparrow\ V_\tau
\qquad\text{a.s.}
\end{equation}
Fix $n$ and let $Q^{*}_n\in\cM(X)$ be the pasting of $Q$ and $Q_n$ at time
$\tau$ in the sense of Definition~\ref{def:mstable}; thus $Q^{*}_n=Q$ on
$\cF_\tau$ and $\E^{Q^{*}_n}[\,\cdot\mid\cF_\tau]
=\E^{Q_n}[\,\cdot\mid\cF_\tau]$, and $Q^*_n\in\cM(X)$ by
Proposition~\ref{pro:mstable}.  In particular
\begin{equation}\label{eq:paste1}
\E^{Q^{*}_n}\bigl[\eta_{\theta_n}\mid\cF_\tau\bigr]=\Xi_n .
\end{equation}
Since $\Xi_n$ is $\cF_\tau$-measurable and $Q^{*}_n=Q$ on $\cF_\tau\supset
\cF_\sigma$, we have $\E^{Q^{*}_n}[\ind_B\Xi_n]=\E^{Q}[\ind_B\Xi_n]$ for every
$B\in\cF_\sigma$, i.e.
\begin{equation}\label{eq:paste2}
\E^{Q^{*}_n}\bigl[\Xi_n\mid\cF_\sigma\bigr]
=\E^{Q}\bigl[\Xi_n\mid\cF_\sigma\bigr].
\end{equation}
Combining \eqref{eq:paste1}, \eqref{eq:paste2} and the tower property under
$Q^{*}_n$,
\[
\E^{Q}\bigl[\Xi_n\mid\cF_\sigma\bigr]
=\E^{Q^{*}_n}\bigl[\Xi_n\mid\cF_\sigma\bigr]
=\E^{Q^{*}_n}\Bigl[\E^{Q^{*}_n}\bigl[\eta_{\theta_n}\mid\cF_\tau\bigr]
\ \Big|\ \cF_\sigma\Bigr]
=\E^{Q^{*}_n}\bigl[\eta_{\theta_n}\mid\cF_\sigma\bigr]\le V_\sigma ,
\]
the last inequality holding by definition of $V_\sigma$, because
$Q^{*}_n\in\cM(X)$ and $\theta_n\in\cT_\tau\subset\cT_\sigma$.  Letting
$n\to\infty$ in \eqref{eq:Xn}, the $\Xi_n$ being increasing and signed, we apply
conditional monotone convergence to the non-negative sequence
$\Xi_n-\Xi_1$ and add back $\E^{Q}[\Xi_1\mid\cF_\sigma]$, which is finite
because $|\Xi_1|\le R_\tau$ under Assumption~\ref{assum:E} --- respectively
$|\Xi_1|\le M_\tau$ under Assumption~\ref{assum:G} --- and those bounds are
$Q$-integrable.  This gives
\[
\E^{Q}\bigl[V_\tau\mid\cF_\sigma\bigr]
=\lim_{n}\E^{Q}\bigl[\Xi_n\mid\cF_\sigma\bigr]\le V_\sigma .
\]
No limit is interchanged with an essential supremum over $\cM(X)$: every limit
above is taken under a single fixed measure.
\end{proof}

\begin{thm}\label{thm:snell-aggreg}
Under Assumption~\ref{assum:E} or Assumption~\ref{assum:G} there exists an
optional \ladlag{} process,
still denoted $V=(V_t)_{0\le t\le T}$ and called the \emph{generalised Snell
envelope of $\eta$}, unique up to indistinguishability, such that
\[
V_\theta=\esssup_{Q\in\cM(X)}\ \esssup_{\tau\in\cT_\theta}\
\E^{Q}\bigl[\eta_\tau\mid\cF_\theta\bigr]\qquad\text{a.s.\ for every }
\theta\in\cT .
\]
It satisfies $V\ge\eta$, $V_T=\eta_T$, and it is a strong
optional $Q$-supermartingale of class~(D) for every $Q\in\cM(X)$; under
Assumption~\ref{assum:E} also $|V|\le R$, where $R:=V(|\eta|)$ is the
generalised Snell envelope of $|\eta|$, the inequality $|V_\theta|\le R_\theta$
of Proposition~\ref{pro:envelope-UI} holding at every stopping time and hence
up to evanescence; and under Assumption~\ref{assum:G} also $|V|\le M$.  In
particular its right jumps are non-positive,
\begin{equation}\label{eq:rightjumpV}
\Dp V_t=V_{t+}-V_t\le0\qquad\text{for all }t\in[0,T)\text{ simultaneously,
a.s.}
\end{equation}
The quantifier matters: \eqref{eq:rightjumpV} is used below at stopping times
--- in the proofs of Lemma~\ref{lem:Y-max} and
Proposition~\ref{pro:C-skorokhod} --- and since $V$ is \ladlag{} an inequality
holding at each deterministic time separately would not license that.
\end{thm}

\begin{proof}
Theorem~\ref{thm:snell-supermart} and Theorem~\ref{thm:aggregation} applied
under $Q=\Prob$ produce an optional \ladlag{} aggregator $V$, unique up to
indistinguishability by the optional section theorem.  For any other
$Q\in\cM(X)$ the supermartingale inequalities of
Theorem~\ref{thm:snell-supermart} read $V_\sigma\ge\E^{Q}[V_\tau\mid
\cF_\sigma]$, i.e.\ $V$ is a strong optional $Q$-supermartingale; it is of
class~(D) under every $Q$ by Theorem~\ref{thm:snell-supermart}, hence by
Proposition~\ref{pro:envelope-UI} under Assumption~\ref{assum:E} and by
Corollary~\ref{cor:envelope-common-domination} under
Assumption~\ref{assum:G}.  For \eqref{eq:rightjumpV} we do not argue
along a sequence $t_n\downarrow t$, which would give the inequality at each
deterministic time separately and no more.  Instead we read it off the Mertens
decomposition.  Being a strong optional $\Prob$-supermartingale of class~(D),
$V$ admits the decomposition \eqref{eq:mertens} under $\Prob$, and
\eqref{eq:rightjump} is an identity between processes, valid for all $t\in[0,T)$
simultaneously outside a single null set:
$\Dp V_t=-\Dm C^{\Prob}_t$.  Since $C^{\Prob}$ is non-decreasing,
$\Dm C^{\Prob}\ge0$ everywhere, and \eqref{eq:rightjumpV} follows up to
evanescence.  In particular it may be evaluated at any stopping time.
\end{proof}

\begin{rem}
Property \eqref{eq:rightjumpV} is what makes the process $C$ of
Lemma~\ref{lem:aggreg-C} non-decreasing when applied to $V$: the right jumps of
a strong supermartingale are automatically downward.  Note that no
semicontinuity assumption on $\eta$ was needed to obtain
Theorem~\ref{thm:snell-aggreg}; in particular we do not need the notion of right
upper semicontinuity in expectation, nor a pasting theorem for
r.u.s.c.\ $\cT$-systems.  This is a genuine simplification with respect to the
usual route, and it is available precisely because we aggregate a
supermartingale system rather than an arbitrary one.
\end{rem}

\subsection{Minimality and decomposition}

\begin{thm}[Minimality]\label{thm:snell-minimal}
Under Assumption~\ref{assum:E} or Assumption~\ref{assum:G}, $V$ is the
smallest optional process which
dominates $\eta$ and is a strong optional $Q$-supermartingale for every
$Q\in\cM(X)$.  Precisely, if $U$ is optional, $U\ge\eta$ up to evanescence, and
$U$ is a strong optional $Q$-supermartingale for every $Q\in\cM(X)$, then
$U_\theta\ge V_\theta$ a.s.\ for every $\theta\in\cT$, hence $U\ge V$ up to
evanescence.
\end{thm}

\begin{proof}
Let $\theta\in\cT$, $Q\in\cM(X)$ and $\tau\in\cT_\theta$.  The strong
supermartingale property of $U$ under $Q$ and $U\ge\eta$ give
\[
U_\theta\ \ge\ \E^{Q}\bigl[U_\tau\mid\cF_\theta\bigr]
\ \ge\ \E^{Q}\bigl[\eta_\tau\mid\cF_\theta\bigr].
\]
Taking the essential supremum over $\tau\in\cT_\theta$ and over $Q\in\cM(X)$
yields $U_\theta\ge V_\theta$.  Since $V$ itself has the two stated properties
by Theorem~\ref{thm:snell-aggreg}, it is the smallest such process.  The
statement up to evanescence follows from the optional section theorem.
\end{proof}

Combining Theorem~\ref{thm:snell-aggreg} with Theorem~\ref{thm:uniform} gives
the structure of the envelope.

\begin{cor}\label{cor:snell-decomp}
Under Assumption~\ref{assum:E} or Assumption~\ref{assum:G} there exist a
predictable $X$-integrable process
$H$, an adapted right-continuous non-decreasing process $A$ with $A_0=0$, and
the purely discontinuous non-decreasing process
$C_t=-\sum_{s\le t,\,s<T}\Dp V_s$, none depending on $Q$, such that
\begin{equation}\label{eq:snell-decomp}
V_t=V_0+\int_0^tH_s\,dX_s-A_t-C_{t-},\qquad 0\le t\le T,\quad\text{a.s.},
\end{equation}
with $\int H\,dX$ a uniformly integrable $Q$-martingale and
$\E^{Q}[A_T+C_T]=V_0-\E^{Q}[\eta_T]$ for every $Q\in\cM(X)$.  If moreover
Assumption~\ref{assum:H1} holds, $A$ may be taken predictable and is then
unique.
\end{cor}

\begin{proof}
Immediate from Theorems~\ref{thm:snell-aggreg}, \ref{thm:uniform}
and~\ref{thm:predictable}: the envelope $V$ is a strong optional
supermartingale of class~(D) under every $Q\in\cM(X)$.  No bound on $V$ is
used: this is exactly where the removal of the lower bound from
Theorem~\ref{thm:uniform} does its work.  Under Assumption~\ref{assum:E} there
is in addition the explicit estimate
\[
\sup_{Q\in\cM(X)}\E^{Q}[A_T+C_T]\ \le\ 2B ,
\qquad\text{respectively}\ \le\ 2K_p^{1/p} ,
\]
the second whenever $K_p<\infty$ for some $p>1$.  Both follow from
$|V_0|\le B$ and $|\E^{Q}[\eta_T]|\le B$, respectively from the
$L^{p}$ bound, in Proposition~\ref{pro:envelope-UI}.
\end{proof}

The envelope of the absolute value now supplies a common traded bound, which is
how Assumption~\ref{assum:E} implies Assumption~\ref{assum:G}.

\begin{cor}[{\bf E} implies {\bf G}]\label{cor:ui-traded-domination}
Let Assumption~\ref{assum:E} hold and let $R$ be the generalised Snell envelope
of $|\eta|$, with decomposition
$R=R_0+\int_0^{\cdot}H^{R}\,dX-A^{R}-C^{R}_{-}$ given by
Corollary~\ref{cor:snell-decomp}.  Then
\[
M:=R_0+\int_0^{\cdot}H^{R}\,dX=R+A^{R}+C^{R}_{-}
\]
is a non-negative uniformly integrable $Q$-martingale for every $Q\in\cM(X)$,
of the traded form of Assumption~\ref{assum:G}, and $|\eta|\le R\le M$.  Thus
$\eta$ satisfies Assumption~\ref{assum:G}.
\end{cor}

\begin{proof}
Assumption~\ref{assum:E} is invariant under $\eta\mapsto|\eta|$, so
Corollary~\ref{cor:snell-decomp} applies to $|\eta|$ and makes
$\int H^{R}\,dX$ a uniformly integrable $Q$-martingale for every $Q$; hence so
is $M$, and $M$ is traded by construction.  As $A^{R},C^{R}_{-}\ge0$ we get
$M\ge R\ge|\eta|\ge0$.
\end{proof}

\begin{eg}[{\bf G} does not imply {\bf E}]\label{eg:G-not-UI}
The implication in Corollary~\ref{cor:ui-traded-domination} is strict.  Let
$W=(W^{1},W^{2})$ be a two-dimensional Brownian motion, $\F$ its augmented
filtration, $T=2$, $\nu:=e^{W^{2}_1}$, and
\[
X_t:=1\ \ (t\le1),\qquad
X_t:=\exp\Bigl(\nu\bigl(W^{1}_t-W^{1}_1\bigr)-\tfrac12\nu^{2}(t-1)\Bigr)
\ \ (1<t\le2).
\]
For $Q\in\cM(X)$ the Girsanov drift on $W^{1}$ must vanish after time $1$, so
conditionally on $\cF_1$ the increments $W^{1}_t-W^{1}_1$ remain standard under
$Q$ and $X$ is a true uniformly integrable $Q$-martingale.  Hence $\eta:=X$
satisfies Assumption~\ref{assum:G}, with $M=X$.

It fails \eqref{eq:UI}.  Conditionally on $\nu$, $\log X_2$ is Gaussian with
mean $-\nu^{2}/2$ and variance $\nu^{2}$, so for $L>1$ the Cameron--Martin shift
gives
\[
\E\bigl[X_2\ind_{\{X_2>L\}}\bigm|\nu\bigr]
=\Phi\Bigl(\frac\nu2-\frac{\log L}{\nu}\Bigr)\ \longrightarrow\ 1
\qquad(\nu\to\infty),
\]
where $\Phi$ is the standard normal distribution function.  Changes of measure
with $\cF_1$-measurable densities leave $X$ a local martingale and can
concentrate the law of $\nu$ on $\{\nu\ge n\}$; along such a sequence
$\E^{Q}[X_2\ind_{\{X_2>L\}}]\to1$.  Hence $\varepsilon(L)=1$ for every $L>1$ and
\eqref{eq:UI} fails, although Assumption~\ref{assum:G} holds.  The mechanism is
that the \emph{size} of the traded bound is not controlled uniformly over the
family, only its martingale property.
\end{eg}

\begin{rem}[The increasing process is not carried by the contact set]
\label{rem:no-skorokhod}
For a classical Snell envelope whose reward has finite left limits, the
predictable increasing part obeys the usual contact conditions, using left
contact at its jumps (Remark~\ref{rem:complete-case}).  For a completely
irregular reward those conditions use the left upper envelope instead.
Even within the class of rewards with left limits, the usual contact
conditions fail for some rewards as soon as $\cM(X)$ is not a singleton: for a
European-type reward, Proposition~\ref{pro:no-skorokhod} below shows that the
Skorokhod conditions on the total regulator hold only if the terminal claim is
attainable by a bounded wealth process.  The increasing part of the minimal
value superposes the reflection of $V$ on $\eta$ and the premium for
non-attainability, and the latter may increase on $\{V>\eta\}$.  This is why
the reflected BSDE of Section~\ref{sec:bsde} is formulated through minimality;
Corollary~\ref{cor:dyn-minimality} gives an equivalent condition at each
stopping time, and Remark~\ref{rem:complete-case} recovers the Skorokhod
conditions when $\cM(X)=\{\Prob\}$.
\end{rem}

\section{Reflected BSDEs under the family \texorpdfstring{$\cM(X)$}{M(X)}}\label{sec:bsde}

\subsection{An incomplete It\^o market}\label{subsec:market}

From now on we place ourselves in the setting of Example~\ref{eg:H1-brownian},
which we restate for convenience.  Let $W$ be an $n$-dimensional Brownian motion
on $(\Omega,\cF,\Prob)$, let $\F$ be the usual augmentation of its natural
filtration, and let
\begin{equation}\label{eq:market}
X_t=\int_0^t\sigma_s\,dW_s,\qquad 0\le t\le T,
\end{equation}
where $\sigma$ is a predictable $d\times n$ matrix-valued process with
$\int_0^T\|\sigma_s\|^2ds<\infty$ a.s.; $X$ is then continuous, hence locally
bounded, as Section~\ref{sec:prelim} requires.  No bound and no local
boundedness of $\sigma$ is assumed.  The components of $X$ are the discounted
prices of $d$ traded assets driven by $n$ sources of risk; some assets may be
redundant, and $d$ may exceed $n$.  The market is complete if and only if
$\mathrm{rank}\,\sigma_s=n$ $ds\otimes d\Prob$-a.e., which requires $d\ge n$,
and \emph{incomplete}, with $\cM(X)$ infinite,
as soon as $\mathrm{rank}\,\sigma_s<n$ on a set of positive
$ds\otimes d\Prob$-measure.  By
Example~\ref{eg:H1-brownian}, Assumption~\ref{assum:H1} holds in this market, so
that the predictable form of the decomposition
(Theorem~\ref{thm:predictable}) is available throughout this section.

For a predictable $X$-integrable $\R^{d}$-valued process $H$, that is, one
with $\int_0^T|\sigma_s\trans H_s|^2ds<\infty$ a.s., we write
\begin{equation}\label{eq:ZH}
Z_t:=\sigma_t\trans H_t\in\R^{n},\qquad\text{so that}\qquad
\int_0^{t}H_s\,dX_s=\int_0^{t}Z_s\trans\,dW_s ,
\end{equation}
and we note that $Z$ is constrained to take its values in the predictable family
of subspaces $\cK_t:=\mathrm{Im}\,\sigma_t\trans\subset\R^{n}$.  For
$Q\in\cM(X)$ with $\lambda^{Q}=\int\ell\trans dW$ and $\sigma\ell=0$, Girsanov's
theorem makes $W^{Q}:=W-\int\ell\,ds$ a $Q$-Brownian motion, and
$Z\trans\ell=H\trans\sigma\ell=0$ gives
$\int H\,dX=\int Z\trans dW=\int Z\trans dW^{Q}$: the same integrand represents
the gain process under every $Q\in\cM(X)$.  This constraint
is exactly the market incompleteness, and it is the reason why the increasing
process of the decomposition cannot be dispensed with.  No bound on $\sigma$
is used in the rest of this section: the estimates are stated for $Z$ and for
the gain $\int Z\trans dW$, not for $\sigma$ or $H$ separately.

We use the following norms, in which the supremum over measures is a genuine
supremum of real numbers:
\[
\begin{gathered}
\|\varphi\|_{\Lbb^2}^2:=\sup_{Q\in\cM(X)}\E^{Q}\bigl[|\varphi|^2\bigr],
\qquad
\|\varphi\|_{\Hbb^2}^2:=\sup_{Q\in\cM(X)}
\E^{Q}\Bigl[\int_0^T|\varphi_s|^2ds\Bigr],\\
\|\varphi\|_{\Sbb^2}^2:=\sup_{Q\in\cM(X)}
\E^{Q}\Bigl[\sup_{0\le t\le T}|\varphi_t|^2\Bigr],
\end{gathered}
\]
and we write $\Sbb^{\infty}$ for the space of bounded optional processes,
identified up to evanescence, with
\[
\|\varphi\|_{\infty}:=\Bigl\|\sup_{0\le t\le T}|\varphi_t|\Bigr\|_{L^{\infty}(\Prob)}
=\sup_{\tau\in\cT}\|\varphi_\tau\|_{L^{\infty}(\Prob)},
\]
the equality following from the optional section theorem.  The supremum is
over stopping times and not over deterministic times: for an atomless stopping
time $\tau$ the optional process $\ind_{\{t=\tau\}}$ vanishes a.s.\ at every
deterministic $t$ but has norm one.  With this norm $\Sbb^{\infty}$ is a Banach
space, and $\Sbb^{\infty}\subset\Sbb^2\subset\Hbb^2$.

\begin{assum}[Data]\label{assum:data}
\begin{itemize}\itemsep2pt
\item[{\bf (A1)}] The obstacle $\xi=(\xi_t)_{0\le t\le T}$ is optional and
bounded: $|\xi_t|\le c_\xi$ for all $(t,\omega)$.  No left or right limits,
and no semicontinuity, are assumed.
\item[{\bf (A2)}] The driver $f:[0,T]\times\Omega\times\R\to\R$ is
$\cP\otimes\cB(\R)$-measurable, $L$-Lipschitz in $y$ uniformly in $(t,\omega)$,
and $|f(t,\omega,0)|\le c_f$ for all $(t,\omega)$.
\end{itemize}
\end{assum}

\begin{rem}[The market is continuous; the obstacle is not]\label{rem:where-ladlag}
It is worth saying at once where the irregularity of this section comes from,
since the market of \eqref{eq:market} supplies none of it.  There $X$ is
continuous, $\F$ is the usual augmentation of a Brownian filtration, so every
$\Prob$-martingale is continuous and the optional and predictable
$\sigma$-fields coincide up to evanescence
\cite[Theorem~4.11 and Example~4.12]{Nikeghbali2006}; in particular the wealth process
$\int H\,dX=\int Z\trans\,dW$ of \eqref{eq:ZH} is continuous.  None of this
regularises the value process.

The irregularity is contractual, not dynamic.  {\bf (A1)} allows $\xi$ to be an
arbitrary bounded optional process, and the generalised Snell
envelope of an irregular obstacle can remain irregular however smooth the
filtration and the traded assets may be: nothing in the dynamics forces $Y$ to
be right-continuous.  Example~\ref{eg:windowed} is the case to keep in
mind: if exercise is permitted only on a finite union $D$ of closed intervals
and isolated dates, the effective payoff is optional and \ladlag{} but not
right-continuous, and the value process can have right jumps, which can occur
only at right endpoints of windows and at isolated dates --- in a Brownian
filtration, with a continuous $X$.

The solution $Y$ is nevertheless \ladlag: it is an optional Snell envelope
minus a continuous time integral.  This regularity is a conclusion about $Y$,
not an assumption on $\xi$.  Finite left limits of $\xi$ are imposed locally
only when writing the classical conditions involving $\xi_{-}$, in
Remark~\ref{rem:why-minimal}, Corollary~\ref{cor:dichotomy},
Remark~\ref{rem:complete-case} and Problem~P3.
\end{rem}

\begin{rem}
The driver does not depend on the control variable $Z$.  The standard route to
a $Z$-dependent driver passes through Gal'chouk--Lenglart estimates in which
the reflection terms are signed by means of contact conditions; here the
canonical $C$ satisfies its contact condition
(Proposition~\ref{pro:C-skorokhod}) but the total predictable $A$ need not, and
we do not have the required difference estimate for it (Problem~P2).
Boundedness of the data is assumed until Section~\ref{subsec:weighted}, which
replaces it by domination by a traded martingale, measuring the solution in a
norm weighted by that martingale.  No maximal inequality for the sublinear
expectation $\esssup_{Q}\E^{Q}[\,\cdot\mid\cF_t]$ is needed (Problem~P1), but
the square-integrable estimates of Theorem~\ref{thm:apriori} are lost
(Example~\ref{eg:weighted-not-L2}).
\end{rem}

\subsection{Supersolutions and the definition of a solution}

\begin{defn}\label{def:rbsde}
A quadruple $(Y,H,A,C)$ is a \emph{supersolution} of the reflected BSDE with
data $(f,\xi)$ if:
\begin{enumerate}\itemsep2pt
\item[(i)] $Y\in\Sbb^{\infty}$ is optional and \ladlag; $H$ is predictable and
$X$-integrable, identified modulo equality of $\int H\,dX$; $A$ is adapted,
right-continuous, non-decreasing with $A_0=0$; $C$ is adapted,
right-continuous, purely discontinuous, non-decreasing with $C_{0-}=0$ and
$C_T=C_{T-}$; and $\E^{Q}[A_T+C_T]<\infty$ for every $Q\in\cM(X)$;
\item[(ii)] for every $\tau\in\cT$,
\begin{equation}\label{eq:rbsde}
Y_\tau=\xi_T+\int_\tau^T f(s,Y_s)\,ds-\int_\tau^T H_s\,dX_s
+A_T-A_\tau+C_{T-}-C_{\tau-}\qquad\text{a.s.};
\end{equation}
\item[(iii)] $Y_t\ge\xi_t$ for every $t\in[0,T]$, a.s., equivalently (by the
optional section theorem) $Y_\tau\ge\xi_\tau$ a.s.\ for every $\tau\in\cT$.
\end{enumerate}
A supersolution $(Y,H,A,C)$ is a \emph{solution} if it is minimal in the first
component: $Y_t\le Y'_t$ for all $t$, a.s., for every supersolution
$(Y',H',A',C')$.
\end{defn}

Besides Remark~\ref{rem:where-ladlag}, the predictability of $A$ in
Definition~\ref{def:rbsde} can be seen directly from its compensator.  Indeed, $\E^{\Prob}[A_T]<\infty$ by (i), so
let $A^p$ be its predictable compensator, with $A^p_0=0$.  Then $A-A^p$ is a
$\Prob$-martingale and is continuous in the Brownian filtration.  It also has
finite variation, so it vanishes identically.  Hence $A=A^p$ is predictable;
no additional assumption on a supersolution is needed.

Without the normalisation $C_T=C_{T-}$ uniqueness would fail trivially: with $f=0$,
$\xi=0$ and $Y=H=A=0$, both $C=0$ and $C=b\ind_{\{T\}}$, $b>0$, satisfy
\eqref{eq:rbsde}, which only sees $C_{T-}$.

\begin{rem}[Why minimality replaces the Skorokhod conditions]
\label{rem:why-minimal}
For this remark assume additionally that $\xi$ has finite left limits on
$(0,T]$, outside one null set; no right limits are needed.  In the classical
theory the increasing process is then pinned down by the
requirements $\int_0^T\ind_{\{Y_t>\xi_t\}}\,dA^{c}_t=0$,
$(Y_{\tau-}-\xi_{\tau-})\Dm A^{d}_\tau=0$ for predictable $\tau$, and
$(Y_\tau-\xi_\tau)\Dm C_\tau=0$.  By Proposition~\ref{pro:no-skorokhod} the
first two cannot be imposed on the single regulator of \eqref{eq:rbsde} for all
data once the market is incomplete, because that regulator also carries the
premium for claims that are not attainable by bounded wealth.  A formulation
carrying two increasing processes is not subject to this, and
Remark~\ref{rem:pengxu} discusses one.  For \eqref{eq:rbsde} as written we
select the solution by minimality, as in the constrained-BSDE literature; it
reduces to the Skorokhod conditions when the market is complete
(Remark~\ref{rem:complete-case}).  For arbitrary optional obstacles, the
classical conditions on $A$ use the left upper envelope instead of
$\xi_{-}$; the condition on $C$ still uses $\xi$ itself
\cite{Grigorova2020}.  No left-limit assumption enters
Definition~\ref{def:rbsde} or the existence, comparison and dynamic-minimality
results below.
\end{rem}

The following proposition makes the obstruction precise.  We call a bounded
$\cF_T$-measurable random variable $\zeta$ \emph{boundedly attainable} if there
are $x\in\R$ and a predictable $X$-integrable $H$ such that $x+\int_0^{\cdot}H\,dX$
is bounded on $[0,T]$ and $\zeta=x+\int_0^TH_s\,dX_s$.  Such a wealth process
is a bounded martingale under every $Q\in\cM(X)$, so
$\E^{Q}[\zeta]=x$ for all $Q$.  This is stronger than representability by a
wealth process that is only bounded below, which is the admissibility notion
of Definition~\ref{def:admissible}: in a two-Brownian market trading only
$W^{1}$, the indicator of $\{W^{2}_{1/2}>0\}$ has different expectations under
different elements of $\cM(X)$, hence is not boundedly attainable, yet (for
$T>1/2$) it is the terminal value of a non-negative wealth process started at
$1$ that, after time $1/2$ and on $\{W^{2}_{1/2}\le0\}$, runs the gain
$\int_{1/2}^{t}(T-s)^{-1}dW^{1}_s$ until it first reaches $-1$, which happens
before $T$ almost surely.

\begin{pro}[Skorokhod conditions and bounded attainability]\label{pro:no-skorokhod}
Let $\zeta$ be $\cF_T$-measurable and bounded, let $c>\|\zeta\|_{\infty}$, and
consider the European-type obstacle
\[
\xi_t:=-c\ \ (0\le t<T),\qquad \xi_T:=\zeta ,
\]
which is \cadlag{} and satisfies {\bf (A1)}, with $f\equiv0$.  A supersolution of
\eqref{eq:rbsde} satisfying the three Skorokhod conditions of
Remark~\ref{rem:why-minimal} exists if and only if $\zeta$ is boundedly
attainable.  Any such supersolution has $A\equiv0$, $C\equiv0$ and
$Y=Y_0+\int_0^{\cdot}H\,dX$ bounded with $Y_T=\zeta$; in particular
\begin{equation}\label{eq:unique-price}
\E^{Q}[\zeta]=Y_0\qquad\text{for every }Q\in\cM(X).
\end{equation}
Moreover the solution of Theorem~\ref{thm:existence} for these data has
$C\equiv0$, and it has $A\equiv0$ if and only if $\zeta$ is boundedly
attainable.
\end{pro}

\begin{proof}
The proof uses results stated later: Lemma~\ref{lem:super-is-supermart} and
the differential form \eqref{eq:rbsde-diff} of the equation for the first two
assertions, and Theorem~\ref{thm:existence} (representation and minimality)
and Proposition~\ref{pro:C-skorokhod} for the last one.  None of them depends
on the present subsection.  What is taken from
Lemma~\ref{lem:super-is-supermart} is only that the first component of a
supersolution is a strong optional $Q$-supermartingale for every
$Q\in\cM(X)$.

Let $(Y,H,A,C)$ be a supersolution.  Then $Y$ is a strong optional
$\Prob$-supermartingale with $Y_T=\zeta$, so $Y_\tau\ge\E^{\Prob}[\zeta\mid
\cF_\tau]\ge-\|\zeta\|_{\infty}$ for every $\tau\in\cT$; by the optional section
theorem $Y\ge-\|\zeta\|_{\infty}$ up to evanescence, and the left limits obey
the same bound.  Since $\xi=\xi_-=-c$ on $[0,T)$ and $\xi_{T-}=-c$, we get
$Y_t>\xi_t$ for $t<T$ and $Y_{t-}>\xi_{t-}$ for $0<t\le T$.

Assume the Skorokhod conditions.  The first forces $A^{c}\equiv0$.  The second
forces $\Dm A_\tau=0$ at every predictable $\tau\in(0,T]$, including $\tau=T$,
and $A$ is predictable by the compensator argument following
Definition~\ref{def:rbsde}; hence
$A^{d}\equiv0$.  The third forces $\Dm C_\tau=0$ for $\tau<T$, and
$C_T=C_{T-}$ removes the remaining atom; hence $C\equiv0$.  Then
\eqref{eq:rbsde-diff} reads $Y=Y_0+\int_0^{\cdot}H\,dX$, which is bounded since
$Y\in\Sbb^{\infty}$, and $Y_T=\zeta$; so $\zeta$ is boundedly attainable and
\eqref{eq:unique-price} holds, a bounded $Q$-local martingale being a
$Q$-martingale.

Conversely, if $\zeta=x+\int_0^TH\,dX$ with $x+\int_0^{\cdot}H\,dX$ bounded,
this wealth process is a bounded $\Prob$-martingale, equal to
$\E^{\Prob}[\zeta\mid\cF_\cdot]\ge-\|\zeta\|_\infty>-c$.  With $A=C=0$ it is a
supersolution, and it satisfies the three conditions because $A$ and $C$ do not
charge anything.

For the solution, $C\equiv0$ follows from Proposition~\ref{pro:C-skorokhod},
since the gap $Y-\xi$ is strictly positive before $T$ and $C$ has no atom at
$T$.  If $\zeta$ is boundedly attainable, the wealth process of the previous
paragraph is a supersolution and, being a martingale under every $Q$, equals
$\esssup_Q\E^{Q}[\zeta\mid\cF_\tau]\le Y_\tau$ at every $\tau$ by the
representation \eqref{eq:representation}; minimality then identifies it with
$Y$, so $A\equiv0$.  Conversely, if $A\equiv0$, then $Y=Y_0+\int_0^{\cdot}H\,dX$
is bounded and ends at $\zeta$.
\end{proof}

\begin{cor}[Dichotomy]\label{cor:dichotomy}
Exactly one of the following holds.
\begin{enumerate}\itemsep2pt
\item[(i)] $\cM(X)=\{\Prob\}$.  Then for all data satisfying
Assumption~\ref{assum:data} whose obstacle has finite left limits on $(0,T]$,
the solution of \eqref{eq:rbsde} satisfies the
Skorokhod conditions of Remark~\ref{rem:why-minimal}
(Remark~\ref{rem:complete-case}).
\item[(ii)] $\cM(X)$ has at least two elements.  Then there exist data
satisfying Assumption~\ref{assum:data}, even with a \cadlag{} obstacle, for
which no supersolution of
\eqref{eq:rbsde} satisfies the Skorokhod conditions.
\end{enumerate}
\end{cor}

\begin{proof}
Suppose $\cM(X)$ contains $Q_1\neq Q_2$, choose $B\in\cF_T$ with
$Q_1(B)\neq Q_2(B)$ and put $\zeta:=\ind_B$.  If a supersolution for the data
of Proposition~\ref{pro:no-skorokhod} built on $\zeta$ satisfied the Skorokhod
conditions, \eqref{eq:unique-price} would give $Q_1(B)=Y_0=Q_2(B)$.
\end{proof}

The quantifier in (ii) matters: in an incomplete market many data admit
Skorokhod solutions --- the zero obstacle, for instance --- and what fails is
the classical characterisation \emph{for all data with left limits}.  The
existence and minimality theorems themselves require only optionality of the
obstacle.  Since $\cM(X)$ is convex,
case (ii) means that $\cM(X)$ is infinite.

\begin{rem}[Constrained reflected BSDEs]\label{rem:pengxu}
Proposition~\ref{pro:no-skorokhod} concerns a \emph{single} regulator carrying
both the reflection and the premium for non-attainability.  The
constrained-BSDE literature, to which the minimality formulation of
Definition~\ref{def:rbsde} belongs, avoids it by carrying two increasing
processes.  Cvitani\'c, Karatzas and Soner \cite{CvitanicKaratzasSoner1998}
select the solution of a constrained backward equation by minimality and treat a
reflecting barrier in their Section~6; Peng \cite{Peng1999} obtains smallest
constrained $g$-supersolutions from his monotonic limit theorem.  Closest to us,
Peng and Xu \cite{PengXu2010} study reflected BSDEs with a constraint
$z_t\in\Gamma_t$ on the integrand --- for us $\Gamma_t=\operatorname{Im}
\sigma_t^{\top}$ --- whose solution is a quadruple $(y,z,A,\bar A)$: $A$ keeps
$z$ in the constraint and $\bar A$ reflects $y$ and satisfies the generalised
Skorokhod condition $\int_0^T(y_{s-}-L^{*}_{s-})\,d\bar A_s=0$ for every
\cadlag{} $L^{*}$ with $y\ge L^{*}\ge L$.  Their solution is the smallest
constrained supersolution, so minimality characterises it as it does here.

Two differences of setting remain.  Their value process is right-continuous
and the obstacle is dominated $d\Prob\otimes dt$-a.e., which for a
right-continuous obstacle is equivalent on $[0,T)$ to domination at every time
but for an irregular obstacle is strictly weaker than
$Y_\tau\ge\xi_\tau$ for all $\tau\in\cT$; the latter is what an American claim
exercised at a stopping time requires.  The Snell-envelope construction makes
$Y$ \ladlag{} and the process $C$ records its right jumps; $\xi$ itself need
not have either left or right limits.  The left limits in the Peng--Xu
condition above belong to their \cadlag{} comparison process $L^{*}$, not to
our general obstacle.  Our aggregated formulation carries a single regulator, which
yields the gross hedge $V+A+C_-$ of Theorem~\ref{thm:robust-american} directly,
at the price quantified by Proposition~\ref{pro:no-skorokhod}.
\end{rem}

We first record that supersolutions are uniform supermartingales; this is the
link with Section~\ref{sec:snell}.

\begin{lem}\label{lem:super-is-supermart}
Let $(Y,H,A,C)$ be a supersolution and set
$\Gamma_t:=\int_0^t f(s,Y_s)\,ds$.  Then, for every $Q\in\cM(X)$,
\begin{enumerate}\itemsep2pt
\item[(i)] $\int_0^{\cdot}H\,dX$ is a uniformly integrable $Q$-martingale;
\item[(ii)] $Y+\Gamma$ is a strong optional $Q$-supermartingale of class~(D)
which dominates the bounded optional process $\xi+\Gamma$.
\end{enumerate}
\end{lem}

\begin{proof}
Apply \eqref{eq:rbsde} at an arbitrary $\tau\in\cT$ and at $\tau=0$, which is
in $\cT$, and subtract.  Since $A_0=0$ and $C_{0-}=0$ by
Definition~\ref{def:rbsde}(i), and $\Gamma_0=0$, this gives, for each
$\tau\in\cT$,
\[
Y_\tau+\Gamma_\tau=Y_0+\int_0^\tau H_s\,dX_s-A_\tau-C_{\tau-}
\qquad\text{a.s.}
\]
(Subtracting instead the instance at $\tau=T$ would return \eqref{eq:rbsde}
itself, that instance reducing to $Y_T=\xi_T$.)  Both sides are optional
processes: $Y$ is optional by hypothesis, $A$ is adapted and right-continuous
hence optional, $C_-$ is adapted and left-continuous hence predictable,
$\Gamma$ is continuous and adapted, and $\int H\,dX$ is \cadlag{} and adapted.
They agree almost surely at every stopping time, so by the optional section
theorem they are indistinguishable:
\begin{equation}\label{eq:rbsde-diff}
Y_t+\Gamma_t=Y_0+\int_0^tH_s\,dX_s-A_t-C_{t-},\qquad 0\le t\le T,
\end{equation}
outside an evanescent set.  Agreement at deterministic times would not have
sufficed: for a process that is not right-continuous it does not imply
indistinguishability.  What the section theorem asks for --- agreement at
every stopping time --- is exactly what \eqref{eq:rbsde} supplies.
Write $N:=\int_0^{\cdot}H\,dX$ and
$\kappa:=2\|Y\|_{\infty}+2T(c_f+L\|Y\|_{\infty})$, so that
$|\Gamma_t|\le T(c_f+L\|Y\|_{\infty})$.  Since $A,C\ge0$, the process $N$ is
bounded below by $-\kappa$; being a stochastic integral with respect to the
locally bounded $Q$-local martingale $X$ and bounded below, it is a $Q$-local
martingale by the Ansel--Stricker lemma \cite{AnselStricker1994}.

For (i), observe that \eqref{eq:rbsde-diff} gives
$N_t=Y_t+\Gamma_t-Y_0+A_t+C_{t-}$, whence, $A$ and $C$ being
non-decreasing and non-negative,
\[
\sup_{0\le t\le T}|N_t|\ \le\ \kappa+A_T+C_T\ =:\ G .
\]
By Definition~\ref{def:rbsde}(i), $G\in L^{1}(Q)$.  Thus the family
$\{N_\tau:\tau\in\cT\}$ is dominated by an integrable random variable and is in
particular uniformly integrable, i.e.\ $N$ is a $Q$-local martingale of
class~(D); such a process is a uniformly integrable $Q$-martingale.  Indeed, if
$(\rho_n)$ localises $N$, then $N_{t\wedge\rho_n}\to N_t$ almost surely and, by
uniform integrability, in $L^1(Q)$, so that
$\E^{Q}[N_t\mid\cF_s]=\lim_n\E^{Q}[N_{t\wedge\rho_n}\mid\cF_s]
=\lim_n N_{s\wedge\rho_n}=N_s$.

For (ii), $N$ being in particular a $Q$-martingale, for
$\sigma\le\tau$ in $\cT$,
\[
\E^{Q}\bigl[Y_\tau+\Gamma_\tau\mid\cF_\sigma\bigr]
=Y_\sigma+\Gamma_\sigma
-\E^{Q}\bigl[A_\tau-A_\sigma+C_{\tau-}-C_{\sigma-}\mid\cF_\sigma\bigr]
\le Y_\sigma+\Gamma_\sigma .
\]
Class~(D) follows from boundedness, and $Y+\Gamma\ge\xi+\Gamma$ from (iii).
\end{proof}

\subsection{Reduction to a monotone driver}

\begin{lem}\label{lem:monotone-driver}
Define, for $(t,\omega,y)\in[0,T]\times\Omega\times\R$,
\[
\bar\xi_t:=e^{-Lt}\xi_t,\qquad
\bar f(t,\omega,y):=e^{-Lt}f\bigl(t,\omega,e^{Lt}y\bigr)+Ly .
\]
Then $\bar\xi$ and $\bar f$ satisfy Assumption~\ref{assum:data} with constants
$c_\xi$, $2L$, $c_f$, and $y\mapsto\bar f(t,\omega,y)$ is non-decreasing.
Moreover $(Y,H,A,C)$ is a supersolution (resp.\ a solution) for the data
$(f,\xi)$ if and only if $(\bar Y,\bar H,\bar A,\bar C)$ is a supersolution
(resp.\ a solution) for the data $(\bar f,\bar\xi)$, where
\[
\bar Y_t:=e^{-Lt}Y_t,\quad \bar H_t:=e^{-Lt}H_t,\quad
\bar A_t:=\int_{(0,t]}e^{-Ls}dA_s,\quad
\bar C_t:=\int_{[0,t]}e^{-Ls}dC_s .
\]
\end{lem}

\begin{proof}
Monotonicity: for $y\ge y'$,
$\bar f(t,y)-\bar f(t,y')\ge-e^{-Lt}L\,e^{Lt}(y-y')+L(y-y')=0$ by the Lipschitz
property.  The Lipschitz constant of $\bar f$ is at most $2L$ and
$|\bar f(t,0)|=e^{-Lt}|f(t,0)|\le c_f$.  For the equivalence, note that
$\bar A$ and $\bar C$ are non-decreasing, $\bar C$ is purely discontinuous, and
that \eqref{eq:rbsde-diff} transforms, by the product rule applied to the
\cadlag{} part and by $\Dp\bar Y_t=e^{-Lt}\Dp Y_t$ for the right jumps, into
\[
\bar Y_t+\int_0^t\bar f(s,\bar Y_s)\,ds
=\bar Y_0+\int_0^t\bar H_s\,dX_s-\bar A_t-\bar C_{t-} ,
\]
which is \eqref{eq:rbsde-diff} for the barred data.  The obstacle condition and
the ordering of the first components are preserved because $e^{-Lt}>0$.
\end{proof}

Monotonicity of the driver,
\begin{equation}\label{eq:monotone}
y\longmapsto f(t,\omega,y)\ \text{is non-decreasing,}
\end{equation}
is used only for minimality and comparison.  We prove the existence theorem
first under \eqref{eq:monotone} and then transfer it, representation included,
to a general Lipschitz driver through Lemma~\ref{lem:monotone-driver}.

\subsection{The zero-driver case}

\begin{thm}\label{thm:zero-driver}
Let $\zeta$ be a bounded optional process and consider the data
$(f,\xi)=(0,\zeta)$.  Then the reflected BSDE
\eqref{eq:rbsde} has a unique solution $(Y,H,A,C)$ in the sense of
Definition~\ref{def:rbsde}, and
\[
Y=V(\zeta),
\]
the generalised Snell envelope of $\zeta$ constructed in
Theorem~\ref{thm:snell-aggreg}.  The processes $H$, $A$, $C$ are those of
Corollary~\ref{cor:snell-decomp}; $A$ is predictable,
$C_t=-\sum_{s\le t,\,s<T}\Dp Y_s$, and $A$ and $\int H\,dX$ are unique.
\end{thm}

\begin{proof}
\emph{$V(\zeta)$ is a supersolution.}  By Theorem~\ref{thm:snell-aggreg},
$V:=V(\zeta)$ is optional, \ladlag, bounded by $c_\zeta$, dominates $\zeta$ and
satisfies $V_T=\zeta_T$.  By Corollary~\ref{cor:snell-decomp} and
Assumption~\ref{assum:H1} (valid here by Example~\ref{eg:H1-brownian}) there are
$H$ predictable and $X$-integrable, $A$ predictable right-continuous
non-decreasing with $A_0=0$, and $C_t=-\sum_{s\le t,\,s<T}\Dp V_s$ purely
discontinuous non-decreasing with $C_{0-}=0$ and $C_T=C_{T-}$, such that
$V_t=V_0+\int_0^tH\,dX-A_t-C_{t-}$ with $\E^{Q}[A_T+C_T]<\infty$ for every $Q$.
This identity is one between processes, valid up to indistinguishability by
Theorem~\ref{thm:uniform} and Lemma~\ref{lem:aggreg-C}, so it may be evaluated
directly at any $\tau\in\cT$; doing so at $\tau$ and at $T$ and subtracting
gives \eqref{eq:rbsde} with $f=0$ and $\xi=\zeta$ for every $\tau\in\cT$.  We
do not obtain it at deterministic times and then invoke the optional section
theorem: for a \ladlag{} process that upgrade is exactly what is unavailable,
as the discussion before Lemma~\ref{lem:Y-max} records.

\emph{Minimality.}  Let $(Y',H',A',C')$ be any supersolution.  By
Lemma~\ref{lem:super-is-supermart} with $f=0$, $Y'$ is a strong optional
$Q$-supermartingale for every $Q\in\cM(X)$ and $Y'\ge\zeta$.  By
Theorem~\ref{thm:snell-minimal}, $Y'\ge V$.  Hence $V$ is minimal.

\emph{Uniqueness.}  If $(Y,H,A,C)$ and $(Y',H',A',C')$ are two solutions then
each is minimal, so $Y=Y'$ up to indistinguishability.  The right jumps of $Y$
then determine $C=C'$ through $C_t=-\sum_{s\le t,\,s<T}\Dp Y_s$, which is
forced by \eqref{eq:rbsde-diff} since $\int H\,dX$ and $A$ are
right-continuous, together with the normalisation $C_T=C_{T-}$.  Finally
$\int_0^{\cdot}H\,dX-A=\int_0^{\cdot}H'\,dX-A'$ with $A,A'$ predictable, and the
uniqueness statement in Theorem~\ref{thm:predictable} gives $A=A'$ and
$\int H\,dX=\int H'\,dX$.
\end{proof}

\subsection{Existence, uniqueness and representation}

\begin{defn}\label{def:Phi}
For $Y\in\Sbb^{\infty}$ define $\Gamma^{Y}_t:=\int_0^tf(s,Y_s)\,ds$,
$\zeta^{Y}:=\xi+\Gamma^{Y}$ and
\[
\Phi(Y):=V\bigl(\zeta^{Y}\bigr)-\Gamma^{Y},
\]
where $V(\cdot)$ is the generalised Snell envelope of
Theorem~\ref{thm:snell-aggreg}.  Equivalently, since $\Gamma^{Y}_\theta$ is
$\cF_\theta$-measurable, at every stopping time
\begin{equation}\label{eq:Phi-explicit}
\Phi(Y)_\theta=\esssup_{Q\in\cM(X)}\ \esssup_{\tau\in\cT_\theta}\
\E^{Q}\Bigl[\xi_\tau+\int_\theta^{\tau}f(s,Y_s)\,ds\ \Big|\ \cF_\theta\Bigr],
\qquad \theta\in\cT .
\end{equation}
Since $\Gamma^{Y}$ is a Lebesgue integral, $\Phi(Y)$ depends on $Y$ only through
its class in $\Sbb^{\infty}$.
\end{defn}

The map $\Phi$ is well defined: $\zeta^{Y}$ is optional and bounded by
$c_\xi+T(c_f+L\|Y\|_{\infty})$, so Assumption~\ref{assum:E} holds for it and
$V(\zeta^{Y})\in\Sbb^{\infty}$ is optional and \ladlag; hence
$\Phi(Y)\in\Sbb^{\infty}$.

For $\beta\ge0$ we use on $\Sbb^{\infty}$ the equivalent norms
\[
\|U\|_{\beta}:=\Bigl\|\sup_{0\le t\le T}e^{-\beta(T-t)}|U_t|\Bigr\|_{L^{\infty}(\Prob)}
=\sup_{\tau\in\cT}\bigl\|e^{-\beta(T-\tau)}U_\tau\bigr\|_{L^{\infty}(\Prob)},
\qquad e^{-\beta T}\|U\|_{\infty}\le\|U\|_{\beta}\le\|U\|_{\infty},
\]
the second expression again by the optional section theorem.  The weight is
attached to the time at which $U$ is evaluated, and the supremum runs over
stopping times; a supremum over deterministic times would not be a norm on
$\Sbb^{\infty}$, by the example following the definition of $\Sbb^{\infty}$.

\begin{lem}\label{lem:Phi-properties}
Let Assumption~\ref{assum:data} hold.
\begin{enumerate}\itemsep2pt
\item[(i)] If \eqref{eq:monotone} holds, $\Phi$ is monotone: $U^{1}\le U^{2}$
implies $\Phi(U^{1})\le\Phi(U^{2})$.
\item[(ii)] For every $\beta>0$ and $U^{1},U^{2}\in\Sbb^{\infty}$,
\[
\|\Phi(U^{1})-\Phi(U^{2})\|_{\beta}\ \le\ \frac{L}{\beta}\,\|U^{1}-U^{2}\|_{\beta}.
\]
In particular, with $\beta:=2L+1$, $\Phi$ is a strict contraction on
$(\Sbb^{\infty},\|\cdot\|_{\beta})$, whether or not \eqref{eq:monotone} holds
and including the case $L=0$, and has a unique fixed point.
\end{enumerate}
\end{lem}

\begin{proof}
(i) is immediate from \eqref{eq:Phi-explicit} and \eqref{eq:monotone}, the
integral being taken with respect to a positive measure and the essential
suprema being monotone.

(ii) Put $b:=\|U^{1}-U^{2}\|_{\beta}$.  Outside an evanescent set,
$|U^{1}_s-U^{2}_s|\le b\,e^{\beta(T-s)}$ for all $s$.  Fix $\theta\in\cT$.  The
difference of two essential suprema is bounded by the essential supremum of the
differences, so by \eqref{eq:Phi-explicit} and the Lipschitz property,
\[
\begin{aligned}
\bigl|\Phi(U^{1})_\theta-\Phi(U^{2})_\theta\bigr|
&\le\esssup_{Q}\esssup_{\tau\in\cT_\theta}
\E^{Q}\Bigl[\int_\theta^{\tau}L\,|U^1_s-U^2_s|\,ds\Big|\cF_\theta\Bigr]\\
&\le Lb\int_\theta^{T}e^{\beta(T-s)}ds
=\frac{Lb}{\beta}\bigl(e^{\beta(T-\theta)}-1\bigr),
\end{aligned}
\]
the pathwise bound holding outside one null set and its right-hand side being
$\cF_\theta$-measurable.  Multiplying by $e^{-\beta(T-\theta)}$ and taking the
supremum over $\theta\in\cT$ gives the claim.  Banach's fixed point theorem
applies because $\Sbb^{\infty}$ is complete.
\end{proof}

\begin{lem}[Domination]\label{lem:domination}
Assume \eqref{eq:monotone} and let $U\in\Sbb^{\infty}$ satisfy $U\ge\Phi(U)$.
Then $U\ge Y$, where $Y$ is the fixed point of $\Phi$.
\end{lem}

\begin{proof}
From $\Phi(U)\le U$ and Lemma~\ref{lem:Phi-properties}(i), by induction
$\Phi^{n+1}(U)\le\Phi^{n}(U)\le U$ for every $n$.  By
Lemma~\ref{lem:Phi-properties}(ii), $\|\Phi^{n}(U)-Y\|_{\beta}\le
(L/\beta)^{n}\|U-Y\|_{\beta}\to0$, so $\Phi^{n}(U)\to Y$ uniformly outside an
evanescent set, and passing to the limit in $\Phi^{n}(U)\le U$ gives $Y\le U$.
\end{proof}

\begin{thm}[Existence, uniqueness, representation]\label{thm:existence}
Under Assumption~\ref{assum:data} the reflected BSDE \eqref{eq:rbsde} has a
unique solution $(Y,H,A,C)$ in the sense of Definition~\ref{def:rbsde}, with $A$
predictable, $C_t=-\sum_{s\le t,\,s<T}\Dp Y_s$, and $Y,A,C$ and $\int H\,dX$
unique.  Its first component is the unique fixed point of $\Phi$ in
$\Sbb^{\infty}$ and satisfies the dynamic programming principle
\begin{equation}\label{eq:representation}
Y_\tau=\esssup_{Q\in\cM(X)}\ \esssup_{\theta\in\cT_\tau}\
\E^{Q}\Bigl[\xi_\theta+\int_\tau^{\theta}f(s,Y_s)\,ds\ \Big|\ \cF_\tau\Bigr]
\qquad\text{a.s., for every }\tau\in\cT ,
\end{equation}
and in particular at every deterministic time $\tau=t$.
\end{thm}

\begin{thm}[A priori bounds]\label{thm:apriori}
The solution of Theorem~\ref{thm:existence} satisfies
\begin{equation}\label{eq:apriori}
\|Y\|_{\infty}\le\bigl(c_\xi+Tc_f\bigr)e^{LT}=:c_Y,
\qquad
\|Z\|^2_{\Hbb^2}\le 5c_V^{2},
\qquad
\sup_{Q\in\cM(X)}\E^{Q}\bigl[(A_T+C_T)^{2}\bigr]\le21\,c_V^{2},
\end{equation}
where $Z=\sigma\trans H$ as in \eqref{eq:ZH} and
$c_V:=c_Y+T(c_f+Lc_Y)$.  In particular $\int_0^{\cdot}H\,dX$ is a
square-integrable $Q$-martingale for every $Q\in\cM(X)$, and $A,C\in\Sbb^2$.
\end{thm}

Before proving these two theorems we recall the change-of-variable formula for
optional semimartingales, in the form we use it.  Let
\[
Y=Y_0+Y^{r}+Y^{g},
\]
where $Y^{r}$ is a \cadlag{} semimartingale with $Y^{r}_0=0$, and $Y^{g}$ is an
adapted left-continuous purely discontinuous process of finite variation with
$Y^{g}_0=0$; write $Y^{r}=M^{r}+B^{r}$ with $M^{r}$ a \cadlag{} local
martingale and $B^{r}$ a \cadlag{} adapted process of finite variation, both
null at zero.  Left continuity alone would not make the Stieltjes integral
$\int F'(Y_s)\,dY^{g}_s$ below meaningful.  In our application $Y^{g}=-C_{-}$.

\begin{thm}[Gal'chouk--Lenglart
\cite{Galchouk1981,Lenglart1980}]\label{thm:galchouk}
For every $F\in C^{2}(\R)$ and $t\in[0,T]$,
\[
\begin{aligned}
F(Y_t)=F(Y_0)
&+\int_{(0,t]}F'(Y_{s-})\,dY^{r}_s
+\tfrac12\int_{(0,t]}F''(Y_{s-})\,d\langle M^{r,c}\rangle_s\\[2pt]
&+\sum_{0<s\le t}\bigl[F(Y_s)-F(Y_{s-})-F'(Y_{s-})\Dm Y_s\bigr]\\[2pt]
&+\int_{[0,t)}F'(Y_s)\,dY^{g}_s
+\sum_{0\le s<t}\bigl[F(Y_{s+})-F(Y_s)-F'(Y_s)\Dp Y_s\bigr].
\end{aligned}
\]
Here the integral $\int F'(Y_{s-})dY^{r}_s$ is taken over the \emph{whole}
\cadlag{} part $Y^{r}$, jumps included, and $\langle M^{r,c}\rangle$ is the
bracket of the continuous local martingale part.  When $F$ does not have bounded
derivatives the formula is applied along a localising sequence.
\end{thm}

\begin{proof}[Proof of Theorem~\ref{thm:existence}]
\emph{Step 1: fixed point.}  By Lemma~\ref{lem:Phi-properties}(ii), $\Phi$ has
a unique fixed point $Y\in\Sbb^{\infty}$, and $\Phi^{m}(U)\to Y$ for every
$U\in\Sbb^{\infty}$; monotonicity of $f$ is not needed for this.  The identity
$Y=\Phi(Y)$ holds at every stopping time, $\Phi$ being defined through the
generalised Snell envelope, which is specified at every $\tau\in\cT$; this is
exactly \eqref{eq:representation}.

\emph{Step 2: $Y$ is the first component of a supersolution.}  Put
$\Gamma:=\Gamma^{Y}$, $\zeta:=\xi+\Gamma$ and $V:=V(\zeta)$, so that
$Y=V-\Gamma$.  By Theorem~\ref{thm:zero-driver} applied to the obstacle
$\zeta$ there are $H$ predictable, $A$ predictable non-decreasing and $C$
purely discontinuous non-decreasing and normalised, with
$V_t=V_0+\int_0^tH\,dX-A_t-C_{t-}$ and $\E^{Q}[A_T+C_T]<\infty$ for all $Q$.
Evaluating at $\tau$ and at $T$, subtracting, and using
$V_T=\zeta_T=\xi_T+\Gamma_T$ and $Y=V-\Gamma$,
\[
Y_\tau=V_\tau-\Gamma_\tau
=\xi_T+\Gamma_T-\Gamma_\tau-\int_\tau^TH\,dX+A_T-A_\tau+C_{T-}-C_{\tau-}
\]
which is \eqref{eq:rbsde} because $\Gamma_T-\Gamma_\tau=\int_\tau^Tf(s,Y_s)ds$.
Moreover $Y=V-\Gamma\ge\zeta-\Gamma=\xi$, and $Y\in\Sbb^{\infty}$ is
\ladlag{} because $V$ is \ladlag{} and $\Gamma$ is continuous.  Only
optionality and boundedness of $\zeta$ were needed to construct $V$; no path
limits of $\xi$ are used.  Finally
$\Dp Y=\Dp V$ because $\Gamma$ is continuous, so
$C_t=-\sum_{s\le t,\,s<T}\Dp Y_s$.  Steps~1 and~2 do not use
\eqref{eq:monotone}.

\emph{Step 3: minimality under \eqref{eq:monotone}.}  Let $(Y',H',A',C')$ be
any supersolution, and set $\Gamma':=\Gamma^{Y'}$.  By
Lemma~\ref{lem:super-is-supermart}, $Y'+\Gamma'$ is a strong optional
$Q$-supermartingale for every $Q\in\cM(X)$ dominating $\xi+\Gamma'$; by
Theorem~\ref{thm:snell-minimal}, $Y'+\Gamma'\ge V(\xi+\Gamma')$, that is,
$Y'\ge\Phi(Y')$.  Lemma~\ref{lem:domination} then gives $Y'\ge Y$.  Hence $Y$ is
minimal and $(Y,H,A,C)$ is a solution.

\emph{Step 4: uniqueness.}  Two solutions have equal first components by
minimality; the identification of $C$, $A$ and $\int H\,dX$ then proceeds
exactly as in the proof of Theorem~\ref{thm:zero-driver}.

\emph{Step 5: a general Lipschitz driver.}  Let $f$ be $L$-Lipschitz, not
necessarily monotone, and let $(\bar f,\bar\xi)$ be the data of
Lemma~\ref{lem:monotone-driver}.  Steps~1--4 give a unique solution
$(\bar Y,\bar H,\bar A,\bar C)$ for the barred data, satisfying the barred
representation, and Lemma~\ref{lem:monotone-driver} transforms it into a
unique solution $(Y,H,A,C)$ for $(f,\xi)$.  It remains to show that this $Y$
satisfies \eqref{eq:representation}, i.e.\ is the fixed point of the unbarred
map $\Phi$; the transformation does not transport essential suprema directly,
so we argue through the gaps.  For any supersolution, write
$K:=A+C_-$ and $G:=Y-\xi\ge0$.  Evaluating \eqref{eq:rbsde-diff} at
$\tau\le\theta$ and conditioning, Lemma~\ref{lem:super-is-supermart}(i) gives
\begin{equation}\label{eq:gap-identity}
\E^{Q}\bigl[K_\theta-K_\tau+G_\theta\mid\cF_\tau\bigr]
=Y_\tau-\E^{Q}\Bigl[\xi_\theta+\int_\tau^{\theta}f(s,Y_s)\,ds\Bigm|\cF_\tau\Bigr],
\qquad Q\in\cM(X),\ \theta\in\cT_\tau .
\end{equation}
Applied to the barred solution, the barred representation says that the
essential infimum over $(Q,\theta)$ of the barred left-hand side vanishes.  The
barred integrand is
\[
\int_{(\tau,\theta]}e^{-Ls}\,dA_s+\int_{[\tau,\theta)}e^{-Ls}\,dC_s
+e^{-L\theta}G_\theta\ \ge\ e^{-LT}\bigl(K_\theta-K_\tau+G_\theta\bigr)\ \ge0 .
\]
Hence the essential infimum over $(Q,\theta)$ of the unbarred left-hand side of
\eqref{eq:gap-identity} also vanishes, and \eqref{eq:gap-identity} yields
\eqref{eq:representation}.  By Lemma~\ref{lem:Phi-properties}(ii), which does
not need \eqref{eq:monotone}, $Y$ is then the unique fixed point of $\Phi$.
\end{proof}

\begin{proof}[Proof of Theorem~\ref{thm:apriori}]
\emph{Bound on $Y$.}  Put $U^{0}:=0$, $U^{n+1}:=\Phi(U^{n})$, and
$h_0:=0$, $h_{n+1}(t):=c_\xi+\int_t^T\bigl(c_f+Lh_n(s)\bigr)ds$.  If
$|U^{n}_s|\le h_n(s)$ for all $s$ outside an evanescent set, then
\eqref{eq:Phi-explicit} and $|f(s,y)|\le c_f+L|y|$ give
$|U^{n+1}_\theta|\le h_{n+1}(\theta)$ a.s.\ for every $\theta\in\cT$, hence,
by the optional section theorem, outside an evanescent set.  The functions
$h_n$ increase to the solution $h$ of $h(t)=c_\xi+\int_t^T(c_f+Lh(s))ds$, and
$h(t)\le(c_\xi+Tc_f)e^{L(T-t)}\le c_Y$ by Gronwall's lemma.  Since
$U^{n}\to Y$ in $\Sbb^{\infty}$ (Step~1 of the proof of
Theorem~\ref{thm:existence}), $|Y|\le c_Y$.  Consequently
$V:=Y+\Gamma^{Y}$ satisfies $\|V\|_{\infty}\le c_Y+T(c_f+Lc_Y)=c_V$.

\emph{Bound on $Z$.}  Fix $Q\in\cM(X)$ and recall
$V_t=V_0+\int_0^tZ_s\trans dW^{Q}_s-A_t-C_{t-}$, where $W^{Q}$ is a
$Q$-Brownian motion and $Z=\sigma\trans H$.  Here $Y^{r}=\int Z\trans dW^{Q}-A$
and $Y^{g}=-C_{-}$, with $\langle M^{r,c}\rangle=\int|Z|^{2}ds$,
$\Dm V=-\Dm A$ and $\Dp V=-\Dm C$.  Applying Theorem~\ref{thm:galchouk} with
$F(y)=y^{2}$ on $[0,T]$ along a localising sequence and rearranging,
\begin{equation}\label{eq:GL-identity}
\begin{aligned}
V_0^{2}+\int_0^T|Z_s|^{2}ds
={}&V_T^{2}-2\int_{(0,T]}V_{s-}Z_s\trans dW^{Q}_s
+2\int_{(0,T]}V_{s-}\,dA_s+2\int_{[0,T)}V_s\,dC_s\\
&-\sum_{0<s\le T}(\Dm A_s)^{2}-\sum_{0\le s<T}(\Dm C_s)^{2}.
\end{aligned}
\end{equation}
Note the signs: \eqref{eq:GL-identity} is Theorem~\ref{thm:galchouk} solved for
$V_0^{2}$, so the two reflection integrals, which enter the formula itself with
a negative sign, appear here with a \emph{positive} one; it is the quadratic
jump sums that are favourably signed.  Note also the convention in the last
integral: the left-continuous part is $Y^{g}=-C_{-}$, and the Lebesgue--Stieltjes
measure of $t\mapsto C_{t-}$ assigns to $\{s\}$ the same mass $\Dm C_s$ as $dC$
does, so $\int_{[0,T)}V_s\,d(C_-)_s=\int_{[0,T)}V_s\,dC_s
=\sum_{s<T}V_s\,\Dm C_s$; we write the latter throughout.  Since $|V|\le c_V$, taking expectations (the stochastic
integral being a true martingale along the localisation, and Fatou's lemma
allowing the passage to the limit) gives
\[
\E^{Q}\Bigl[\int_0^T|Z_s|^{2}ds\Bigr]
\le c_V^{2}+2c_V\,\E^{Q}\bigl[A_T+C_{T}\bigr].
\]
By Lemma~\ref{lem:super-is-supermart}, $\int_0^{\cdot}H\,dX$ is a uniformly
integrable $Q$-martingale, and $C_T=C_{T-}$, so
\begin{equation}\label{eq:AC-L1}
\E^{Q}\bigl[A_T+C_{T}\bigr]
=V_0-\E^{Q}[V_T]\le2c_V .
\end{equation}
Hence $\E^{Q}\bigl[\int_0^T|Z|^{2}ds\bigr]\le c_V^{2}+4c_V^{2}=5c_V^{2}$,
uniformly in $Q$.  In particular $Z\in\Hbb^{2}$ and
$\int_0^{\cdot}H\,dX=\int_0^{\cdot}Z\trans dW^{Q}$ is a square-integrable
$Q$-martingale.

\emph{Bound on $A$ and $C$.}  From $A_T+C_{T}=V_0-V_T+\int_0^TH\,dX$ and
$(a+b+c)^{2}\le3(a^{2}+b^{2}+c^{2})$,
\[
\E^{Q}\bigl[(A_T+C_{T})^{2}\bigr]
\le3\Bigl(c_V^{2}+c_V^{2}+\E^{Q}\Bigl[\int_0^T|Z|^{2}ds\Bigr]\Bigr)
\le3\bigl(2c_V^{2}+5c_V^{2}\bigr)=21c_V^{2},
\]
uniformly in $Q$.  Since $A$ and $C$ are non-decreasing and non-negative, this
gives $A,C\in\Sbb^{2}$; here the normalisation $C_T=C_{T-}$ is what makes a
bound on $C_{T-}$ control $\sup_t C_t$.
\end{proof}

\subsection{Comparison}\label{subsec:comparison}

Reflected BSDEs are expected to come with a comparison theorem, and here it is a
corollary of the same monotonicity that produced minimality.  No linearisation
of the driver is involved, and --- in contrast with the classical proofs --- no
control whatever of the reflection terms is needed, since the argument never
leaves the level of the map $\Phi$.

\begin{thm}[Comparison]\label{thm:comparison}
For $i=1,2$ let the data $(f^{i},\xi^{i})$ satisfy
Assumption~\ref{assum:data}, with $f^{i}$ Lipschitz in $y$ with constant
$L_i$, and let
$(Y^{i},H^{i},A^{i},C^{i})$ be the corresponding solutions of \eqref{eq:rbsde}
given by Theorem~\ref{thm:existence}.  Assume
\[
\xi^{1}_t\le\xi^{2}_t\quad\text{for all }t,\ \text{a.s.},
\qquad\text{and}\qquad
f^{1}(t,\omega,y)\le f^{2}(t,\omega,y)\quad
\text{for all }(t,\omega,y).
\]
Then $Y^{1}_t\le Y^{2}_t$ for all $t\in[0,T]$, a.s.
\end{thm}

\begin{proof}
Set $L:=\max(L_1,L_2)$, so that both drivers are $L$-Lipschitz in $y$, and apply
Lemma~\ref{lem:monotone-driver} to the two sets of data \emph{with this same
constant}: the transformation is $\bar Y_t=e^{-Lt}Y_t$ and
$\bar f^{i}(t,y)=e^{-Lt}f^{i}(t,e^{Lt}y)+Ly$.  Since the same $L$ is used for
$i=1$ and $i=2$ and $e^{-Lt}>0$, the pointwise inequalities
$\bar f^{1}\le\bar f^{2}$ and $\bar\xi^{1}\le\bar\xi^{2}$ are preserved, as is
the inequality $\bar Y^{1}\le\bar Y^{2}$ to be proved; and each $\bar f^{i}$ is
non-decreasing in $y$ with Lipschitz constant at most $2L$, by
Lemma~\ref{lem:monotone-driver} applied to the $L$-Lipschitz driver $f^{i}$.  We
may therefore assume that $f^{1}$ and $f^{2}$ are non-decreasing in $y$, i.e.\
that \eqref{eq:monotone} holds for both.

Let $\Phi^{i}$ denote the map of Definition~\ref{def:Phi} built from the data
$(f^{i},\xi^{i})$.  For any $U\in\Sbb^{\infty}$ we have, pointwise in
$(s,\omega)$, $f^{1}(s,U_s)\le f^{2}(s,U_s)$, and $\xi^{1}_\theta\le
\xi^{2}_\theta$ for every $\theta\in\cT$; since the essential suprema in
\eqref{eq:Phi-explicit} are monotone and $\int_t^{\theta}\cdot\,ds$ integrates
against a positive measure,
\begin{equation}\label{eq:Phi-order}
\Phi^{1}(U)\ \le\ \Phi^{2}(U),\qquad U\in\Sbb^{\infty}.
\end{equation}
Taking $U=Y^{2}$ in \eqref{eq:Phi-order} and using that $Y^{2}$ is a fixed point
of $\Phi^{2}$ (Theorem~\ref{thm:existence}),
\[
\Phi^{1}(Y^{2})\ \le\ \Phi^{2}(Y^{2})\ =\ Y^{2},
\]
so $Y^{2}\in\Sbb^{\infty}$ dominates its own image under $\Phi^{1}$.  By
Lemma~\ref{lem:domination} applied to $\Phi^{1}$, whose unique fixed point is
$Y^{1}$, we conclude $Y^{2}\ge Y^{1}$.
\end{proof}

\begin{rem}[Comparison and minimality are one argument]\label{rem:comp-min}
Lemma~\ref{lem:domination} says that any bounded optional process dominating its
own image under $\Phi$ dominates the solution.  Minimality is that statement
applied to a supersolution of the \emph{same} problem, which satisfies
$Y'\ge\Phi(Y')$ by Lemma~\ref{lem:super-is-supermart} and
Theorem~\ref{thm:snell-minimal}; comparison is the same statement applied to the
solution of a \emph{dominating} problem, which satisfies $Y^{2}\ge\Phi^{1}(Y^{2})$
by \eqref{eq:Phi-order}.  The two theorems are therefore one argument run from
two different starting points.  Note also that
Lemma~\ref{lem:domination} requires of $U$ only that it lie in $\Sbb^{\infty}$:
the map $\Phi$ of Definition~\ref{def:Phi} is defined on all bounded optional
processes and not merely on those dominating the obstacle, so no admissibility
check on the starting point of the iteration is needed.  Boundedness of the data,
or more generally domination by a traded martingale
(Section~\ref{subsec:weighted}), is what makes this free.
\end{rem}

\subsection{A dynamic characterisation of the solution}

Minimality is a global requirement, and it is natural to ask for a condition on
$(Y,H,A,C)$ stated at each stopping time and equivalent to it --- the analogue
of the Skorokhod conditions in the classical theory.  Such a condition exists,
and it has the shape of the minimum condition used in the theory of
second-order (reflected) BSDEs
\cite{SonerTouziZhang2012,MatoussiPossamaiZhou2013,MPZcorrigendum}, where the
same difficulty is met in the non-dominated setting.  The condition is
dynamic, not pathwise: the essential infimum ranges over the whole family
$\cM(X)$.  It is the dynamic programming principle \eqref{eq:representation}
rewritten in the variables $(A,C,Y-\xi)$, and we count it as a
reformulation, not as a separate result; its use is to make a candidate
quadruple verifiable (as in Theorem~\ref{thm:nonlinear-wealth}) and to state
Problem~P3.

\begin{cor}[Dynamic minimality condition]\label{cor:dyn-minimality}
Let $(Y,H,A,C)$ be a supersolution.  Then $(Y,H,A,C)$ is the solution of
\eqref{eq:rbsde} if and only if
\begin{equation}\label{eq:min-condition}
\essinf_{Q\in\cM(X),\ \theta\in\cT_\tau}\
\E^{Q}\Bigl[\bigl(A_\theta-A_\tau\bigr)
+\bigl(C_{\theta-}-C_{\tau-}\bigr)
+\bigl(Y_\theta-\xi_\theta\bigr)\ \Big|\ \cF_\tau\Bigr]=0
\qquad\text{for every }\tau\in\cT .
\end{equation}
The three terms inside the conditional expectation are non-negative, and the
family over which the essential infimum is taken is downward directed, so that
\eqref{eq:min-condition} is attained as the decreasing limit of a sequence.
\end{cor}

\begin{proof}
\emph{Directedness.}  Given $(Q_1,\theta_1)$ and $(Q_2,\theta_2)$, let $B$ be
the $\cF_\tau$-measurable set on which the first conditional expectation does
not exceed the second, let $\theta:=\theta_1\ind_B+\theta_2\ind_{B^{c}}
\in\cT_\tau$, and let $Q$ paste $Q_1$ on $B$ and $Q_2$ on $B^{c}$ at time
$\tau$, which lies in $\cM(X)$ by Lemma~\ref{lem:event-paste}.  The pair $(Q,\theta)$ realises the minimum of the
two conditional expectations, whence the directedness.

\emph{The identity.}  Evaluating \eqref{eq:rbsde-diff} at the stopping times
$\theta$ and $\tau$ with $\tau\le\theta$ and subtracting,
\[
Y_\theta-Y_\tau+\Gamma_\theta-\Gamma_\tau
=\int_\tau^{\theta}H_s\,dX_s-\bigl(A_\theta-A_\tau\bigr)
-\bigl(C_{\theta-}-C_{\tau-}\bigr),
\]
where $\Gamma_t=\int_0^tf(s,Y_s)ds$.  Taking $\E^{Q}[\,\cdot\mid\cF_\tau]$ and using
Lemma~\ref{lem:super-is-supermart}(i), by which
$\E^{Q}\bigl[\int_\tau^{\theta}H\,dX\mid\cF_\tau\bigr]=0$, we obtain
\begin{equation}\label{eq:min-identity}
\E^{Q}\Bigl[\bigl(A_\theta-A_\tau\bigr)+\bigl(C_{\theta-}-C_{\tau-}\bigr)
+\bigl(Y_\theta-\xi_\theta\bigr)\Big|\cF_\tau\Bigr]
=Y_\tau-\E^{Q}\Bigl[\xi_\theta+\int_\tau^{\theta}f(s,Y_s)\,ds\Big|\cF_\tau\Bigr].
\end{equation}
Taking the essential infimum over $(Q,\theta)$ on both sides,
\eqref{eq:min-condition} holds if and only if
\begin{equation}\label{eq:DPP}
Y_\tau=\esssup_{Q\in\cM(X)}\ \esssup_{\theta\in\cT_\tau}\
\E^{Q}\Bigl[\xi_\theta+\int_\tau^{\theta}f(s,Y_s)\,ds\Big|\cF_\tau\Bigr]
\qquad\text{for every }\tau\in\cT ,
\end{equation}
that is, if and only if $Y_\tau=\Phi(Y)_\tau$ for every $\tau\in\cT$, the map
$\Phi$ of Definition~\ref{def:Phi} being defined at all stopping times through
the generalised Snell envelope.  By the optional section theorem this means
$Y=\Phi(Y)$, and by Lemma~\ref{lem:Phi-properties} the fixed point of $\Phi$ is
unique and is, by Theorem~\ref{thm:existence}, the first component of the
solution.  The remaining components are then identified as in the proof of
Theorem~\ref{thm:existence}.
\end{proof}

\begin{rem}[Where the content sits]\label{rem:regularity-role}
Both directions of Corollary~\ref{cor:dyn-minimality} are read off the single
identity \eqref{eq:min-identity}.  What makes it an identity rather than an
inequality is that $\int H\,dX$ is a true $Q$-martingale, which follows from
the integrability clause $\E^{Q}[A_T+C_T]<\infty$ in
Definition~\ref{def:rbsde}(i) through the domination in the proof of
Lemma~\ref{lem:super-is-supermart}; that clause is the place to look when
transporting the statement to unbounded data.  For each fixed $\tau$,
\eqref{eq:min-condition} is an almost-sure identity between
$\cF_\tau$-measurable variables, so no aggregation is needed.
\end{rem}

\begin{rem}
Condition \eqref{eq:min-condition} says that one may choose a measure and a
stopping time making \emph{simultaneously} small the increment of the increasing
process and the gap to the obstacle.  It replaces the Skorokhod conditions,
which, for obstacles with left limits, impose current contact for the
continuous part and left contact for predictable jumps.  By
Proposition~\ref{pro:no-skorokhod}, those conditions fail for some data in an
incomplete market; \eqref{eq:min-condition} demands only that flatness and contact be
approached along a common near-optimal pair $(Q,\theta)$.  When
$\cM(X)=\{\Prob\}$ the essential infimum is over $\theta$ alone and one recovers
the classical characterisation of the Snell envelope.  Its contact conditions
use the left upper envelope in general and reduce to the formulas with
$\xi_{-}$ when the obstacle has left limits
(Remark~\ref{rem:complete-case}).
\end{rem}

\subsection{The Skorokhod condition for \texorpdfstring{$C$}{C} holds automatically}

Proposition~\ref{pro:no-skorokhod} shows that in an incomplete market the
Skorokhod conditions for $A$ fail for some data.  The condition for $C$ holds
for all data, so the obstruction concerns $A$ alone.

\begin{lem}[Dynamic programming at $\tau+$]\label{lem:Y-max}
Let $(Y,H,A,C)$ be the solution of \eqref{eq:rbsde}, and adopt the convention
$Y_{T+}:=Y_T$.  Then
\begin{equation}\label{eq:Y-max}
Y_\tau=\max\bigl(\xi_\tau,\;Y_{\tau+}\bigr)\quad\text{a.s.,
\qquad for every }\tau\in\cT .
\end{equation}
\end{lem}

The statement is at stopping times, not merely at deterministic times, because
that is how Proposition~\ref{pro:C-skorokhod} uses it.  The two are not
interchangeable here: $Y$ is \ladlag\ and need not be right-continuous, so
almost-sure agreement of two processes at every deterministic time does not make
them indistinguishable, and the section theorem --- which would upgrade the
deterministic statement --- requires as its input exactly the inequality at
every stopping time that is to be proved.  The argument below is therefore
carried out at $\tau$ from the start; the only change from the deterministic
case is that the approximating times must be truncated at $T$.

\begin{proof}
``$\ge$''.  Taking $\theta=\tau$ and any $Q$ in \eqref{eq:representation} gives
$Y_\tau\ge\xi_\tau$, and $Y_\tau\ge Y_{\tau+}$ because $\Dp Y=\Dp V\le0$ by
\eqref{eq:rightjumpV}, the drift $\Gamma$ being continuous.  On $\{\tau=T\}$
both assertions are the equality $Y_T=\xi_T=Y_{T+}$.

``$\le$''.  Fix $Q\in\cM(X)$ and $\theta\in\cT_\tau$, and set
$B:=\{\theta>\tau\}$.  Since $\{\theta>\tau\}=\{\theta\le\tau\}^{c}$ and
$\{\theta\le\tau\}\in\cF_\tau$ for $\theta,\tau\in\cT$, we have $B\in\cF_\tau$;
note also that $B\subseteq\{\tau<T\}$.  For $n\ge1$ put
\[
\tau_n:=\bigl(\tau+\tfrac1n\bigr)\wedge T,
\qquad
\sigma_n:=\theta\vee\tau_n ,
\]
both stopping times, the first because $\tau\in\cT$ and the second as the
maximum of two stopping times; and $\sigma_n\in\cT_{\tau_n}$.  On $B$ one has
$\sigma_n=\theta$ for every $n$ with $\tfrac1n<\theta-\tau$, since then
$\tau_n\le\tau+\tfrac1n<\theta$; the threshold depends on $\omega$, which is
harmless because the passage to the limit below is pathwise.  Applying
\eqref{eq:representation} at the stopping time $\tau_n$, with $\sigma_n$ as the
competing stopping time, and then conditioning on $\cF_\tau$,
\[
\E^{Q}\Bigl[\xi_{\sigma_n}+\int_{\tau_n}^{\sigma_n}f(s,Y_s)ds\;\Big|\;\cF_\tau\Bigr]
\ \le\ \E^{Q}\bigl[Y_{\tau_n}\mid\cF_\tau\bigr],
\]
so that, adding $\E^{Q}[\int_\tau^{\tau_n}f(s,Y_s)ds\mid\cF_\tau]$ to both
sides, bounding that term on the right by
$\E^{Q}[\int_\tau^{\tau_n}|f(s,Y_s)|ds\mid\cF_\tau]$, and multiplying by
$\ind_B$,
\[
\ind_B\,\E^{Q}\Bigl[\xi_{\sigma_n}+\int_{\tau}^{\sigma_n}f(s,Y_s)ds\;\Big|\;\cF_\tau\Bigr]
\ \le\ \ind_B\,\E^{Q}\bigl[Y_{\tau_n}\mid\cF_\tau\bigr]
+\ind_B\,\E^{Q}\Bigl[\int_\tau^{\tau_n}|f(s,Y_s)|ds\;\Big|\;\cF_\tau\Bigr].
\]
Now let $n\to\infty$.  On $B$ one has $\sigma_n=\theta$ for $n$ large, so
$\xi_{\sigma_n}\to\xi_\theta$ and $\int_\tau^{\sigma_n}f\to\int_\tau^{\theta}f$
pointwise by eventual equality, without any path regularity of $\xi$;
$\xi$ is bounded by {\bf (A1)} and $|f(s,Y_s)|\le c_f+Lc_Y$ by
Theorem~\ref{thm:apriori}, so conditional dominated convergence applies to every
term, and the last term tends to $0$ because $\tau_n-\tau\to0$.  For the
right-hand side the chain is: $\tau_n\downarrow\tau$ with $\tau_n>\tau$ on
$\{\tau<T\}$ and $\tau_n=T$ on $\{\tau=T\}$, so $Y$ being \ladlag\ gives
$Y_{\tau_n}\to Y_{\tau+}$ a.s.\ on all of $\Omega$, the convention $Y_{T+}=Y_T$
covering $\{\tau=T\}$; $Y$ is bounded, so
$\E^{Q}[Y_{\tau_n}\mid\cF_\tau]\to\E^{Q}[Y_{\tau+}\mid\cF_\tau]$; and $Y_{\tau+}$
is $\cF_{\tau+}$-measurable with $\cF_{\tau+}=\cF_\tau$, $\F$ being
right-continuous, so $\E^{Q}[Y_{\tau+}\mid\cF_\tau]=Y_{\tau+}$.
Hence
\[
\ind_B\,\E^{Q}\Bigl[\xi_{\theta}+\int_\tau^{\theta}f(s,Y_s)ds\;\Big|\;\cF_\tau\Bigr]
\ \le\ \ind_B\,Y_{\tau+},
\]
while on $B^{c}=\{\theta=\tau\}$ the left-hand side equals $\ind_{B^{c}}\xi_\tau$.
Adding the two gives, for every such pair $(Q,\theta)$,
\[
\E^{Q}\Bigl[\xi_{\theta}+\int_\tau^{\theta}f(s,Y_s)ds\;\Big|\;\cF_\tau\Bigr]
\ \le\ \max\bigl(\xi_\tau,Y_{\tau+}\bigr),
\]
and taking the essential supremum over $(Q,\theta)$ in
\eqref{eq:representation} gives $Y_\tau\le\max(\xi_\tau,Y_{\tau+})$.
\end{proof}

\begin{pro}[$C$ reflects on the contact set]\label{pro:C-skorokhod}
Let $(Y,H,A,C)$ be the solution of \eqref{eq:rbsde}.  Then
\[
(Y_\tau-\xi_\tau)\,\Dm C_\tau=0\qquad\text{for every }\tau\in\cT ,
\]
that is, $C$ increases only on the contact set $\{Y=\xi\}$.
\end{pro}

\begin{proof}
By Theorem~\ref{thm:existence}, $\Dm C_\tau=-\Dp Y_\tau=Y_\tau-Y_{\tau+}$.  If
$\Dm C_\tau>0$ then $Y_\tau>Y_{\tau+}$, so the maximum in \eqref{eq:Y-max} is
not attained by its second argument; hence $Y_\tau=\xi_\tau$.
\end{proof}

\begin{rem}\label{rem:Y-max}
Identity \eqref{eq:Y-max} is the dynamic programming principle read at $\tau+$:
\emph{stop now, or continue, the value of continuing being the value at
$\tau+$}.  No left or right limits of $\xi$ are used --- only its boundedness,
through dominated convergence, or its domination by a traded martingale
(Corollary~\ref{cor:weighted-comparison}) --- so
Proposition~\ref{pro:C-skorokhod} holds in the full generality of
Section~\ref{sec:bsde}.
\end{rem}

\begin{rem}[The asymmetry, and what it does to Problem~P3]\label{rem:asymmetry}
Propositions~\ref{pro:no-skorokhod} and \ref{pro:C-skorokhod} together say that
of the three classical Skorokhod conditions listed in
Remark~\ref{rem:why-minimal}, the two governing $A$ fail for some data in an
incomplete market while the one governing $C$ holds for all data.  Reflection
and the premium for non-attainability are confounded in the predictable part
only; the purely
discontinuous part remains characterised pathwise by the right jumps of $Y$,
exactly as in the complete-market theory.  Problem~P3 should therefore be posed
for $A$ alone --- $C$ already factors.
\end{rem}

\begin{rem}[Comparison with second-order reflected BSDEs]\label{rem:2bsde}
In the non-dominated setting of Matoussi, Possama\"i and Zhou
\cite{MatoussiPossamaiZhou2013} the solution of a second-order reflected BSDE
consists of a value $Y$ and an integrand $Z$ common to all priors, aggregated
quasi-surely, together with a family of increasing processes $\{K^{\Prob}\}$
indexed by the priors.  They observe \cite[Remark~3.2, p.~2429]{MatoussiPossamaiZhou2013}
that $K^{\Prob}$ does not satisfy the Skorokhod condition, being formally the
sum of an increasing process due to the second-order structure and one due to
the reflection, of which only the second is bound to satisfy it, and they leave
open the corresponding decomposition of $K^{\Prob}$
\cite[Remark~3.5, p.~2433]{MatoussiPossamaiZhou2013}.  The uniqueness of their
solution rests on a minimum condition across priors, corrected in
\cite{MPZcorrigendum}; when the family of priors is a singleton that condition
forces $K^{\Prob}$ to coincide with the reflection process of the classical
reflected BSDE, and under an assumption on the oscillation of the obstacle it
is equivalent to a Skorokhod-type condition with an essential infimum over
priors \cite{MPZcorrigendum}.

Proposition~\ref{pro:no-skorokhod} and Corollary~\ref{cor:dichotomy} are
therefore the dominated counterpart of a phenomenon identified there, not a new
phenomenon.  What the dominated setting adds is that the measures are
equivalent: the regulator $A$ of \eqref{eq:rbsde} is a single process, equal
to the predictable Doob--Meyer part of $Y+\Gamma$ under every $Q\in\cM(X)$
simultaneously (by the remark following Theorem~\ref{thm:predictable} and the
uniform decomposition), and the obstruction can be stated as an exact criterion
(bounded attainability) rather than as a remark.  Equivalence does not by
itself prevent splitting $A$ into two increasing processes with different
roles, as Peng and Xu do for a right-continuous value
(Remark~\ref{rem:pengxu}); Problem~P3 asks for such a split in the \ladlag{}
case.
\end{rem}

\begin{rem}[The complete-market case]\label{rem:complete-case}
Suppose $\mathrm{rank}\,\sigma_s=n$ a.e., so that $\cM(X)=\{\Prob\}$ and the
constraint $\cK_t=\R^{n}$ is vacuous.  Then \eqref{eq:representation} says that
$V:=Y+\Gamma^{Y}$ is the classical $\Prob$-Snell envelope of the bounded
optional reward $\eta:=\xi+\Gamma^{Y}$, and $(H,A,C)$ come from its Mertens
decomposition.  For $t>0$ put
$\bar\eta_t:=\limsup_{s\uparrow t,\ s<t}\eta_s$, the predictable left upper
envelope.  The classical contact conditions are
\[
\int_{(0,T]}\ind_{\{V_{t-}>\bar\eta_t\}}\,dA^{c}_t=0,\qquad
(V_{\tau-}-\bar\eta_\tau)\Dm A_\tau=0,\qquad
(V_\theta-\eta_\theta)\Dm C_\theta=0,
\]
for predictable stopping times $0<\tau\le T$ and stopping times $0\le\theta\le T$
\cite{ElKaroui1981,Grigorova2020}.

Assume now, for the reduction to the formulas of
Remark~\ref{rem:why-minimal}, that $\xi$ has finite left limits on $(0,T]$.
Then $\bar\eta=\eta_-$, since $\Gamma^{Y}$ is continuous.  A path with finite
left limits differs from its left-limit path at at most countably many times;
the same is true of $V$.  The continuous measure $dA^{c}$ ignores these
times, so its contact condition can also be written using $V_t$ and $\eta_t$.
Subtracting $\Gamma^{Y}$ gives the three conditions of
Remark~\ref{rem:why-minimal}.  Conversely, under this left-limit assumption,
a supersolution $Y'$ satisfying those conditions is the Snell envelope of
$\xi+\Gamma^{Y'}$ minus $\Gamma^{Y'}$, hence $Y'=\Phi(Y')$ and is the
solution.  This is the singleton case of Corollary~\ref{cor:dichotomy}.

The left-limit qualification cannot be omitted from the formulas with current
contact.  For $X=W$, $T=1$, $f=0$ and
$\xi_t=(1-t)\ind_{\mathbb Q\cap[0,1]}(t)$, the value is $Y_t=1-t$, with
$H=C=0$ and $A_t=t$: dyadic exercise times decreasing to any stopping time
give the lower bound, and no later reward can exceed $1-t$.  At an irrational
initial time $t<1$ the supremum is not attained by any stopping time.
Nevertheless
$\int_0^1\ind_{\{Y_t>\xi_t\}}\,dA_t=1$.  The left upper envelope is
$\bar\xi_t=1-t$, so the general classical contact condition holds.  This
bounded optional example is covered by Theorem~\ref{thm:existence}, although
its obstacle has neither left nor right limits at any $t\in(0,1)$.

The classical argument uses that the supermartingale
property and the optimality are expressed under one and the same measure; when
$\cM(X)$ is not a singleton, $\Prob$-a.s.\ vanishing of increments would have to
be deduced from $\essinf_{Q\in\cM(X)}\E^{Q}[\,\cdot\mid\cF_t]=0$, which is
weaker.  This is Problem~P3.
\end{rem}

\subsection{Unbounded data: a weighted solution class}\label{subsec:weighted}

Boundedness of the data enters the existence and minimality argument above in
three places: through the space $\Sbb^{\infty}$ on which $\Phi$ contracts,
through Theorem~\ref{thm:zero-driver} in Step~2 of the proof of
Theorem~\ref{thm:existence}, and through Lemma~\ref{lem:super-is-supermart},
where it gives $\int H\,dX$ a constant lower bound.  All three can be replaced
by domination by a single traded martingale, and no maximal inequality is
needed.  The square-integrable bounds of Theorem~\ref{thm:apriori} rest on
$|V|\le c_V$ and do not survive (Example~\ref{eg:weighted-not-L2}).

\begin{assum}[Weighted data]\label{assum:weighted}
A \emph{weight} is a \cadlag{} process $M\ge1$ which is a uniformly integrable
$Q$-martingale for every $Q\in\cM(X)$.  For a weight $M$ and constants
$c_\xi,c_f\ge0$, the obstacle $\xi$ satisfies {\bf (A1)} with the bound
replaced by $|\xi_t|\le c_\xi M_t$, and the driver $f$ satisfies {\bf (A2)}
with the bound on $f(t,\omega,0)$ replaced by
$|f(t,\omega,0)|\le c_fM_t(\omega)$.  In particular, the obstacle is only
optional; no path limits are assumed in the weighted theory either.
\end{assum}

By Corollary~\ref{cor:uniform-ui-representation}, $M-M_0$ is a stochastic
integral of $X$, so a weight is of the traded form of
Assumption~\ref{assum:G}; in the present Brownian filtration it is continuous.
Assumption~\ref{assum:data} is the case $M\equiv1$.  For a weight $M$ and
$\beta\ge0$ let $\Sbb^{M}$ be the space of optional processes $Y$ with
\[
\|Y\|_{M}:=\Bigl\|\sup_{0\le t\le T}\frac{|Y_t|}{M_t}\Bigr\|_{L^{\infty}(\Prob)}<\infty,
\qquad
\|Y\|_{\beta,M}:=\Bigl\|\sup_{0\le t\le T}e^{-\beta(T-t)}\,
\frac{|Y_t|}{M_t}\Bigr\|_{L^{\infty}(\Prob)} ,
\]
and let $\Sbb^{\rm tr}:=\bigcup_M\Sbb^{M}$, the union over all weights: the
optional processes dominated in absolute value by a traded martingale.  The map
$Y\mapsto Y/M$ is an isometry of $(\Sbb^{M},\|\cdot\|_{\beta,M})$ onto
$(\Sbb^{\infty},\|\cdot\|_{\beta})$, so $\Sbb^{M}$ is complete, and as there the
suprema may be taken over stopping times.  If $M\le M'$ then
$\Sbb^{M}\subset\Sbb^{M'}$; and $\Sbb^{\rm tr}$ is a vector space, a sum of two
weights being a weight.

A \emph{supersolution of class~(D)} is a quadruple satisfying
Definition~\ref{def:rbsde} with ``$Y\in\Sbb^{\infty}$'' in (i) replaced by
``$Y$ is of class~(D) under every $Q\in\cM(X)$'', and a \emph{solution of
class~(D)} is one whose first component is minimal among these.  Every process
in $\Sbb^{\rm tr}$ is of class~(D) under every $Q$, a weight being a uniformly
integrable martingale; Lemma~\ref{lem:classD-tr} below shows that, for weighted
data, the first component of every supersolution of class~(D) lies in
$\Sbb^{\rm tr}$.  The device that replaces a maximal inequality is the
following identity.

\begin{lem}[Time integrals of a weight]\label{lem:weight-integral}
Let $M$ be a weight, $Q\in\cM(X)$, $\theta\in\cT$ and $\varphi\ge0$ a bounded
Borel function on $[0,T]$.  Then
\begin{equation}\label{eq:weight-integral}
\E^{Q}\Bigl[\int_\theta^T\varphi(s)M_s\,ds\Bigm|\cF_\theta\Bigr]
=M_\theta\int_\theta^T\varphi(s)\,ds .
\end{equation}
Consequently $J_t:=\int_0^tM_s\,ds+(T-t)M_t$ is a non-negative uniformly
integrable $Q$-martingale for every $Q\in\cM(X)$, with
$J_T=\int_0^TM_s\,ds$, $\E^{Q}[J_T]=TM_0$ and $J\ge\int_0^{\cdot}M_s\,ds$; and
$J=TM_0+\int_0^{\cdot}(T-s)\,dM_s$ is again traded.
\end{lem}

\begin{proof}
Let $B\in\cF_\theta$ and $s\in[0,T]$.  Since $B\cap\{\theta\le s\}\in\cF_\theta$
and $\E^{Q}[M_s\mid\cF_\theta]=M_{\theta\wedge s}$ by optional sampling,
$\E^{Q}[\ind_{B\cap\{\theta\le s\}}M_s]=\E^{Q}[\ind_{B\cap\{\theta\le s\}}M_\theta]$.
Multiplying by $\varphi(s)$ and integrating over $s$, Tonelli's theorem, which
applies because every term is non-negative, gives
$\E^{Q}[\ind_B\int_\theta^T\varphi(s)M_s\,ds]
=\E^{Q}[\ind_BM_\theta\int_\theta^T\varphi(s)\,ds]$; the right-hand side of
\eqref{eq:weight-integral} is $\cF_\theta$-measurable.  With $\varphi\equiv1$
and $\theta=t$ this reads $\E^{Q}[J_T\mid\cF_t]=J_t$, and with $t=0$ it gives
$\E^{Q}[J_T]=TM_0$.  The last identity is integration by parts for the product
of $M$ and the deterministic function $T-t$.
\end{proof}

\begin{lem}[Supersolutions of class~(D)]\label{lem:weighted-super}
Let Assumption~\ref{assum:weighted} hold for a weight $M$, and let $(Y,H,A,C)$
be a supersolution of class~(D).  Put $\Gamma_t:=\int_0^tf(s,Y_s)\,ds$ and
$\overline\Gamma:=c_fJ_T+L\int_0^T|Y_s|\,ds$.  Then $|\Gamma_t|\le\overline\Gamma$
for all $t$, $\overline\Gamma\in L^{1}(Q)$, and for every $Q\in\cM(X)$:
\begin{enumerate}\itemsep2pt
\item[(i)] $\int_0^{\cdot}H\,dX$ is a uniformly integrable $Q$-martingale;
\item[(ii)] $Y+\Gamma$ is a strong optional $Q$-supermartingale of class~(D)
which dominates $\xi+\Gamma$.
\end{enumerate}
If moreover $Y\in\Sbb^{M}$ and $b:=\|Y\|_M$, then
$\overline\Gamma\le c_bJ_T$ with $c_b:=c_f+Lb$.
\end{lem}

\begin{proof}
The identity \eqref{eq:rbsde-diff} was derived from \eqref{eq:rbsde} and the
optional section theorem alone, and holds here.  Since
$|f(s,Y_s)|\le c_fM_s+L|Y_s|$, we have $|\Gamma_t|\le\overline\Gamma$, and
$\E^{Q}[\overline\Gamma]\le c_fTM_0+LT\sup_{\tau\in\cT}\E^{Q}|Y_\tau|<\infty$, a
class-(D) family being bounded in $L^{1}(Q)$.  Fix $Q$ and put
$N:=\int_0^{\cdot}H\,dX$, a continuous $Q$-local martingale because $X$ is one.
By \eqref{eq:rbsde-diff}, with $A,C\ge0$ non-decreasing, for every $\tau\in\cT$
\begin{equation}\label{eq:weighted-N}
|N_\tau|\ \le\ |Y_\tau|+|Y_0|+\overline\Gamma+A_T+C_T .
\end{equation}
The family $\{Y_\tau\}_{\tau\in\cT}$ is uniformly $Q$-integrable, and
$\overline\Gamma$ and $A_T+C_T$ lie in $L^{1}(Q)$, the latter by
Definition~\ref{def:rbsde}(i).  Hence $N$ is a $Q$-local martingale of
class~(D), that is, a uniformly integrable $Q$-martingale, which is (i).  The
supermartingale inequality in (ii) then follows exactly as in
Lemma~\ref{lem:super-is-supermart}, class~(D) from
$|Y_\tau+\Gamma_\tau|\le|Y_\tau|+\overline\Gamma$, and the domination from
Definition~\ref{def:rbsde}(iii).  If $|Y|\le bM$, then
$L\int_0^T|Y_s|\,ds\le LbJ_T$.
\end{proof}

Estimate \eqref{eq:weighted-N} takes the place of the constant lower bound
used in Lemma~\ref{lem:super-is-supermart}, and it uses the integrability of the
supersolution at hand, not the a~priori bound of the solution.  In a market
with jumps, where $N$ is a priori only a sigma-martingale, the same identity
gives $N^{-}_\tau\le|Y_\tau|+|Y_0|+\overline\Gamma$, uniformly $Q$-integrable in
$\tau$, and Lemma~\ref{lem:BS} supplies the local martingale property before
\eqref{eq:weighted-N} is used.

\begin{lem}[Supersolutions of class~(D) are dominated by a traded martingale]
\label{lem:classD-tr}
Let Assumption~\ref{assum:weighted} hold for a weight $M$, and let $(Y,H,A,C)$
be a supersolution of class~(D).  Then there is a non-negative \cadlag{}
process $\widetilde M$ which is a uniformly integrable $Q$-martingale for every
$Q\in\cM(X)$ and satisfies $|Y|\le e^{LT}\widetilde M$.  In particular
$Y\in\Sbb^{1+\widetilde M}\subset\Sbb^{\rm tr}$.
\end{lem}

\begin{proof}
Put $b_s:=(f(s,Y_s)-f(s,0))/Y_s$ on $\{Y_s\neq0\}$ and $b_s:=0$ otherwise, so
that $|b|\le L$ and $f(s,Y_s)=f(s,0)+b_sY_s$, and put
$R_t:=\exp(\int_0^tb_s\,ds)$, so that $e^{-LT}\le R\le e^{LT}$.  The process
$R$ is continuous, adapted and of finite variation.  Applying the product rule
to $R$ and \eqref{eq:rbsde-diff} --- to the \cadlag{} part in the usual way and
to the left-continuous part $C_-$ as a Stieltjes integral, with no bracket
terms since $R$ is continuous and of finite variation --- gives
\begin{equation}\label{eq:classD-linear}
R_tY_t=Y_0+\widetilde N_t-\int_{(0,t]}R_s\,dA_s-\int_{[0,t)}R_s\,dC_s
-\int_0^tR_sf(s,0)\,ds,
\qquad \widetilde N:=\int_0^{\cdot}R_sH_s\,dX_s .
\end{equation}
Fix $Q$.  The process $\widetilde N$ is a continuous $Q$-local martingale and
is traded.  By \eqref{eq:classD-linear},
$|\widetilde N_\tau|\le e^{LT}(|Y_\tau|+A_T+C_T+c_fJ_T)+|Y_0|$.
Here $\{Y_\tau:\tau\in\cT\}$ is uniformly $Q$-integrable by hypothesis,
$A_T+C_T\in L^1(Q)$ by Definition~\ref{def:rbsde}(i), which is retained in
the definition of a supersolution of class~(D), and $J_T\in L^1(Q)$ by
Lemma~\ref{lem:weight-integral}.  The displayed bound therefore makes
$\widetilde N$ of class~(D), hence a uniformly integrable $Q$-martingale.
Dropping the
two increasing terms in \eqref{eq:classD-linear} gives
$RY\le P:=|Y_0|+\widetilde N+e^{LT}c_fJ$, while
$RY\ge R\xi\ge-e^{LT}c_\xi M$ by Definition~\ref{def:rbsde}(iii).  Hence
$P+e^{LT}c_\xi M\ge0$, and
\[
\widetilde M:=P+2e^{LT}c_\xi M
\]
is a non-negative traded uniformly integrable $Q$-martingale for every $Q$,
with $\widetilde M\ge P\ge RY$ and $\widetilde M\ge e^{LT}c_\xi M\ge-RY$.  So
$|RY|\le\widetilde M$, and $|Y|\le e^{LT}\widetilde M$.
\end{proof}

\begin{thm}[Data dominated by a traded martingale]\label{thm:weighted}
Let Assumption~\ref{assum:weighted} hold for a weight $M$.  Then
\eqref{eq:rbsde} has a unique solution of class~(D), $(Y,H,A,C)$: its first
component is minimal among those of all supersolutions of class~(D), and in
particular among all supersolutions whose first component lies in
$\Sbb^{\rm tr}$.  Moreover:
\begin{enumerate}\itemsep2pt
\item[(i)] $Y$ is the unique fixed point of $\Phi$ in $\Sbb^{M'}$ for every
weight $M'\ge M$, and satisfies \eqref{eq:representation} at every stopping
time;
\item[(ii)] $|Y_t|\le c_YM_t$ for all $t$, a.s., with
$c_Y=(c_\xi+Tc_f)e^{LT}$ as in \eqref{eq:apriori};
\item[(iii)] $A$ is predictable, $C_t=-\sum_{s\le t,\,s<T}\Dp Y_s$, and $Y$, $A$,
$C$ and $\int H\,dX$ are unique; for every $Q\in\cM(X)$, $\int H\,dX$ is a
uniformly integrable $Q$-martingale and $Y+\Gamma^{Y}$ is of class~(D); and
\begin{equation}\label{eq:weighted-AC}
\sup_{Q\in\cM(X)}\E^{Q}\bigl[A_T+C_T\bigr]\ \le\ 2c_VM_0,
\qquad c_V=c_Y+T(c_f+Lc_Y).
\end{equation}
\end{enumerate}
When the data are bounded, $M\equiv1$ is a weight and $Y$ is the solution of
Theorem~\ref{thm:existence}; that solution is therefore minimal among all
supersolutions of class~(D), not only among bounded ones.
\end{thm}

\begin{proof}
Fix a weight $M'\ge M$.  The data satisfy Assumption~\ref{assum:weighted} for
$M'$ with the same constants; let $J'$ be the process of
Lemma~\ref{lem:weight-integral} built on $M'$.

\emph{$\Phi$ maps $\Sbb^{M'}$ into itself.}  Let $U\in\Sbb^{M'}$ and
$b:=\|U\|_{M'}$.  Since $|f(s,U_s)|\le(c_f+Lb)M'_s$, we have
$|\xi+\Gamma^{U}|\le c_\xi M'+(c_f+Lb)J'$, a non-negative traded uniformly
integrable martingale under every $Q$; so $\xi+\Gamma^{U}$ satisfies
Assumption~\ref{assum:G}, Theorem~\ref{thm:snell-aggreg} defines
$\Phi(U)=V(\xi+\Gamma^{U})-\Gamma^{U}$ as an optional \ladlag{} process, and
\eqref{eq:Phi-explicit} holds at every stopping time, $\Gamma^{U}_\theta$ being
$\cF_\theta$-measurable and $Q$-integrable.  For $\theta\in\cT$,
$\tau\in\cT_\theta$ and $Q\in\cM(X)$, optional sampling for $M'$ and
\eqref{eq:weight-integral} give
\[
\E^{Q}\Bigl[|\xi_\tau|+\int_\theta^{\tau}|f(s,U_s)|\,ds\Bigm|\cF_\theta\Bigr]
\le c_\xi M'_\theta+(c_f+Lb)(T-\theta)M'_\theta ,
\]
so $|\Phi(U)_\theta|\le(c_\xi+T(c_f+Lb))M'_\theta$ and $\Phi(U)\in\Sbb^{M'}$.
The domination by $J'$ alone would not give this, $J'/M'$ being unbounded in
general.

\emph{Contraction.}  Let $U^{1},U^{2}\in\Sbb^{M'}$, $\beta>L$ and
$b:=\|U^{1}-U^{2}\|_{\beta,M'}$, so that
$|U^{1}_s-U^{2}_s|\le b\,e^{\beta(T-s)}M'_s$ for all $s$ outside an evanescent
set.  The first display in the proof of Lemma~\ref{lem:Phi-properties}(ii),
with \eqref{eq:weight-integral} for $\varphi(s)=e^{\beta(T-s)}$ in place of the
deterministic bound, gives
\[
\bigl|\Phi(U^{1})_\theta-\Phi(U^{2})_\theta\bigr|
\le Lb\,M'_\theta\int_\theta^Te^{\beta(T-s)}ds
=\frac{Lb}{\beta}\,M'_\theta\bigl(e^{\beta(T-\theta)}-1\bigr),
\]
whence $\|\Phi(U^{1})-\Phi(U^{2})\|_{\beta,M'}\le(L/\beta)\|U^{1}-U^{2}\|_{\beta,M'}$.
Lemma~\ref{lem:Phi-properties}(i) and Lemma~\ref{lem:domination} hold in
$\Sbb^{M'}$ with the same proofs.  So $\Phi$ has a unique fixed point in
$\Sbb^{M'}$, and since $\Sbb^{M}\subset\Sbb^{M'}$ it is the fixed point $Y$ in
$\Sbb^{M}$.  This is (i), the representation being the identity $Y=\Phi(Y)$ at
stopping times.

\emph{Supersolution.}  Step~2 of the proof of Theorem~\ref{thm:existence} goes
through with Corollary~\ref{cor:snell-decomp}, applied under
Assumption~\ref{assum:G} to $\zeta:=\xi+\Gamma^{Y}$, in place of
Theorem~\ref{thm:zero-driver}: it gives $H$, a predictable $A$ (by
Example~\ref{eg:H1-brownian}) and the right-jump process $C$ with
$V(\zeta)=V_0+\int H\,dX-A-C_-$ and $\E^{Q}[A_T+C_T]<\infty$ for every $Q$, and
the rest of that step is unchanged.  As $Y\in\Sbb^{M}$, it is a supersolution
of class~(D).

\emph{Minimality, uniqueness, general driver.}  Let $(Y',H',A',C')$ be a
supersolution of class~(D).  By Lemma~\ref{lem:classD-tr},
$Y'\in\Sbb^{M''}$ for some weight $M''$; take $M':=M+M''$, so that
$Y'\in\Sbb^{M'}$ and the data are weighted by $M'$.  Under
\eqref{eq:monotone}, Lemma~\ref{lem:weighted-super} and
Theorem~\ref{thm:snell-minimal}, which applies because $\xi+\Gamma^{Y'}$
satisfies Assumption~\ref{assum:G} as above, give $Y'\ge\Phi(Y')$, and
Lemma~\ref{lem:domination} in $\Sbb^{M'}$ gives $Y'\ge Y$.  Steps~4 and~5 of the
proof of Theorem~\ref{thm:existence} then apply verbatim, with
Lemma~\ref{lem:weighted-super} in place of Lemma~\ref{lem:super-is-supermart};
the transformation of Lemma~\ref{lem:monotone-driver} preserves
Assumption~\ref{assum:weighted}, with the same constants and Lipschitz constant
$2L$, and maps each $\Sbb^{M'}$ onto itself and supersolutions of class~(D)
onto supersolutions of class~(D).

\emph{Bounds.}  For (ii), run the iteration of the proof of
Theorem~\ref{thm:apriori}: if $|U^{n}_s|\le h_n(s)M_s$, then optional sampling
and \eqref{eq:weight-integral} with $\varphi=c_f+Lh_n$ give
$|U^{n+1}_\theta|\le h_{n+1}(\theta)M_\theta$ at every $\theta\in\cT$, and the
rest of that argument is unchanged.  In (iii), uniqueness is Step~4, and the
martingale and class~(D) assertions are Lemma~\ref{lem:weighted-super}.  For
\eqref{eq:weighted-AC}, taking $Q$-expectations in
$A_T+C_T=V_0-V_T+\int_0^TH\,dX$, with $V=Y+\Gamma^{Y}$ and
$V_T=\xi_T+\Gamma^{Y}_T$,
\[
\E^{Q}[A_T+C_T]=Y_0-\E^{Q}\bigl[\xi_T+\Gamma^{Y}_T\bigr]
\le c_YM_0+c_\xi M_0+(c_f+Lc_Y)TM_0\le2c_VM_0 ,
\]
using $\E^{Q}[M_T]=M_0$, $\E^{Q}[J_T]=TM_0$ and $c_\xi\le c_Y$.  Finally, for
the constant weight $M\equiv1$, the weighted space $\Sbb^{M}$ and its norm
$\|\cdot\|_M$ coincide with $\Sbb^{\infty}$ and $\|\cdot\|_\infty$.
For bounded data, uniqueness of the fixed point therefore identifies $Y$
with the solution of Theorem~\ref{thm:existence}.
\end{proof}

\begin{cor}[The integrability hypotheses of the linear theory]
\label{cor:weighted-EG}
Let $\xi$ be optional and satisfy Assumption~\ref{assum:E} or
Assumption~\ref{assum:G}, and let $f$ satisfy {\bf (A2)}.  Then the
conclusions of Theorem~\ref{thm:weighted} hold.
\end{cor}

\begin{proof}
Let $M^{0}$ be the non-negative traded martingale dominating $|\xi|$, given by
Assumption~\ref{assum:G} or constructed from Assumption~\ref{assum:E} by
Corollary~\ref{cor:ui-traded-domination}.  Then $M:=1+M^{0}$ is a weight,
$|\xi|\le M$, and $|f(t,0)|\le c_f\le c_fM_t$.
\end{proof}

For drivers with $f(\cdot,0)$ bounded, the non-linear equation is therefore
well posed under exactly the integrability hypotheses under which
Section~\ref{sec:snell} constructs the linear envelope.  The remaining results
of this section carry over as well.

\begin{cor}[Comparison and the dynamic characterisation]
\label{cor:weighted-comparison}
Under Assumption~\ref{assum:weighted}, Theorem~\ref{thm:comparison},
Corollary~\ref{cor:dyn-minimality}, Lemma~\ref{lem:Y-max} and
Proposition~\ref{pro:C-skorokhod} hold with ``supersolution'' and ``solution''
read as ``supersolution of class~(D)'' and ``solution of class~(D)''.  In
Theorem~\ref{thm:comparison} the two data may be weighted by different weights.
\end{cor}

\begin{proof}
\emph{Comparison.}  If the data are weighted by $M^{1}$ and $M^{2}$, both are
weighted by the weight $M^{1}+M^{2}$ with unchanged constants, and the proof of
Theorem~\ref{thm:comparison} runs in $\Sbb^{M^{1}+M^{2}}$, where $Y^{1}$ and
$Y^{2}$ lie and Lemma~\ref{lem:domination} holds.

\emph{Dynamic characterisation.}  The identity \eqref{eq:min-identity} uses
\eqref{eq:rbsde-diff}, the martingale property of $\int H\,dX$, now
Lemma~\ref{lem:weighted-super}(i), and the integrability of the terms
conditioned, which follows from Definition~\ref{def:rbsde}(i), class~(D) of $Y$,
$|\xi|\le c_\xi M$ and $\overline\Gamma\in L^{1}(Q)$
(Lemma~\ref{lem:weighted-super}).  The concluding step uses uniqueness of the fixed point
of $\Phi$, now Theorem~\ref{thm:weighted}(i), applied in a space $\Sbb^{M'}$
with $M'\ge M$ containing $Y$, which exists by Lemma~\ref{lem:classD-tr}.

\emph{Right jumps.}  The proof of Lemma~\ref{lem:Y-max} uses boundedness only
to apply conditional dominated convergence.  Here the $f$-terms are dominated
by $(c_f+Lc_Y)J_T\in L^{1}(Q)$, while $|\xi_{\sigma_n}|\le c_\xi M_{\sigma_n}$
and $|Y_{\tau_n}|\le c_YM_{\tau_n}$ are uniformly $Q$-integrable families;
Vitali's theorem, along a common subsequence, gives the same limits.
Proposition~\ref{pro:C-skorokhod} follows from Lemma~\ref{lem:Y-max} as before.
\end{proof}

\begin{eg}[No unweighted square-integrable estimate]\label{eg:weighted-not-L2}
The market of Example~\ref{eg:G-not-UI} is of the form \eqref{eq:market} up to
the additive constant $X_0=1$, which changes neither $\cM(X)$ nor any
stochastic integral, with $n=2$, $d=1$, $T=2$ and
$\sigma_t=(\nu X_t\ind_{\{t>1\}},0)$, which is adapted, left-continuous and has
right limits, hence is predictable and locally bounded.  There $M:=1+X$ is a
weight.  In (a) and (b) take $f\equiv0$; then $\Phi(U)=V(\xi)$ for every $U$.

(a) For $\xi:=M$, $V(\xi)=M$, since $\E^{Q}[M_\tau\mid\cF_\theta]=M_\theta$ for
every $Q$ and $\tau\ge\theta$.  So the solution of Theorem~\ref{thm:weighted} is
$Y=M$, with $\int H\,dX=X-X_0$, $Z=\sigma\trans H=(\nu X_t\ind_{\{t>1\}},0)$ and
$A=C=0$.  Conditionally on $\nu$, $X_t$ is lognormal for $t>1$ with
$\E[X_t^{2}\mid\nu]=e^{\nu^{2}(t-1)}$.  Hence
$\E^{\Prob}[Y_2^{2}]\ge\E^{\Prob}[e^{\nu^{2}}]=\E[\exp(e^{2W^{2}_1})]=\infty$,
and
\[
\E^{\Prob}\int_0^2|Z_t|^{2}\,dt=\int_0^1\E^{\Prob}\bigl[\nu^{2}e^{\nu^{2}s}\bigr]\,ds
=\infty .
\]

(b) For $\xi':=M\ind_{[0,T)}$, the same computation gives
$V(\xi')_\theta=M_\theta$ on $\{\theta<T\}$ and $V(\xi')_T=0$.  So
$Y'=M\ind_{[0,T)}$, with $\int H'\,dX=X-X_0$, $C'=0$ and
$A'=M_T\ind_{\{T\}}$, which is predictable because $\cF_{T-}=\cF_T$ in the
Brownian filtration.  Thus $\E^{Q}[A'_T]=M_0$ for every $Q$, in accordance with
\eqref{eq:weighted-AC}, while $\E^{\Prob}[(A'_T)^{2}]=\E^{\Prob}[M_2^{2}]=\infty$.

(c) The running cost alone has the same effect.  Take $\xi'':=0$ and
$f(t,y):=M_t$, which satisfy Assumption~\ref{assum:weighted} with $c_\xi=0$,
$c_f=1$ and $L=0$.  Since $M>0$, the time $\tau=T$ is optimal in
\eqref{eq:Phi-explicit}, and \eqref{eq:weight-integral} gives
$\Phi(U)_\theta=(T-\theta)M_\theta$ for every $U$ and every $Q$.  So
$Y''=(T-\cdot)M$, with $Y''+\Gamma^{Y''}=J$, $\int H''\,dX=J-J_0$ and
$A''=C''=0$.  The obstacle vanishes, yet
$\E^{\Prob}[(Y''_t)^{2}]\ge(T-t)^{2}\,\E^{\Prob}\bigl[e^{\nu^{2}(t-1)}\bigr]=\infty$
for $1<t<2$.

So the bounds of Theorem~\ref{thm:apriori} on $Z$ in $\Hbb^{2}$ and on
$A_T+C_T$ in $L^{2}(Q)$ have no counterpart for data dominated by a traded
martingale, and $Y$ need not lie in $\Sbb^{2}$; by (c), this persists with a
zero obstacle as soon as $f(\cdot,0)$ is unbounded.  The weighted bound (ii) of
Theorem~\ref{thm:weighted} and the first-moment bound \eqref{eq:weighted-AC}
are what survive.
\end{eg}

\begin{rem}[Scope of the weighted theory]\label{rem:weighted-scope}
Minimality and uniqueness hold in the natural class, the supersolutions of
class~(D) under every $Q$, and Lemma~\ref{lem:classD-tr} is what makes the
weighted norm no restriction there.  Two limits remain.  The regulators obey
the first-moment bound \eqref{eq:weighted-AC} and, by
Example~\ref{eg:weighted-not-L2}, no square-integrable one.  And the driver
still does not depend on the control (Problem~P2).
\end{rem}

\section{Robust superhedging of American options}\label{sec:financial}

In this section we return to the general market of Section~\ref{sec:prelim}:
$X$ is a locally bounded \cadlag{} semimartingale with $\Prob\in\cM(X)$.  The
discounted payoff of the American claim is an optional process $\xi$ satisfying
Assumption~\ref{assum:E} or Assumption~\ref{assum:G}; no path regularity is
required for the duality, and no two-sided bound.  The holder may exercise at
any stopping time $\tau$ and then receives $\xi_\tau$.

\subsection{Superhedging}

\begin{defn}\label{def:admissible}
A pair $(x,H)$ with $x\in\R$ and $H$ predictable and $X$-integrable is an
\emph{admissible strategy} if the gross wealth process
\[
\cS^{x,H}_t:=x+\int_0^tH_s\,dX_s,\qquad 0\le t\le T,
\]
is bounded below by a constant, i.e.\ $\cS^{x,H}\ge-K$ for some $K\ge0$.  It
\emph{superhedges} the American claim $\xi$ if
\[
\cS^{x,H}_t\ \ge\ \xi_t\qquad\text{for all }t\in[0,T],\ \text{a.s.}
\]
The \emph{superhedging price} of $\xi$ is
\[
p_0(\xi):=\inf\bigl\{x\in\R:\ \exists\,H\ \text{such that }(x,H)\ \text{is
admissible and superhedges }\xi\bigr\}.
\]
\end{defn}

Two remarks on the lower bound.  It is not what keeps the infimum finite:
taking $t=0$ in the superhedging requirement already gives $x\ge\xi_0$, so
$p_0(\xi)>-\infty$ for every optional $\xi$.  Nor is it a restriction when
$\xi$ is bounded below: a wealth process dominating $\xi$ at \emph{all} times
then satisfies $\cS^{x,H}\ge\xi\ge-c$ of its own accord, so every superhedging
pair is automatically admissible, and the clause is redundant wherever
Theorem~\ref{thm:robust-american} applies.  What a constant lower bound does
is make the dual inequality available: it is through the Ansel--Stricker lemma
\cite{AnselStricker1994} that $\cS^{x,H}$ is a $Q$-supermartingale, and hence
that $x\ge\E^{Q}[\xi_\tau]$ for every $Q\in\cM(X)$ and every $\tau\in\cT$.
For a payoff that is not bounded below the domination $\cS^{x,H}\ge\xi$
supplies no such bound, and restoring one is precisely the role of the credit line in
Theorem~\ref{thm:robust-american-signed}.  Doubling strategies, by contrast,
are the obstruction in the terminal-only problem, where nothing at all is
required of the wealth before $T$.

\begin{thm}\label{thm:robust-american}
Let $\xi$ be an optional process which is bounded below, $\xi\ge-c$, and
satisfies Assumption~\ref{assum:E} or Assumption~\ref{assum:G}.  Let
$V=V(\xi)$ be its generalised
Snell envelope (Theorem~\ref{thm:snell-aggreg}), and let
$V=V_0+\int H\,dX-A-C_-$ be its decomposition (Corollary~\ref{cor:snell-decomp}).
Then the infimum in Definition~\ref{def:admissible} is attained and
\begin{equation}\label{eq:duality}
p_0(\xi)=V_0(\xi)
=\sup_{Q\in\cM(X)}\ \sup_{\tau\in\cT}\ \E^{Q}\bigl[\xi_\tau\bigr].
\end{equation}
An optimal superhedge is given by $(V_0,H)$, with gross wealth
\begin{equation}\label{eq:wealth}
\cS^{V_0,H}_t=V_t+A_t+C_{t-},\qquad 0\le t\le T .
\end{equation}
In the It\^o market of Section~\ref{subsec:market} with $\xi$ satisfying
{\bf (A1)}, $V=Y$ is the first component of the solution of \eqref{eq:rbsde}
with $f\equiv0$ (Theorem~\ref{thm:zero-driver}), $A$ is predictable and the
gross wealth is continuous.
\end{thm}

\begin{proof}
\emph{$V_0$ is attained as a superhedging capital.}  By
Corollary~\ref{cor:snell-decomp}, $\cS^{V_0,H}=V+A+C_-$, which is
\eqref{eq:wealth}.  Since $A,C\ge0$ and $V\ge\xi$ we get
$\cS^{V_0,H}\ge V\ge\xi\ge-c$, so $(V_0,H)$ is admissible and
superhedges; hence $p_0(\xi)\le V_0$.  Only the lower bound on $\xi$ is used
here, not a two-sided one: this is why Definition~\ref{def:admissible} needs no
change for payoffs bounded below, however large they are above.

\emph{No smaller capital suffices.}  Let $(x,H')$ be admissible and superhedge
$\xi$, and put $U:=\cS^{x,H'}$.  For every $Q\in\cM(X)$, $U$ is a stochastic
integral with respect to the locally bounded $Q$-local martingale $X$ and is
bounded below, hence a $Q$-local martingale by the Ansel--Stricker lemma
\cite{AnselStricker1994} and therefore a $Q$-supermartingale.  Being \cadlag{}
and bounded below, it satisfies the optional sampling inequality
$U_\sigma\ge\E^{Q}[U_\tau\mid\cF_\sigma]$ for all $\sigma\le\tau$ in $\cT$, so
$U$ is a strong optional $Q$-supermartingale, for every $Q\in\cM(X)$.  Since
$U\ge\xi$, Theorem~\ref{thm:snell-minimal} gives $U\ge V$, and in particular
$x=U_0\ge V_0$.  Hence the infimum equals $V_0$ and is attained.

Finally $V_0=\sup_{Q}\sup_{\tau}\E^{Q}[\xi_\tau]$ by the definition of $V_0$,
$\cF_0$ being trivial.  The statement about the
It\^o market is Theorem~\ref{thm:zero-driver} together with the continuity of
$\int H\,dX$ there.
\end{proof}

\begin{rem}[Gross wealth and consumption]\label{rem:consumption}
The theorem is a statement about the gross self-financing wealth
$V_0+\int H\,dX$.  If the seller is allowed to withdraw funds, one may read
\[
V_t=V_0+\int_0^tH_s\,dX_s-A_t-C_{t-}
\]
as net wealth after cumulative withdrawals $A+C_-$, with the convention that a
jump $\Dm C_t$ is withdrawn immediately after the exercise opportunity at $t$
and therefore enters net wealth only at times strictly later than $t$; $A+C_-$
is then an increasing process that need not be right-continuous.  This reading
requires the convention to be stated, and the duality does not depend on it;
nothing forces the seller to withdraw, and a gross hedge may retain its
surplus.
\end{rem}

\begin{rem}[Relation to Kramkov's American theorem]\label{rem:kramkov-american}
Kramkov \cite[Theorem~3.3]{Kramkov96} states, in the same generality of market,
that for a positive adapted reward $f$ with
$\sup_{\tau}\sup_{Q}\E^{Q}f_\tau<\infty$ there is a minimal hedging portfolio,
whose capital at time $t$ is $\esssup_{Q,\,\tau\ge t}\E^{Q}[f_\tau\mid\cF_t]$.
His proof passes through \cite[Proposition~4.3]{Kramkov96}, stated for a
\cadlag{} $f$, where right-continuity of $f$ yields a \cadlag{} modification of
the value process; capital processes, with right-continuous consumption, are
\cadlag.  Karatzas and Kou \cite{KaratzasKou1998} treat the constrained problem
for continuous payoffs.

The right-continuity cannot be dropped from that conclusion.  Take $\F$
trivial, $X\equiv0$ and $f_t:=\ind_{[0,1/2]}(t)$.  The capitals
$\ind_{[0,1/2+\epsilon)}$ hedge $f$ for every $\epsilon>0$, so a capital
dominated by that of every hedge would equal $\ind_{[0,1/2]}$, which is not
right-continuous and is the capital of no portfolio; the minimal initial
capital, $1$, is attained, but there is no minimal hedging portfolio.  For a
bounded optional payoff that is not right-continuous the value $V(\xi)$ is in
general only \ladlag{}, and the minimality that Kramkov states for the capital
is relocated to the Snell envelope (Theorem~\ref{thm:snell-minimal}); what is a
capital process is $V+A+C_-$.  For right-continuous $\xi$ one has $C\equiv0$
and \eqref{eq:duality} is \cite[Theorem~3.3]{Kramkov96}.
\end{rem}

\begin{rem}[Reading the two increasing processes]\label{rem:reading-AC}
In the It\^o market, $A$ is predictable and $C_t=-\sum_{s\le t,\,s<T}\Dp Y_s$
records the downward right jumps of the value, which occur only at contact
points (Proposition~\ref{pro:C-skorokhod}), that is, where the payoff fails to
be right-continuous.  It is the presence of $C$ that allows an irregular
payoff to be superhedged by a \cadlag{} wealth process while the value is only
\ladlag.  If $\xi$ is right-continuous then $C\equiv0$: at every stopping time
$Y_{\tau+}\ge\xi_{\tau+}=\xi_\tau$, so \eqref{eq:Y-max} gives
$Y_\tau=Y_{\tau+}$.  Then \eqref{eq:duality} reduces to
\cite[Theorem~3.3]{Kramkov96}, and to the superhedging duality of El Karoui and
Quenez \cite{ElKarouiQuenez1995} when in addition $X$ is continuous.  If,
further, the market is complete, $A$ is the regulator of the classical Snell
envelope and $p_0(\xi)=\sup_\tau\E^{\Prob}[\xi_\tau]$.
\end{rem}

For payoffs that are not bounded below, the conclusion of
Theorem~\ref{thm:robust-american} can fail in the class of
Definition~\ref{def:admissible}, and the obstruction is not an artefact of
incompleteness.

\begin{eg}[Attainment can fail for payoffs unbounded below]\label{eg:no-attainment}
Let $X=W$ be a one-dimensional Brownian motion in its augmented filtration, so
that $\cM(X)=\{\Prob\}$ and the market is complete, and let $\xi_t:=W_t$, which
has finite moments of every order at every stopping time and so satisfies the
hypothesis of Lemma~\ref{lem:Lp}.  Here $V_t=W_t$ and $V_0=0$.  Suppose
$(0,H')$ were admissible in the sense of Definition~\ref{def:admissible} and
superhedged $\xi$, and put $U:=\cS^{0,H'}$.  Being a stochastic integral
bounded below, $U$ is a $\Prob$-supermartingale, so $U-W$ is a non-negative
supermartingale with $(U-W)_0=0$, hence $U=W$ --- which has no deterministic
lower bound.  So no admissible strategy attains $V_0=0$.

The value is nevertheless still $0$: the wealth processes
$U^{k}_t:=\E[\max(W_T,-k)\mid\cF_t]$ are stochastic integrals of $W$, satisfy
$U^{k}\ge-k$ and $U^{k}_t\ge W_t$ by conditional Jensen, and
$U^{k}_0=\E[\max(W_T,-k)]\downarrow0$ as $k\to\infty$.  Thus
$p_0(\xi)=0=V_0$ while the infimum is not attained.  Equality of values and
attainment in a given class are different assertions, and weakening the payoff
assumption without touching Definition~\ref{def:admissible} can destroy the
second.
\end{eg}

\begin{thm}[Signed payoffs, with a common credit line]
\label{thm:robust-american-signed}
Let $\xi$ satisfy Assumption~\ref{assum:E} or Assumption~\ref{assum:G}, and
let $M$ be the common non-negative traded uniformly integrable martingale
dominating $|\xi|$ --- given by Assumption~\ref{assum:G}, or constructed from
Assumption~\ref{assum:E} by Corollary~\ref{cor:ui-traded-domination}.  Call
$(x,H)$ \emph{$M$-admissible} if $\cS^{x,H}\ge-M$.  Then
\[
\inf\bigl\{x:\ \exists H,\ (x,H)\ M\text{-admissible and superhedging }\xi
\bigr\}
=V_0(\xi)=\sup_{Q\in\cM(X)}\sup_{\tau\in\cT}\E^{Q}[\xi_\tau],
\]
and the infimum is attained, by $(V_0,H)$ with $H$ from
Corollary~\ref{cor:snell-decomp}.
\end{thm}

\begin{proof}
\emph{Attainment.}  As in Theorem~\ref{thm:robust-american},
$\cS^{V_0,H}=V+A+C_-\ge V\ge\xi\ge-M$, so $(V_0,H)$ is $M$-admissible and
superhedges.

\emph{Lower bound.}  Let $(x,H')$ be $M$-admissible and superhedge $\xi$, and
put $U:=\cS^{x,H'}$.  Fix $Q\in\cM(X)$.  Then
$U+M-(x+m)=\int_0^{\cdot}(H'+G)\,dX\ge-(x+m)$ is a stochastic integral of $X$
bounded below by a constant, hence a $Q$-local martingale by the
Ansel--Stricker lemma \cite{AnselStricker1994}; so $U+M$ is a non-negative
$Q$-local martingale and therefore a $Q$-supermartingale.  Since $M$ is a uniformly integrable $Q$-martingale,
$U=(U+M)-M$ is a $Q$-supermartingale as well, and optional sampling gives
$x=U_0\ge\E^{Q}[U_\tau]\ge\E^{Q}[\xi_\tau]$ for every $\tau\in\cT$.  Taking
suprema over $\tau$ and $Q$ gives $x\ge V_0$ by
Theorem~\ref{thm:snell-aggreg}.
\end{proof}

\begin{rem}[Which admissibility class]\label{rem:admissibility-choice}
Definition~\ref{def:admissible} is unchanged, and
Theorem~\ref{thm:robust-american} continues to use it: for payoffs bounded
below --- including non-negative payoffs unbounded above --- nothing needs to be
said.  Theorem~\ref{thm:robust-american-signed} is a separate statement with a
separate class, and Example~\ref{eg:no-attainment} shows that the change of
class is necessary rather than cosmetic.  An alternative formulation, with the
same proof, takes as admissible those strategies whose gross wealth is a
$Q$-supermartingale for every $Q\in\cM(X)$; the optimal wealth
\eqref{eq:wealth} lies in that class because, by
Corollary~\ref{cor:snell-decomp}, it \emph{is} the uniformly integrable
$Q$-martingale $V_0+\int H\,dX$, for every $Q\in\cM(X)$.  That $V$ is a
uniform supermartingale and that $A$ and $C$ are increasing would not give
this: adding an increasing process to a supermartingale need not leave a
supermartingale.
\end{rem}

\begin{rem}[Two components of $A$]\label{rem:two-sources}
By Proposition~\ref{pro:no-skorokhod}, $A$ does not in general increase only on
the contact set: part of it reflects the value on the payoff, and part of it is
the premium for non-attainability, already present for European claims.  We do
not identify the second part with the regulator of the European claim $\xi_T$;
the two are generated by different value functions and there is no reason for
them to agree.  In the Markovian specification of Example~\ref{eg:stochvol}
studied in \cite{BFOAmerican}, the American value $w$ is strictly convex in the
continuation region, the maximising volatility there is $\varsigma_{\max}$,
and the consumption rate of the American superhedge in the continuation region
is \emph{formally} $\tfrac12[\varsigma_{\max}^{2}-\varsigma(v_t)^{2}]\,w_{xx}$;
we give no verification argument for this Markovian formula.  Problem~P3 asks
for the split of $A$ in general.
\end{rem}

\subsection{Examples}

\begin{eg}[American put with an irregular exercise rule]\label{eg:windowed}
In the It\^o market let $P=(P_t)$ be a positive asset price, which is
continuous, and consider the put payoff $\xi^{0}_t=(K-P_t)^{+}$, so that
$0\le\xi^{0}\le K$ and {\bf (A1)} holds with $c_\xi=K$.  Payoffs unbounded
above, such as the American call, are not covered by {\bf (A1)}.  They
\emph{are} covered by the linear results of
Sections~\ref{sec:snell}--\ref{sec:financial}: a
non-negative payoff dominated by a traded non-negative asset that is a true
uniformly integrable $Q$-martingale for every $Q\in\cM(X)$ satisfies
Assumption~\ref{assum:G}, and being bounded below it is covered by
Theorem~\ref{thm:robust-american} in the unchanged admissible class.  The
qualification matters: an asset which is only a local martingale under some
pricing measure does not satisfy Assumption~\ref{assum:G}.  For the non-linear
equation of Section~\ref{sec:bsde} such payoffs are covered, for drivers with
$f(\cdot,0)$ bounded, by Corollary~\ref{cor:weighted-EG}.

Suppose the contract permits exercise only on a set $D\subset[0,T]$ which is a
finite union of closed intervals and isolated dates, with $T\in D$.  Put
\[
\xi_t:=\xi^{0}_t\ind_{D}(t)-K\ind_{D^{c}}(t).
\]
Then $\xi$ is optional and \ladlag, and the penalty $-K$ lies strictly below
$\xi_T\ge0$: a stopping time taking values outside $D$ may be replaced by $T$
on that event, which does not decrease the payoff, so the superhedging price of
the restricted contract is $\sup_Q\sup_\tau\E^{Q}[\xi_\tau]$ and
Theorem~\ref{thm:robust-american} applies.  The value process can have right
jumps only at right endpoints of windows and at isolated dates: elsewhere
either $\xi$ is right-continuous at the point or $Y>\xi$ there (outside $D$,
$Y_t\ge\E^{\Prob}[\xi_T\mid\cF_t]\ge0>-K=\xi_t$), and in both cases
\eqref{eq:Y-max} gives $Y_\tau=Y_{\tau+}$.  A right jump at such a point is
possible but not automatic.

More generally, let $D\subset[0,T]$ be any deterministic Borel set containing
$T$.  The same formula defines a bounded optional reward, so both
Theorem~\ref{thm:robust-american} and the zero-driver reflected-BSDE
representation of Theorem~\ref{thm:zero-driver} apply.  The replacement of an
exercise time outside $D$ by $T$ is still a stopping time and still improves
the payoff.  Taking $T=1$ and $D=\{0\}\cup\{1/n:n\ge1\}$ gives an effective
reward with no right limit at $0$, so this extension includes contracts
excluded by a \ladlag{} assumption on the obstacle.  The value nevertheless
has left and right limits, and its right-jump regulator satisfies
Proposition~\ref{pro:C-skorokhod}.
\end{eg}

\begin{eg}[Stochastic volatility]\label{eg:stochvol}
Let $n=2$, $d=1$ and $dX_t=\varsigma(v_t)\,dW^{1}_t$, where $v$ solves an
autonomous SDE driven by $W^{2}$ and $\varsigma$ is bounded and strictly
positive.  Then $\sigma_t=(\varsigma(v_t),0)$, $\cK_t=\R\times\{0\}$, and the
elements of $\cM(X)$ are obtained by changing the drift of $W^{2}$ alone; the
market is incomplete and $\cM(X)$ is infinite.  For a European payoff that is
not attainable by bounded wealth, $A\neq0$ by
Proposition~\ref{pro:no-skorokhod}.

In the specification of \cite{BFOAmerican,BFOEuropean} --- an
Ornstein--Uhlenbeck factor and
the smooth increasing volatility function $\varsigma$ given there, with range
$(\varsigma_{\min},\varsigma_{\max})$, $\varsigma_{\min}>0$ --- the
superhedging price is identified explicitly.  For a bounded Lipschitz European payoff $g(X_T)$ it is
the value of the Black--Scholes--Barenblatt equation with volatility bounds the
range of $\varsigma$, a known identification
\cite{AvellanedaLevyParas1995,Lyons1995,CvitanicPhamTouzi1999}.  For the
American claim $g(X_t)$ with zero interest rate and time-independent payoff it
is the value of the Barenblatt obstacle problem, which coincides with the
American Bachelier price at the maximal volatility $\varsigma_{\max}$
\cite{BFOAmerican}; the same three-way identification holds at a positive interest
rate when $g\ge0$, with both values discounted, on the funding convention used
there, in which $X$ is a zero-carry quotation --- a futures or forward price,
not the discounted price of \eqref{eq:market} --- and the cash balance earns
the rate $r$.  For
European claims, \cite{BFOEuropean} quantifies the approximation of this value
by the measures of $\cM(X)$ whose drift change on $W^{2}$ is bounded by
$\kappa$: for bounded Lipschitz payoffs the error is
$O(\kappa^{-1/2}(1+\log\kappa)^{3/2})$, uniformly over payoffs with given
bounds on $\|g\|_\infty$ and $\mathrm{Lip}(g)$, and the worst-case error over
such a class is at least of order $\kappa^{-1/2}$; for an explicit payoff with
a single sign change of gamma the error is of exact order $\kappa^{-1}$.  Nothing
in Sections~\ref{sec:decomposition}--\ref{sec:financial} depends on these
results.
\end{eg}

\subsection{Non-linear wealth dynamics}

A driver $f$ satisfying {\bf (A2)} can model a wealth equation in which
borrowing and lending rates differ, or in which the cash account earns a
state-dependent rate; $f(t,y)=-r_ty$ with $r$ bounded and predictable is
ordinary discounting.  We call $(x,H)$ an admissible non-linear strategy if
the solution of
\begin{equation}\label{eq:nonlinear-S}
\cS^{x,H}_t=x-\int_0^tf(s,\cS^{x,H}_s)\,ds+\int_0^tH_s\,dX_s,
\qquad 0\le t\le T,
\end{equation}
is bounded below by a constant, and say that it superhedges $\xi$ if
$\cS^{x,H}\ge\xi$.  For a given continuous gain process $\int H\,dX$,
\eqref{eq:nonlinear-S} has a unique continuous adapted solution, pathwise, by
the Lipschitz property of $f$.

\begin{thm}\label{thm:nonlinear-wealth}
In the It\^o market, let $(f,\xi)$ satisfy Assumption~\ref{assum:data} and let
$(Y,H,A,C)$ be the solution of \eqref{eq:rbsde}.  The infimum of the initial
capitals of admissible non-linear strategies superhedging $\xi$ equals $Y_0$,
and it is attained by $(Y_0,H)$.  In general $\cS^{Y_0,H}\neq Y+A+C_-$.
\end{thm}

\begin{proof}
Write $K:=A+C_-$ and, for a continuous adapted $\cS$, put $D:=\cS-Y$,
\[
b_t:=\frac{f(t,\cS_t)-f(t,Y_t)}{\cS_t-Y_t}\ \ \text{on}\ \{\cS_t\neq Y_t\},
\qquad b_t:=0\ \ \text{otherwise},
\qquad R_t:=\exp\Bigl(\int_0^tb_s\,ds\Bigr),
\]
so that $|b|\le L$ and $e^{-LT}\le R\le e^{LT}$.

\emph{Upper bound.}  Let $\cS:=\cS^{Y_0,H}$.  By \eqref{eq:rbsde-diff},
$D_t=-\int_0^tb_sD_s\,ds+K_t$, a linear equation driven by the finite-variation
process $K$, whose solution is
\begin{equation}\label{eq:nonlinear-domination}
R_tD_t=\int_{(0,t]}R_s\,dA_s+\int_{[0,t)}R_s\,dC_s\ \ge0 .
\end{equation}
Hence $\cS\ge Y\ge\xi\ge-c_\xi$, so $(Y_0,H)$ is admissible and superhedges.

\emph{Lower bound.}  Let $\cS=\cS^{x,H'}$ be an admissible superhedge, bounded
below by $-c_{\cS}$.  Now
\[
R_tD_t=x-Y_0+N_t+J_t,\qquad
N_t:=\int_0^tR_s(H'_s-H_s)\,dX_s,\qquad
J_t:=\int_{(0,t]}R_s\,dA_s+\int_{[0,t)}R_s\,dC_s .
\]
Under each $Q\in\cM(X)$, $N$ is a local martingale, and
$N\ge-e^{LT}(c_{\cS}+c_Y)-|x-Y_0|-e^{LT}K_T$, a lower bound that is
$Q$-integrable because $\E^{Q}[A_T+C_T]<\infty$.  By Fatou's lemma $N$ is a
$Q$-supermartingale.  For $\theta\in\cT$, $D_\theta\ge\xi_\theta-Y_\theta=-G_\theta$
with $G:=Y-\xi\ge0$, so $R_\theta D_\theta\ge-e^{LT}G_\theta$ and
$J_\theta\le e^{LT}K_\theta$.  Taking expectations,
\[
x-Y_0\ \ge\ \E^{Q}[R_\theta D_\theta]-\E^{Q}[J_\theta]
\ \ge\ -e^{LT}\,\E^{Q}\bigl[K_\theta+G_\theta\bigr].
\]
By Corollary~\ref{cor:dyn-minimality} at $\tau=0$, $\cF_0$ being trivial, the
infimum over $(Q,\theta)$ of $\E^{Q}[K_\theta+G_\theta]$ is zero.  Hence
$x\ge Y_0$.  Finally \eqref{eq:nonlinear-domination}, not the unweighted sum
$K$, describes the surplus of the attaining wealth over $Y$.
\end{proof}

For $f\equiv0$, Theorem~\ref{thm:nonlinear-wealth} is the It\^o-market case of
Theorem~\ref{thm:robust-american}.  By Theorem~\ref{thm:comparison}, the
superhedging price is monotone in the payoff and in the driver.

\section{Concluding remarks and open problems}\label{sec:conclusion}

We have proved a uniform optional decomposition for strong optional
supermartingales with an explicit right-jump part
(Theorem~\ref{thm:uniform}), the superhedging duality for irregular optional
American payoffs under explicit integrability conditions, in a general locally
bounded market
(Theorem~\ref{thm:robust-american}), and, in an incomplete Brownian market,
well-posedness, comparison and a dynamic characterisation of minimal reflected
BSDEs with an optional obstacle, without path regularity assumptions
(Theorems~\ref{thm:existence}
and~\ref{thm:comparison}, Corollary~\ref{cor:dyn-minimality}), for bounded data
and for data dominated by a traded martingale (Theorem~\ref{thm:weighted},
Corollary~\ref{cor:weighted-comparison}), together with an
exact criterion for the classical Skorokhod conditions
(Proposition~\ref{pro:no-skorokhod}).  Three questions remain open within this
framework.

\begin{itemize}
\item[\textbf{P1.}]
\textbf{(Unbounded data.)}  The \emph{linear} part of this question is now
settled, and so is the non-linear part for data dominated by a traded
martingale.  Assumption~\ref{assum:E} is a condition on the data alone, it implies
class~(D) of the robust envelope (Proposition~\ref{pro:envelope-UI}), a uniform
$L^{p}$ bound at stopping times with $p>1$ implies it (Lemma~\ref{lem:Lp}), and
Assumption~\ref{assum:G} is an alternative sufficient condition
(Corollary~\ref{cor:envelope-common-domination}); together with
Theorem~\ref{thm:uniform} at class~(D) these give
Corollary~\ref{cor:snell-decomp} and Theorems~\ref{thm:robust-american}
and~\ref{thm:robust-american-signed} for unbounded payoffs.  None of this uses
a maximal inequality.

For the non-linear equation, Theorem~\ref{thm:weighted} and
Corollary~\ref{cor:weighted-comparison} give existence, uniqueness, the
representation, comparison and the dynamic characterisation for data dominated
by a traded martingale, with minimality among all supersolutions of class~(D)
under every $Q$; by Corollary~\ref{cor:weighted-EG} this covers
Assumptions~\ref{assum:E} and~\ref{assum:G} whenever $f(\cdot,0)$ is bounded.
None of this uses a maximal inequality either.

One question remains, and it concerns quantitative estimates.
Example~\ref{eg:weighted-not-L2} shows that the square-integrable bounds of
Theorem~\ref{thm:apriori} fail for data dominated by a traded martingale, and
its part~(c) shows that square integrability of the obstacle alone does not
restore them once $f(\cdot,0)$ is unbounded.  For square-integrable obstacles,
$\sup_{Q}\E^{Q}[\sup_t|\xi_t|^{2}]<\infty$, while retaining {\bf (A2)}
including boundedness of $f(\cdot,0)$, analogous estimates would follow from
a Doob-type inequality
\[
\sup_{Q\in\cM(X)}\E^{Q}\Bigl[\sup_{t\le T}
\bigl(\esssup_{Q'\in\cM(X)}\E^{Q'}[\rho\mid\cF_t]\bigr)^{2}\Bigr]
\ \le\ c\,\sup_{Q\in\cM(X)}\E^{Q}\bigl[\rho^{2}\bigr],
\qquad \rho\ge0 ,
\]
which would give the a priori estimates of Theorem~\ref{thm:apriori} for such
data.  Such inequalities are
known for sublinear expectations generated by families with uniformly bounded
density generators, or under a reverse H\"older condition on the densities.
Those techniques do not apply to the full set $\cM(X)$, and we do not know
whether the inequality holds there; the failure of a sufficient technique is
not an argument that it fails.  Note that the estimate
$\sup_{Q,\theta}\E^{Q}[R_\theta^{p}]\le K_p$ of
Proposition~\ref{pro:envelope-UI} is \emph{not} of this kind: its supremum over
stopping times is outside the expectation, so it says nothing about
$\E^{Q}[\sup_t R_t^{p}]$.
\item[\textbf{P2.}]
\textbf{(Drivers depending on the control.)}  For a driver $f(t,y,z)$ the
natural fixed-point argument requires an estimate for the difference of the
integrands of two solutions, obtained by applying the Gal'chouk--Lenglart
formula to $|\delta Y|^{2}$.  The term $\int\delta Y_s\,d(\delta C)_s$ can be
signed using Proposition~\ref{pro:C-skorokhod}; the term
$\int\delta Y_{s-}\,d(\delta A)_s$ is usually signed with the Skorokhod
condition on $A$, which is unavailable for some data here.  For
right-continuous values and $d\Prob\otimes dt$-a.e.\ domination, constrained
reflected BSDEs with $z$-dependent drivers are treated in \cite{PengXu2010} by
separating the regulator; the question is whether a difference estimate can be
obtained for the single regulator of \eqref{eq:rbsde} with domination at every
stopping time.
\item[\textbf{P3.}]
\textbf{(Separating reflection from the premium for non-attainability.)}
For this problem retain the bounded data of Assumption~\ref{assum:data} and
assume additionally that $\xi$ has finite left limits on $(0,T]$, outside one
null set.  This assumption is local to the contact-based split below; the
existence and dynamic-minimality theorems do not require it.  Write
$A=A^{c}+A^{d}$ and define
\[
A^{\mathrm{refl}}_t:=\int_0^t\ind_{\{Y_{s-}=\xi_{s-}\}}\,dA^{c}_s
+\sum_{0<s\le t}\ind_{\{Y_{s-}=\xi_{s-}\}}\Dm A_s,
\qquad
A^{\mathrm{cons}}:=A-A^{\mathrm{refl}} .
\]
Both are predictable and increasing, and $A^{\mathrm{refl}}$ satisfies the
first two contact conditions of Remark~\ref{rem:why-minimal} by construction.
The jump of $A$ at a predictable time is classified by \emph{left} contact:
the indicator $\ind_{\{Y_s=\xi_s\}}$, based on post-jump contact, would
misclassify ordinary reflection.  For example, in the complete market $X=W$
with $T=1$, $f=0$ and the deterministic obstacle equal to $1$ on $[0,\frac12)$,
to $0$ on $[\frac12,1)$ and to $\frac12$ at $1$, one has $Y=1$ on $[0,\frac12)$,
$Y=\frac12$ on $[\frac12,1]$, $A=\frac12\ind_{[1/2,1]}$ and $H=C=0$; at
$t=\frac12$ there is left contact, $Y_{1/2-}=\xi_{1/2-}=1$, but post-jump gap
$Y_{1/2}-\xi_{1/2}=\frac12$.  (For the continuous part the two indicators give
the same process: $dA^{c}$ does not charge the countable sets where $Y$ or
$\xi$ differs from its left limit.)  The questions are:
\begin{enumerate}
\item[(a)] does the dynamic minimum condition \eqref{eq:min-condition} imply
\[
\essinf_{Q\in\cM(X)}\ \E^{Q}\bigl[A^{\mathrm{cons}}_T-A^{\mathrm{cons}}_\tau
\ \big|\ \cF_\tau\bigr]=0,\qquad\tau\in\cT,
\]
and do these separated conditions, together with \eqref{eq:rbsde}, conversely
imply minimality?
\item[(b)] does the solution satisfy the Skorokhod-type condition
\[
\essinf_{Q\in\cM(X)}\ \E^{Q}\Bigl[\int_{(\tau,T]}\bigl(Y_{s-}-\xi_{s-}\bigr)\,dA_s
\ \Big|\ \cF_\tau\Bigr]=0,\qquad\tau\in\cT,
\]
the dominated analogue of the condition shown in \cite{MPZcorrigendum} to be
equivalent, under an assumption on the obstacle, to the minimality of
second-order reflected BSDEs?
\end{enumerate}
We assert neither.  The companion paper on the buyer's problem \cite{P4}
establishes, under a structural condition on the payoff, a decomposition of the
sup--inf value with two increasing processes of opposite signs, one of them
supported on the contact set.  That result concerns a different value process:
the order of the optimisation is reversed, so the stopper and the choice of
measure are opponents there and allies here, and the minimality conditions
differ accordingly.  It does not settle (a) or (b).  In the complete case
$A^{\mathrm{cons}}=0$ by
Remark~\ref{rem:complete-case}, and both conditions hold.  For the European
obstacle of Proposition~\ref{pro:no-skorokhod}, $A^{\mathrm{refl}}=0$ and the
representation gives (a); (b) follows there from (a) since $Y-\xi$ is
bounded.  These cases do not establish the general statements.  The
corresponding splitting of the increasing processes of second-order reflected
BSDEs is posed as open in \cite[Remark~3.5, p.~2433]{MatoussiPossamaiZhou2013}; for a
right-continuous obstacle, the two-process formulation of \cite{PengXu2010} is
a natural candidate for the split.  Here the question remains open for
obstacles with left limits, including the \ladlag{} case.
\end{itemize}

Beyond these, two directions seem worth pursuing.  First, a general family of
priors $\cP\subsetneq\cM(X)$, as in the Knightian models of
\cite{CE2002,Riedel2009,Nutz2013}: Section~\ref{sec:snell} carries over, since
only $m$-stability is used, but Theorem~\ref{thm:uniform} does not, because the
identification of $\cP$-uniform local martingales with stochastic integrals of
$X$ is specific to the full set $\cM(X)$ (Remark~\ref{rem:full-family}); what
survives is a decomposition with an abstract martingale part, as in
\cite{Arharas2025}.  Second, volatility uncertainty with mutually singular
priors, where the right-jump regularisation of Lemma~\ref{lem:aggreg-C} would
have to be combined with quasi-sure aggregation; the second-order reflected
BSDEs of \cite{MatoussiPossamaiZhou2013,MPZcorrigendum} and the superhedging
results of \cite{Nutz2015} are the relevant comparison points.

\paragraph{Acknowledgments.}
The third author gratefully acknowledges the support of the Hassan~II Academy of
Sciences and Technology, Rabat, Morocco.


\end{document}